\documentclass[reqno]{amsart}
\usepackage[margin=1.5in,bottom=1.25in,marginpar=.8in]{geometry}		

\usepackage{amsmath}		
\usepackage{amssymb}		
\usepackage{amsfonts}		
\usepackage{amsthm}		
\usepackage[foot]{amsaddr}		

\usepackage{mathtools}		

\mathtoolsset{%
}

\usepackage[utf8]{inputenc}		
\usepackage[T1]{fontenc}		

\usepackage[
cal=cm,
]
{mathalfa}

\usepackage{dsfont}		

\usepackage[proportional,tabular,lining,sf,mono=false]{libertine}

\usepackage[dvipsnames,svgnames]{xcolor}
\colorlet{MyBlue}{MediumBlue}
\colorlet{MyGreen}{DarkGreen!85!Black}

\usepackage[]{titlesec}		
\titleformat{name=\section}{}{\thetitle.}{0.8em}{\centering\scshape}
\titleformat{name=\subsection}[runin]{}{\thetitle.}{0.5em}{\bfseries}[.]
\titleformat{name=\subsubsection}[runin]{}{\thetitle.}{0.5em}{\bfseries}[.]

\titleformat{name=\paragraph,numberless}[runin]{}{}{0em}{\bfseries}[]
\titlespacing{\paragraph}{0em}{\medskipamount}{1em}
\titleformat{name=\subparagraph,numberless}[runin]{}{}{0em}{}[.]
\titlespacing{\subparagraph}{0em}{0em}{0.5em}

\newcommand{\afterhead}{.\;}		

\usepackage[font=small,labelfont=bf]{caption}
\usepackage{subfigure}
\usepackage{tikz}

\usepackage[kerning=true]{microtype}		

\usepackage{acronym}
\usepackage{booktabs}       
\usepackage{latexsym}
\usepackage{paralist}
\usepackage{nicefrac}       
\usepackage{paralist}
\usepackage{wasysym}
\usepackage{xspace}

\usepackage[sort&compress]{natbib}
\bibpunct[, ]{[}{]}{,}{n}{}{,}%
\usepackage{hyperref}
\hypersetup{
colorlinks=true,
linktocpage=true,
pdfstartview=FitH,
breaklinks=true,
pdfpagemode=UseNone,
pageanchor=true,
pdfpagemode=UseOutlines,
plainpages=false,
bookmarksnumbered,
bookmarksopen=false,
bookmarksopenlevel=1,
hypertexnames=true,
pdfhighlight=/O,
urlcolor=MyBlue,linkcolor=MyBlue,citecolor=MyBlue, 
pdftitle={},
pdfauthor={},
pdfsubject={},
pdfkeywords={},
pdfcreator={pdfLaTeX},
pdfproducer={LaTeX with hyperref}
}

\usepackage[sort&compress,capitalize,nameinlink]{cleveref}

\newcommand{\dd}{\:d}

\newcommand{\eps}{\varepsilon}
\newcommand{\from}{\colon}

\newcommand{\by}{y}

\newcommand{\R}{\mathbb{R}}

\newcommand{\N}{\mathcal{N}}

\DeclareMathOperator{\bigoh}{\mathcal O}

\DeclareMathOperator{\ex}{\mathbb{E}}
\DeclareMathOperator{\grad}{\nabla}

\DeclareMathOperator{\prob}{\mathbb{P}}

\DeclareMathOperator{\tr}{tr}

\providecommand\given{} 

\DeclarePairedDelimiter{\braces}{\{}{\}}
\DeclarePairedDelimiter{\bracks}{[}{]}

\DeclarePairedDelimiter{\abs}{\lvert}{\rvert}
\DeclarePairedDelimiter{\norm}{\lVert}{\rVert}
\DeclarePairedDelimiterXPP{\dnorm}[1]{}{\lVert}{\rVert}{_{2}}{#1}

\DeclarePairedDelimiterX{\braket}[2]{\langle}{\rangle}{#1\mathopen{}\delimsize,\mathopen{}#2}

\DeclarePairedDelimiterX{\inner}[2]{\langle}{\rangle}{#1,#2}
\DeclarePairedDelimiterX{\setdef}[2]{\{}{\}}{#1:#2}

\DeclarePairedDelimiterXPP{\exclude}[1]{\mathopen{}\setminus}{\{}{\}}{}{#1}

\DeclarePairedDelimiterXPP{\probof}[1]{\prob}{(}{)}{}{%
	\renewcommand\given{\nonscript\:\delimsize\vert\nonscript\:\mathopen{}}
	#1}

\DeclarePairedDelimiterXPP{\exof}[1]{\ex}{[}{]}{}{%
	\renewcommand\given{\nonscript\:\delimsize\vert\nonscript\:\mathopen{}}
	#1}

\DeclarePairedDelimiterXPP{\trof}[1]{\tr}{[}{]}{}{#1}

\newcommand{\textpar}[1]{\textup(#1\textup)}

\newcommand{\as}{\textup(a.s.\textup)\xspace}

\usepackage[textwidth=30mm]{todonotes}

\newcommand{\debug}[1]{#1}								

\crefname{assumption}{Assumption}{Assumptions}
\crefname{proofstep}{Step}{Steps}
\crefname{lemma}{Lemma}{Lemmas}

\newcommand{\feas}{\mathcal{\debug X}}

\newcommand{\obj}{\debug f}
\newcommand{\sobj}{\debug F}

\newcommand{\sol}{\debug x^{\ast}}

\newcommand{\solset}{\feas^{\ast}}

\newacro{GAN}{generative adversarial network}

\newcommand{\vdim}{\debug d}

\DeclareMathOperator{\mirror}{\mathbf{pr}_{\feas}}

\newcommand{\fench}{\debug E}

\newcommand{\energy}{\fench}

\newcommand{\act}{\debug X}
\newcommand{\actd}{\debug x}
\newcommand{\x}{x}
\newcommand{\y}{y}

\newcommand{\score}{\debug Y}
\newcommand{\scored}{\debug y}

\renewcommand{\v}{\debug \grad\obj}

\newcommand{\start}{\debug 1}
\newcommand{\run}{\debug n}
\newcommand{\iRun}{\debug k}
\newcommand{\step}{\debug\alpha}

\newcommand{\filter}{\mathcal{\debug F}}

\usepackage{algorithm}
\usepackage{algorithmic}

\theoremstyle{plain}
\newtheorem{theorem}{Theorem}		
\newtheorem{lemma}{Lemma}		
\newtheorem{proposition}{Proposition}		

\newtheorem*{corollary*}{Corollary}		

\theoremstyle{definition}
\newtheorem{definition}{Definition}		
\newtheorem{assumption}{Assumption}		

\newtheorem*{definition*}{Definition}		
\newtheorem*{assumption*}{Assumptions}		
\newtheorem*{example*}{Example}		

\theoremstyle{remark}
\newtheorem{remark}{Remark}		

\newtheorem*{remark*}{Remark}		

\newcounter{proofpart}

\newcommand{\flow}{\debug P}

\newcommand{\sample}{\debug\omega}
\newcommand{\samples}{\debug\Omega}

\begin{document}


\newcommand{\longtitle}{\uppercase{Distributed Stochastic Optimization with Large Delays}}		
\title{\longtitle}		

\author[Z.~Zhou]{Zhengyuan Zhou$^{\ast}$}		
\address{$^{\ast}$ Stern School of Business, New York University, NYC, USA.}		
\email{zzhou@stern.nyu.edu}		

\author
[P.~Mertikopoulos]
{Panayotis Mertikopoulos$^{\sharp}$}
\address{$^{\sharp}$\,%
Univ. Grenoble Alpes, CNRS, Inria, Grenoble INP, LIG, 38000, Grenoble, France.}
\email{panayotis.mertikopoulos@imag.fr}

\author[N.~Bambos]{Nicholas Bambos$^{\diamond}$}		
\address{$^{\diamond}$ Department of Management Science and Engineering, Stanford University.}		
\email{bambos@stanford.edu}		

\author[P.~W.~Glynn]{Peter W.~Glynn$^{\diamond}$}		
\email{glynn@stanford.edu}		

\author[N.~Bambos]{Yinyu Ye$^{\diamond}$}		
\email{yinyu-ye@stanford.edu}		



\newcommand{\acdef}[1]{\textit{\acl{#1}} \textup{(\acs{#1})}\acused{#1}}
\newcommand{\acdefp}[1]{\emph{\aclp{#1}} \textup(\acsp{#1}\textup)\acused{#1}}
\newcommand{\acli}[1]{\textit{\acl{#1}}}

\newacro{APT}{asymptotic pseudotrajectory}
\newacroplural{APT}{asymptotic pseudotrajectories}
\newacro{ODE}{ordinary differential equation}
\newacro{FOREL}[FoReL]{follow the regularized leader}
\newacro{iid}[i.i.d.]{independent and identically distributed}
\newacro{KKT}{Karush\textendash Kuhn\textendash Tucker}
\newacro{KL}{Kullback\textendash Leibler}
\newacro{EW}{exponential weights}
\newacro{NE}{Nash equilibrium}
\newacroplural{NE}[NE]{Nash equilibria}
\newacro{RE}{restricted equilibrium}
\newacroplural{RE}[RE]{restricted equilibria}
\newacro{CE}{correlated equilibrium}
\newacroplural{CE}[CE]{correlated equilibria}
\newacro{CCE}{coarse correlated equilibrium}
\newacroplural{CCE}[CCE]{coarse correlated equilibria}
\newacro{GESS}{globally evolutionarily stable state}
\newacro{MSE}{mean squared error}
\newacro{DASGD}[DASGD]{distributed asynchronous stochastic gradient descent}
\newacro{DAGD}[DAGD]{distributed asynchronous gradient descent}
\newacro{VC}[VC]{variational coherence}

\begin{abstract}
One of the most widely used methods for solving large-scale stochastic optimization problems is \ac{DASGD}, a family of algorithms that result from parallelizing stochastic gradient descent on distributed computing architectures (possibly) asychronously. 
However, a key obstacle in the efficient implementation of \ac{DASGD} is the issue of \emph{delays:}
when a computing node contributes a gradient update, the global model parameter may have already been updated by other nodes several times over, thereby rendering this gradient information stale.
These delays can quickly add up if the computational throughput of a node is saturated, so the convergence of \ac{DASGD} may be compromised in the presence of large delays. 
Our first contribution is that, by carefully tuning the algorithm's step-size, convergence to the critical set is still achieved in mean square, even if the delays grow unbounded at a polynomial rate.
We also establish finer results in a broad class of structured optimization problems (called variationally coherent), where we show that \ac{DASGD} converges to a global optimum with probability $1$ under the same delay assumptions.
Together, these results contribute to the broad landscape of large-scale non-convex stochastic optimization by offering state-of-the-art theoretical guarantees and providing insights for algorithm design.
\end{abstract}

\allowdisplaybreaks		
\acresetall		
\maketitle

\section{Introduction}
\label{sec:introduction}

With the advent of high-performance computing infrastructures that are capable of handling massive amounts of data, distributed stochastic optimization has become the predominant paradigm in a broad range of
applications in operations research \citep{dey2017analysis,feng2013online, john, uryasev2010stochastic, shapiro2007tutorial, ruszczynski2003stochastic,lian2016comprehensive}.
Starting with a series of seminal contributions by \citet{TBA86},
recent years have witnessed a commensurate surge of interest in the parallelization of first-order methods, ranging from ordinary (stochastic) gradient descent \citep{NIPS2011_4247, recht2011hogwild, paine2013gpu, chaturapruek2015asynchronous, lian2015asynchronous, feyzmahdavian2016asynchronous,mania2017perturbed},
to coordinate/dual coordinate descent \citep{liu2014asynchronous, avron2015revisiting, liu2015asynchronous,tappenden2017complexity,fercoq2015accelerated,Tran:2015:SUS:2783258.2783412},
randomized Kaczmarz algorithms \citep{liu2014asynchronous},
online methods \citep{QK15,JGS16,HMZ20,HIMM20b},
block coordinate descent \citep{wang2016parallel,mare2015distributed,wright2015coordinate},
ADMM \citep{zhang2014asynchronous,hong2017distributed},
and many others.%

This popularity is a direct consequence of Moore's law of silicon integration and the commensurately increased distribution of computing power.
For instance, in a typical supercomputer cluster, up to several thousands of ``workers'' (or sometimes \emph{tens} of thousands) perform independent computations with little to no synchronization \textendash\ as the cost of such coordination quickly becomes prohibitive in terms of overhead and energy spillage.
Similarly, massively parallel computing grids and data centers (such as those of Google, Amazon, IBM or Microsoft) may house up to several million computing nodes and/or servers, all working asynchronously to execute a variety of different tasks.
Finally, taking the concept of distributed computing to its logical extreme, volunteer computing grids (such as Berkeley's BOINC infrastructure or Stanford's folding@home project) essentially span the entire globe and harness the computing power of a vast, heterogeneous network of non-clustered nodes that receive and process computational requests in a non-concurrent fashion, rendering syncrhonization impossible.
In this way, by eliminating the required coordination overhead, asynchronous operations become simultaneously more appealing (in physically clustered systems) and more scalable (in massively parallel and/or volunteer computing grids).

In this broad context, perhaps the most widely deployed method is \acdef{DASGD} and its variants.
In addition to its long history in mathematical optimization, \ac{DASGD} has also emerged as one of the principal algorithmic schemes for training large-scale machine learning models.
In ``big data'' applications in particular, obtaining first-order information on the underlying learning objective is a formidable challenge, to the extent that the only information that can be readily computed is an imperfect, stochastic gradient thereof \citep{dean2012large, Dean:2012:LSD:2999134.2999271,NIPS2012_4824, zhang2013asynchronous, paine2013gpu, zhang2015deep}.
This information is typically obtained from a group of computing nodes (or processors) working in parallel, and is then leveraged to provide the basis for a distributed descent step.

Depending on the specific computing architecture, the resulting \ac{DASGD} scheme varies accordingly.
More concretely, there are two types of distributed computing architectures that are common in practice:
The first is a cluster-oriented, multi-core, shared memory architecture where different processors independently compute stochastic gradients and update a global model parameter \citep{chaturapruek2015asynchronous,lian2015asynchronous, feyzmahdavian2016asynchronous}.
The second is a ``master-worker'' architecture used predominantly in computing grids (and, especially, volunteer computing grids):
here, each worker node independently \textendash\ and asynchronously \textendash\ computes a stochastic gradient of the objective and sends it to the master;
the master then updates the model's global parameter and sends out new computation requests \citep{NIPS2011_4247,lian2015asynchronous}.
In both cases, \ac{DASGD} is inherently susceptible to \emph{delays}, a key impediment that is usually absent in centralized stochastic optimization settings.
For instance, in a master-worker system,
when a worker sends its gradient update to the master, the master may have already updated the model parameters several times (using updates from other workers), so the received gradient is already stale by the time it is received.
In fact, even in the perfectly synchronized setting where all workers have the same speed and send input to the master in an exact round-robin fashion, there is still a constant delay that grows roughly proportionally to the number of workers in the system \citep{NIPS2011_4247}.

This situation is exacerbated in volunteer computing grids:
here, workers typically volunteer their time and resources following a highly erratic and inconstant update/work schedule, often being turned off and/or being used for different tasks for hours (or even days) on end.
In such cases, there is no lower bound on the fraction of resources used by a worker to compute an update at any given time (this is especially true in heterogeneous computing grids such as BOINC and SimGrid), meaning in turn that there is no upper bound on the induced delays.
This can also happen in parallel computing environments where many tasks with different priorities are executed at the same time across different machines and, likewise, even in multi-core infrastructures with a shared memory, memory-starved processors can become arbitrarily slow in performing gradient computations.

From a theoretical standpoint, the issue of delays and asynchronicities has been studied from the early days of distributed computing \citep{Bertsekas:1997:PDC:548930}, and one of the principal results in the field is the subsequential convergence of \ac{DASGD} with probability $1$ when no constraints are present and when the observed delays grow moderately with time \textendash\ i.e., sublinearly relative to a global clock \citep{Tsi84,TBA86}.
In several current systems (for example, in volunteer computing grids),
as slower workers become saturated and accumulate computation requests over time while new (and possibly faster) workers enter the system, delays can quickly add up and grow at a \emph{superlinear} rate relative to the system's global timer.
Further, in several applications, there are natural constraints imposed on a subset (or all) of the decision variables that represent the model parameters. 
In such contexts, the following questions remain open:
\emph{How robust is the performance envelope of \ac{DASGD} for constrained optimization under large delays and asynchronicities?
Can this robustness be leveraged from a theoretical viewpoint in order to design new and more efficient algorithms?}

\subsection{Our Contributions and Related Work}
Our aim in this paper is to establish the convergence of \ac{DASGD} in the presence of large, superlinear delays, in as wide a class of objectives as possible and in the presence of constraints where efficient projection can be performed.
To that end, we focus on the following classes of problems, where different convergence results can be obtained:

\paragraph{General non-convex objectives\afterhead}
We first consider the class of general smooth non-convex functions, with no structural assumption on the objectives. 
In this (difficult) case, \citet{TBA86} showed that, under \emph{sublinear} delays, \ac{DASGD} converges to a level set of the objective which contains a critical point with probability $1$;
in particular, if every such point is a global minimizer (e.g., if the problem is pseudo-convex) and the method is run with an $\Omega(1/n)$ step-size schedule, \ac{DASGD} converges to the problem's solution set.
More recently, \citet{lian2015asynchronous} derived an estimate for the rate of convergence of the surrogate length $\run^{-1} \sum_{\iRun=1}^{\run} \exof{\norm{\grad \obj(\act_{\iRun})}_{2}^{2}}$ as $\run\to\infty$ under the assumption that the delays affecting the algorithm are \emph{bounded}.
Our first contribution is to show that these assumptions on the delays are not needed:
specifically, as we show in \cref{thm:nonconvex}, by tuning the algorithm's step-size appropriately, it is possible to retain this convergence guarantee,
even if delays grow as polynomials of arbitrary degree.

\paragraph{Variationally coherent objectives\afterhead}
Albeit directly applicable to general non-convex stochastic optimization, \cref{thm:nonconvex} only guarantees convergence to stationarity in the mean square sense;
to ensure global optimality, stronger structural assumptions on the objectives must be imposed.
The ``gold standard'' of such assumptions (and by far the most widely studied one) is convexity:
in the context of distributed stochastic convex optimization, recent works by \citet{NIPS2011_4247} and \citet{recht2011hogwild}, have established convergence for \ac{DASGD} under \emph{bounded} delays for each of the two distributed computing architectures, while \citet{chaturapruek2015asynchronous} and \citet{feyzmahdavian2016asynchronous} extended the bounded delays assumption to a setting with finite-mean \acs{iid} delays. To go beyond this framework, we focus the class of \emph{mean variationally coherent} optimization problems~\cite{zhou2017stochastic, zhou2020}, which includes pseudo-, quasi- and/ star-convex problems, as well as many other classes of functions with non-convex profiles.
Our main result here is that, in such problems, the global state parameter $\act_{\run}$ of \ac{DASGD} converges to a global minimum with probability $1$, even when the delays between gradient updates and requests grow at a polynomial rate (and this, without any distributional assumption on how the underlying delays are generated).

To go beyond this framework, we focus on a class of unimodal problems, which we call \emph{variationally coherent}, and which properly includes all pseudo-, quasi- and/ star-convex problems,
as well as many other classes of functions with highly non-convex profiles.
Our main result here may be stated as follows:
in stochastic variationally coherent problems, provided a lazy descent scheme is used (akin to dual averaging) to mesh with the constraint set in the distributed procedure, the global state parameter $\act_{\run}$ of \ac{DASGD} converges to a global minimum with probability $1$, even when the delays between gradient updates and requests grow at a polynomial rate (and this, without any distributional assumption on how the underlying delays are generated).

This result extends the works mentioned above in several directions:
specifically, it shows that
\begin{enumerate}
[\indent 1.]
\item
Convexity is not required to obtain global convergence results.
\item 
Constraints do not hinder almost sure convergence under a suitable lazy projection scheme.
\item
The robustness of \ac{DASGD} is guaranteed even under large, superlinear delays.
\end{enumerate}
We find these outcomes particularly appealing because, coupled with the existing rich literature on the topic, they help explain and reaffirm the prolific empirical success of \ac{DASGD} in large-scale machine learning problems, and offer concrete design insights for fortifying the algorithm's distributed implementation against delays and asynchronicities.

\paragraph{Techniques\afterhead}
Our analysis relies on techniques and ideas from stochastic approximation and \mbox{(sub-)}mar\-tin\-gale convergence theory.
A key feature of our approach is that, instead of focusing on the discrete-time algorithm, we first establish the convergence of an underlying, deterministic dynamical system by means of a particular energy (Lyapunov) function which is decreasing along continuous-time trajectories and ``quasi-decreasing'' along the iterates of \ac{DASGD}.
To control this gap, we connect the continuous- and discrete-time frameworks via the theory of \acdefp{APT}, as pioneered by \citet{BH96}.
By itself, the \ac{APT} method does not suffice to establish convergence under delays.
However, if the step-size of the method is chosen appropriately (following a quasi-linear decay rate for polynomially growing delays), it is possible to leverage $L^{p}$ martingale tail convergence results to show that the problem's solution set is recurrent under \ac{DASGD}.
This, combined with the above, allows us to prove our core convergence results.

Even though the \ac{ODE} approximation of discrete-time Robbins\textendash Monro algorithms has been widely studied in control and optimization theory \citep{kushner2013stochastic, ghadimi2013stochastic}, transferring the convergence guarantees of an \ac{ODE} solution trajectory to a discrete-time algorithm is a fairly subtle affair that must be done on a case-by-case basis. 
Further, even if this transfer is complete, the results typically have the nature of convergence-in-probability:
almost-sure convergence is usually much harder to obtain \cite{Bor08}.
Specifically, exisiting stochastic approximation results
cannot be applied to our setting because
\begin{inparaenum}
[\itshape a\upshape)]
\item
the non-invertibility of the projection map makes the underlying dynamical system on the problem's feasible region non-autonomous (so \ac{APT} results do not apply);
\item
unbounded delays only serve to aggravate this issue, as they introduce a further disconnect between the \ac{DASGD} algorithm and its continuous-time version.
\end{inparaenum}
To control the discrepancy between discrete and continuous time requires a more fine-grained analysis, for which we resort to a sharper law of large numbers for $L^{p}$-bounded martingales.
Finally, we also mention that the recent work of \citet{lian2016comprehensive} has also considered distributed zeroth-order methods (where only the function value, rather than the gradient is available)
and used techniques that are different from gradient-based analyses. 

\section{Problem Setup}
\label{sec:model}

Let $\feas$ be a subset of $\R^\vdim$ and let $(\samples,\filter,\prob)$ be some underlying (complete) probability space.
Throughout this paper, we focus on the following stochastic optimization problem:
\begin{equation}
\label{eq:opt}
\tag{Opt}
\begin{aligned}
\text{minimize}
&\quad
\obj(\actd)
\\
\text{subject to}
&\quad
\actd\in\feas,
\end{aligned}
\end{equation}
where the objective function $\obj\from\feas\to\R$ is of the form
\begin{equation}\label{eq:obj}
\obj(\actd)
= \exof{\sobj(\actd;\sample)}
	= \int_{\Omega}\sobj(\actd; \sample) \dd\probof{\sample}
\end{equation}
for some random function\footnote{It is understood that a random function here means that $\sobj(\actd;\cdot)\from\samples\to\R$ is a measureable function for each $\actd \in \feas$.} $\sobj\from\feas\times\samples\to\R$.
Using standard optimization terminologies, \eqref{eq:opt} is called an unconstrained stochastic optimization problem if $\feas = \R^\vdim$, and is called a constrained stochastic optimization problem otherwise. 
In this paper, we study both unconstrained and constrained stochastic problems under smooth	objectives.
Specifically, we make the following regularity assumptions for the rest of the paper, which are standard in the literature:
\begin{assumption}\label{asm:basic}
	We assume the following:
	\begin{enumerate}
		\item
		\label{assump:0}
		$\sobj(x;\sample)$ is differentiable in $x$ for $\prob$-almost all $\sample\in\samples$.%
		\item
		\label{assump:0.5}
		$\grad\sobj(x;\sample)$ has\footnote{It is understood here that the gradient $\grad\sobj(x;\sample)$ is only taken with respect to $x$: no differential structure is assumed on $\samples$.} finite second moment, that is,
		$\sup_{\actd \in \feas}\exof{\dnorm{\nabla\sobj(x;\sample)}^{2}} < \infty$.
		\item
		\label{assump:0.6}
		$\nabla\sobj(x;\sample)$ is Lipschitz continuous in the mean:
		$\exof{\nabla\sobj(x;\sample)}$ is Lipschitz on $\feas$.
	\end{enumerate}
\end{assumption}

\begin{remark}\label{rem:basic}
	\quad
Assumptions \ref{assump:0} and \ref{assump:0.5} together imply that $\obj$ is differentiable, because finite second moment (by Statement~\ref{assump:0.5}) implies finite first moment: $\exof{\dnorm{\nabla\sobj(x;\sample)}} < \infty$ for all $x\in\feas$; and hence the expectation $\exof{\nabla\sobj(x;\sample)}$ exists.
By a further application of the dominated convergence theorem, we have $\nabla\obj(x) = \nabla\exof{\sobj(x;\sample)} = \exof{\nabla\sobj(x;\sample)}$. Additonally, Assumption~\ref{assump:0.6}  implies that $\nabla\obj$ is Lipschitz continuous.
In the deterministic optimization literature, $\obj$ is sometimes called $L$-smooth, where $L$ is the Lipschitz constant.
\end{remark}

One important class of motivating applications that can be cast in the current stochastic optimization problem~\eqref{eq:opt} is empirical risk minimization (ERM) in machine learning. 
As is well-known in the distributed optimization/learning literature \cite{NIPS2011_4247,NIPS2012_4824, zhang2013asynchronous, paine2013gpu,lian2015asynchronous}, the expectation in \cref{eq:obj} contains as a special case the common machine learning objectives of the form $\frac{1}{N} \sum_{i=1}^N \obj_i(\actd)$, where each  $\obj_i(\actd)$ is the loss associated with the $i$-th training sample.
 This setup corresponds to ERM without regularization.
With regularization, ERM takes the form $\frac{1}{N} \sum_{i=1}^N \obj_i(\actd) + r(\actd)$,
where $r(\cdot)$ is a regularizer (typically convex and known), which is again a special case of~\eqref{eq:opt}. Other related examples that are also special cases of~\eqref{eq:opt} include the objective $\sum_{i=1}^N v_i \obj_i(\actd)$, which are standard in curriculum learning: $v_i$ are weights (between $0$ and $1$) generated from a learned curriculum that indicates how much emphasis each loss $f_i$ should be given. 

In the large-scale data setting ($N$ is very large), such problems are typically solved in practice using stochastic gradient descent (SGD) on a distributed computing architecture. 
SGD\footnote{In vanilla SGD (i.e. centralized/single-processor setting),  an \textbf{iid} sample of the gradient of $\sobj$ at the current iterate is used to make a descent step (with an appropriate projection made in the constrained optimization case).} is widely used primarily because in many applications, stochastic gradient, computed by first drawing a sub-batch of data and then computing the average of the gradients on that sub-batch, is the only type of information that is computationally feasible and practically convenient to obtain\footnote{For instance, Google's Tensorflow system does automatic differentiation on samples for neural networks (i.e. $f_i(\cdot)$'s are parametrized neural networks).}.
Further, such problems are generally solved on a distributed computing architecture because:
1) Computing gradients is typically the computational bottleneck. Consequently, having multiple processors compute gradients in parallel can harness the available computing power in modern distributed systems.
2) The data (which determine the individual cost functions $f_i$'s) may simply be too large to fit on a single machine; and hence a distributed system is necessary even from a storage perspective.

\begin{figure}[t]
	\centering
	\begin{tabular}{ccc}
		\includegraphics[width=0.45\linewidth]{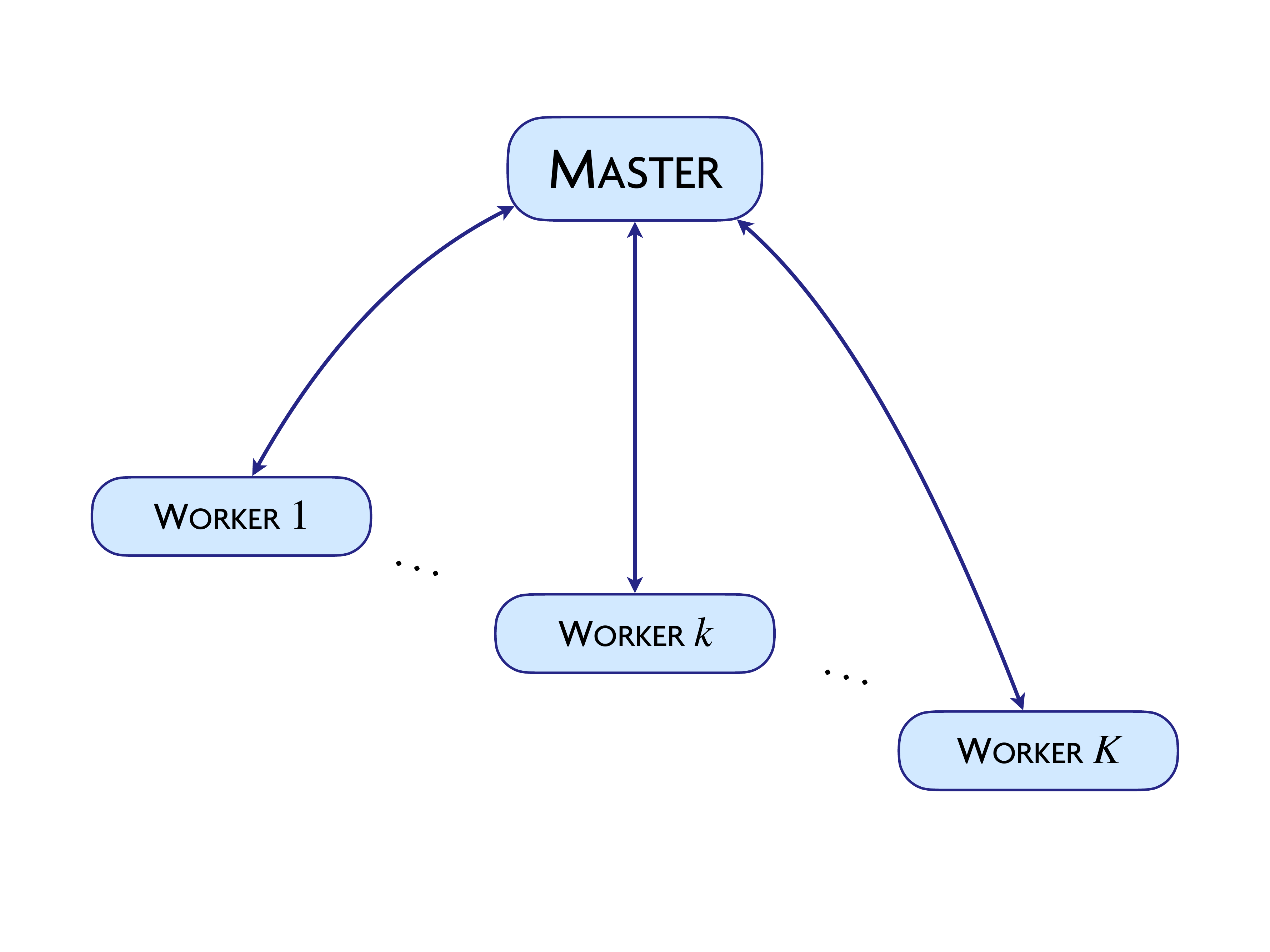}&
		\includegraphics[width=0.45\linewidth]{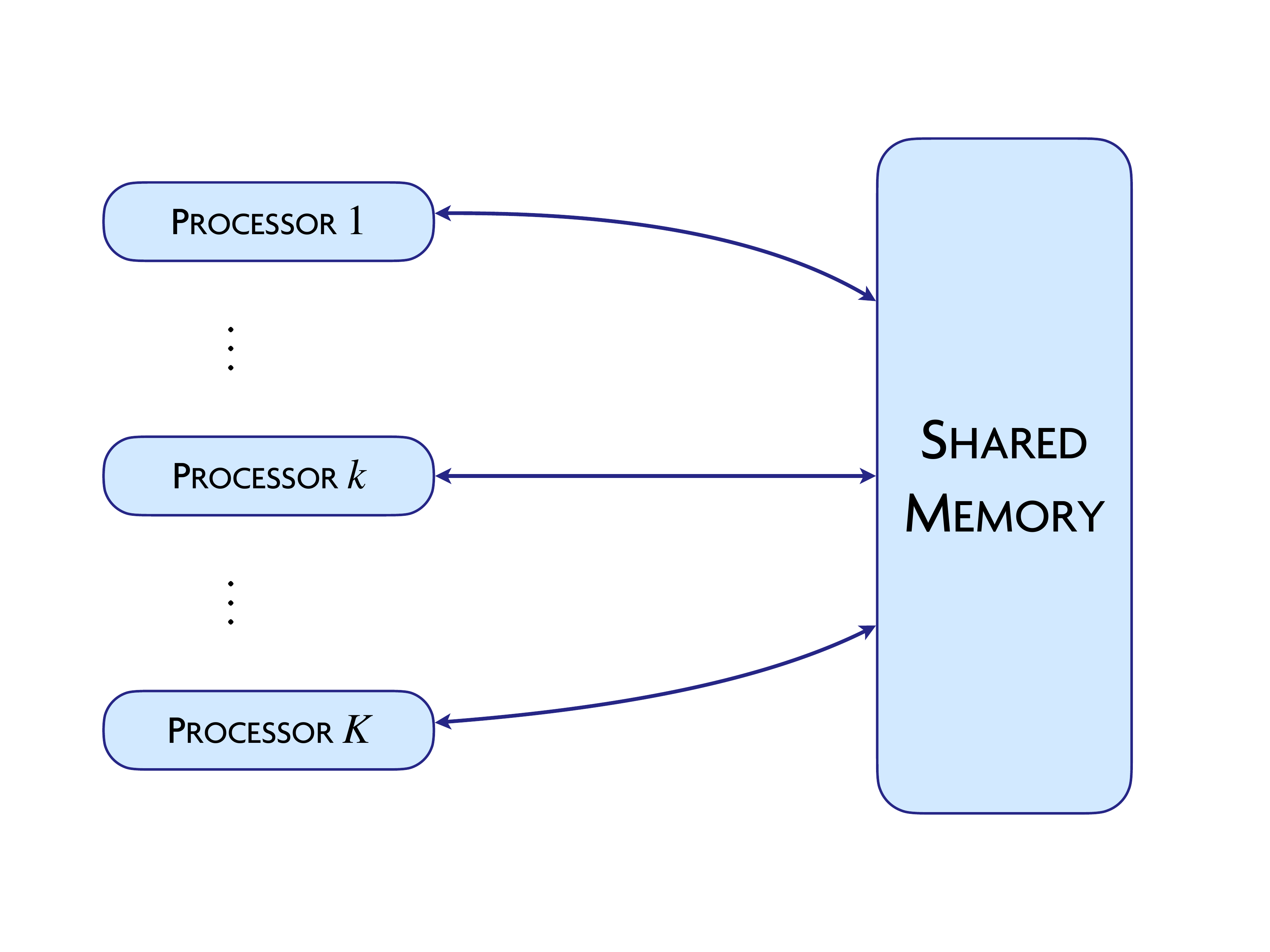}\\
	\end{tabular}
	\vspace{-2ex}
	\caption{Two commonly used distributed computing architecutures: (a) master-worker (left) and multi-processorwith shared memory (right).}
	\label{fig:ms}
	\vspace{-2ex}
\end{figure}

With the above background, our goal in this paper is to study
and establish theoretical convergence guarantees for applying stochastic gradient descent (SGD) to solve~\eqref{eq:opt} on a distributed computing architecture. Two common distributed computing architectures that are widely delpoyed in practice (see also \cref{fig:ms}): 
\begin{enumerate}
	\item \textbf{Master-worker system.}
	This architecture is mostly used in data-centers and parallel computing grids (each computing node is a single machine, virtual or physical).
	\item \textbf{Multi-processor system with shared memory.}
	This architecture is mostly used in multi-core machines or GPU computing: in the former, each processor is a CPU, while in the latter, each processor is a GPU.
\end{enumerate}

In the next two subsections (\cref{subsec:master_slave,subsec:multi_core}), we describe the standard procedure of parallelizing SGD
on each of the two distributed computing architectures.
Although running SGD on these two architectures have some differences, in \cref{subsec:unify} we give a meta algorithmic description, called \acdef{DASGD}  that unifies these two parallelizations on the same footing.

\subsection{SGD on Master-Worker Systems}\label{subsec:master_slave}
Here we consider the first distributed computing architecture: the master-worker system. 
The standard way of deploying stochastic gradient descent in such systems \textendash\ and that which we adopt here \textendash\ is for the workers
to asychronously compute stochastic gradients and then send them to the master,%
\footnote{As alluded to before, in machine learning applications, this is done by sampling a subset of the training data, computing the gradient for each datapoint and averaging over all datapoints in the sample.}
while the master updates the global state of the system and pushes the updated  state back to the workers (\cite{NIPS2011_4247,lian2015asynchronous}).
This process is presented in \cref{alg:MSSGD}.

\begin{algorithm}[htbp]
	\caption{Running SGD on a Master-Worker Architecture} 
	\label{alg:MSSGD}
	\begin{algorithmic}[1]
		\REQUIRE $1$ Master and $K$ workers, $k=1,\dotsc,K$
		\STATE Each worker is seeded with the same inital iterate
		\REPEAT 
		\STATE \textbf{Master}:
		
		(a) Receive a stochastic gradient from worker $k$
		
		(b) Update current iterate.
		
		(c) Send updated iterate to worker $k$
		
		\STATE \textbf{Workers}:
		
		(a) Receive iterate
		
		(b) Compute an i.i.d. stochastic gradient (at the received iterate)
		
		(c) Send the computed gradient to master

		\UNTIL{end}
	\end{algorithmic}
\end{algorithm}

Due to the distributed nature of the master-worker system, 
a gradient received by the master on any given iteration can be stale. As a simple example, consider a fully coordinated update scheme where each worker sends the computed gradient to and receives the updated iterate from the master following a round-robin schedule. In this case, each worker's gradient is received with a delay exactly equal to $K-1$ ($K$ is the number of workers in the system), because by the time the master receives worker $K$'s computed gradient, the master has already applied $K-1$ gradient updates from workers $1$ to $K-1$ (and since the schedule is round-robin, this delay of $K-1$ is true for any one of the $K$ workers).
 
However, delays can be much worse since we allow 
full asynchrony: workers can compute and send (stochastic) gradients to the master without any coordinated schedule.  
In the asynchronous setting, fast workers (workers that are fast in computing gradients) will cause disproprotionately large delays to gradients produced by slow workers (workers that are slow in computing gradients): when a slower worker has finished computing a gradient, a fast worker may have already computed and communicated many gradients to the master. 
Since the master updates the global state of the system (the current iterate), one can gain a clearer representation of this scheme by looking at the master's update. This is given in \cref{subsec:unify}.

\subsection{SGD on Multi-Processor Systems with Shared Memory}\label{subsec:multi_core}

Here we consider the second distributed computing architecture: multi-processor system with shared memory. 
In this architecture, all processors can access a global, shared memory, which holds all the data needed for computing a (stochastic) gradient, as well as the current iterate (the global state of the system).
The standard way of deploying stochastic gradient descent in such systems (\cite{chaturapruek2015asynchronous,lian2015asynchronous}) is for each processor to independently and asychronously read the current global iterate, compute a stochastic gradient~\footnote{This is again done by sampling a subset of the training data in the global memory and computing the gradient at the iterate for each datapoint and averaging over all the comptued gradients in the sample.}, and then update the global iterate in the shared memory.
This process is given \cref{alg:MSSGD1}:

\begin{algorithm}[H]
	\caption{Running SGD on a Multi-Processor System with Shared Memory} 
	\label{alg:MSSGD1}
	\begin{algorithmic}[1]
		\REQUIRE $K$ processors and global (shared) memory.
		\STATE The initial iterate in the global memory. 
		\REPEAT 
		\STATE
		(a) Each proceesor reads the current global iterate.
		
		(b) Each processor reads data from memory and computes a stochastic gradient.
		
		(c) Each processor updates the global iterate.
		
		\UNTIL{end}
	\end{algorithmic}
\end{algorithm}

The key difference from \cref{alg:MSSGD} is that
there is no central entity that updates the global state; instead, each processor 
can both read the global state and update it. Since each processor is performing the operations
asynchronously, different processors may be reading the same global iterate at the same time. Further, the delays
in this case is again caused by the heterogeneity across different processors: if a processor is slow in computing gradients,
then by the time it finishes computing its gradient, the global iterate has been updated by other, faster processors many times over, thereby causing its own gradient stale. Here, we also adopt a common assumption that updating the global iterate is an atomic operation (and hence no two processors will be updating the global iterate at the same time). This is justified\footnote{
	On a related note, we also note that our analysis can be further extended to cases where only one variable or a small block of variables are being updated at a time. We omit this discussion because the resulting notation is quite onerous, and will obscure the main ideas behind the already complex theoretical framework developed here.} because performing gradient update is a simple arithemtic operation, and hence takes negligible time compared to reading data and computing a stochastic gradient, which is the main computational bottleneck in practice. 
However, despite a cheap computation, performing the whole gradient update (typically achieved via locking) does have overhead. In particular, a less stringent model would be only updating one coordinate at a time: this is known as an inconsistent write/read model since different processors are updating different components of the global parameters and hence can read in elements of different ages. For simplicity, we do not consider this case, as our focus in this paper  is on delays. See~\cite{lian2015asynchronous} and \cite{liu2014asynchronous}
for lucid discussions and analyses on this model.
Finally, one can gain a clearer picture of this update scheme by tracking the update to the global iterate in the shared memory. This is given in \cref{subsec:unify}.

\subsection{DASGD: A Unifying Algorithmic Representation}
\label{subsec:unify}
In this subsection, we present a unified algorithmic description,  aptly called \acdef{DASGD}, that formally captures both \cref{alg:MSSGD} and \cref{alg:MSSGD1}, where their differences are reflected in the assumptions of the meta algorithm's parameters. We start with the unconstrained case, see \cref{alg:DASGD_un}.

\begin{algorithm}[htbp]
	\caption{Distributed asynchronous stochastic gradient descent} 
	\label{alg:DASGD_un}
	\begin{algorithmic}[1]
		\REQUIRE Initial state $\act_{0}\in \R^\vdim$, step-size sequence $\step_{\run}$
		\STATE $\run\leftarrow0$;
		\REPEAT
		\STATE $\act_{\run+1} = \act_{\run} - \step_{\run+1} \nabla \sobj(\act_{s(\run)},\sample_{\run+1})$;
		\STATE $\run\leftarrow \run+1$;
		\UNTIL{end}
		\STATE \textbf{return} solution candidate $\act_{\run}$
	\end{algorithmic}
\end{algorithm}

\begin{algorithm}[htbp]
	\caption{Distributed asynchronous stochastic gradient descent with projection} 
	\label{alg:DASGD}
	\begin{algorithmic}[1]
		\REQUIRE Initial state $\score_{0}\in \R^\vdim$, step-size sequence $\step_{\run}$
		\STATE $\run\leftarrow0$;
		\REPEAT
		\STATE $\act_{\run} = \mirror(\score_{\run})$;
		\STATE $\score_{\run+1} = \score_{\run} - \step_{\run+1} \nabla \sobj(\act_{s(\run)},\sample_{\run+1})$;
		\STATE $\run\leftarrow \run+1$;
		\UNTIL{end}
		\STATE \textbf{return} solution candidate $\act_{\run}$
	\end{algorithmic}
\end{algorithm}

In more detail, $\run$ is a global counter and is incremented every time an update occurs to the current solution candidate $\act_{\run}$ (the global iterate): in the master-worker systems, the master updates it; in the multi-processor systems, each processor updates it.
Since there are delays in both systems, the gradient applied to the current iterate $\act_{\run}$ can be a gradient associated with a previous time step. This fact is abstractly captured by Line 3 in \cref{alg:DASGD}.
In full generality, we will write $s(\run)$ for the iteration from which the gradient received at time $\run$ originated.
In other words, the delay associated with iteration $s(\run)$ is $\run - s(\run)$, since it took $\run - s(\run)$ iterations for the gradient computed on iteration $s(\run)$ to be received at stage $\run$. Note that $s(\run)$ is always no larger than $\run$; and if $\run = s(\run)$, then there is no delay in iteration $\run$.

Now, the difference between the two distributed computing archecitures is reflected in the assumption of $s(\run)$.
Specifically, in the master-worker systems, each $s(\cdot)$ is a one-to-one function\footnote{Except initially if all the workers have the same initial point.}, because no two workers will ever receive same iterates from the master per \cref{alg:MSSGD}.
On the other hand, in multi-processor systems, $s(\run)$ can be the same for different $\run$'s (since different processors may read the current iterate at the same time); however, it is easy to observe that the same $s$ will appear at most $K$ times for different $\run$'s, since there are $K$ processors in total.
As an important note, our analysis is \textbf{agnostic} to whether $s(\run)$ is one-to-one or not. Consequently, 
in establishing theoretical guarantees for the meta algorithm \ac{DASGD}, we obtain the same guarantees for both architectures simultaneously.

Notation-wise, we will write $d_{\run}$ for the delay required to compute a gradient requested at iteration $\run$.
This gradient is received at iteration $\run+d_{\run}$. Following this notation, the delay for a gradient received at $\run$ is $d_{s(\run)} = \run - s(\run)$. Note also we have chosen the subscript associated with $\sample$ to be $\run+1$: we can do so because $\sample_{\run}$'s are \textbf{iid} (and hence the indexing is irrelevant). 
Finally, in constrained optimization case (where $\feas$ is a strict subset of $\mathbb{R}^{\vdim}$), projection must be performed. This results in $\ac{DASGD}$ with projection\footnote{This tpe of projection is technically known as lazy projection.}, which is formally given in \cref{alg:DASGD}.

\section{General Nonconvex Objectives}
\label{sec:nonconvex}

In this section, we take $\feas = \mathbb{R}^{\vdim}$ (i.e. unconstrained optimization)
and consider general non-convex objectives. Note that for a general non-convex objective where no further structural assumption (e.g. convexity) is made, convergence to an optimal solution (even a local optimal solution) cannot be expected, and will not hold in general, even in the absence of both noise and delays (i.e. single machine deterministic optimization). 
In such cases, the best one can hope for, which is also the standard metric to determine the stability of the algorithm, is that the gradient vanishes in the limit\footnote{An alternative phrase that is commonly used is that the criticality gap vanishes. This is also colloquially referred to as convergence to a stationary point/critical point in the machine learning community. Note further that convergence to second-order-stationary points can be achieved by stochastic gradient descent (and its variants) under weaker assumptions (than convexity): Lipschitz Hessian and strict saddle point property. See~\cite{ge2015escaping, lee2016gradient,jin2017escape} for this line of work.}.

\subsection{Delay Assumption}
Our goal here is to establish convergence guarantees (in mean square) of \ac{DASGD} for a general non-convex objective in the presence of delays. 
In fact, a large family of unbounded delay processes can be tolerated. We state our main assumption regarding the delays and step-sizes:
\begin{assumption}
	\label{asm:delays}
	\quad
	The gradient delay process $d_{\run}$ and the step-size sequence $\step_{\run}$ of \ac{DASGD} (\cref{alg:DASGD_un,alg:DASGD}) satisfy one of the following conditions:
	\vspace{2ex}
	\begin{enumerate}
		\addtolength{\itemsep}{1ex}
		\item
		\emph{Bounded delays:}
		$\sup_{\run} d_{\run} \le D$ for some positive number $D$
		and
		$\sum_{\run=1}^{\infty} \step_{\run}^{2} < \infty$,
		$\sum_{\run=1}^{\infty} \step_{\run} = \infty$.
		\item
		\emph{Sublinearly growing delays:}
		$d_{\run} = O(\run^p)$ for some $0 < p < 1$
		and
		$\step_{\run} \propto 1/\run$.
		\item
		\emph{Linearly growing delays:}
		$d_{\run} = \bigoh(\run)$
		and
		$\step_{\run} \propto 1/(\run \log\run)$.
		\item
		\emph{Polynomially growing delays:}
		$d_{\run} = \bigoh(\run^{q})$ for some $q\geq1$
		and
		$\step_{\run} \propto 1/(\run \log\run \log\log\run)$.
	\end{enumerate}
\end{assumption}

Note that as delays get larger, we need to use less aggressive step-sizes. This is to be expected, 
because the larger the delays, the more ``averaging" one needs perform in order to remove the staleness that is caused by
the delays; and smaller step-sizes correspond to averaging over a longer horizon. This is a one of the important insights that occur throughout the paper. Another thing to note that is the Assumption~\ref{asm:delays} also highlights the quantitative relationship between the class of delays and the class of step-sizes. For instance, when the delays increase from a linear rate to a polynomial rate, only a factor of $\frac{1}{\log\log \run}$ needs to be added (which is effectively a constant).
From a practical standpoint, this means that a step-size on the order of $1/(\run \log\run)$ will be a good model-agnostic choice and more-or-less sufficient for almost all delay processes.

\subsection{Controlling the Tail Behavior of Second Moments}
We now turn to establish the theoretical convergence guarantees of \ac{DASGD} for general non-convex objectives.
Our first step lies in controlling the tail behavior of the second moments of the gradients that are generated from \ac{DASGD}.
By leveraging the Lipschitz continuity of the gradient, its telescoping sum, appropriate conditionings and a careful analysis of the interplay between delays and step-sizes, 
we show that (next lemma) a particularly weighted tail sum of the second moments are vanishingly small in the limit (see appendix for a detailed proof).

\begin{lemma}\label{lem:tail}\quad
	Under \cref{asm:basic,asm:delays}, if $\inf_{\actd \in \feas}\obj(\actd) > - \infty$, then
	\begin{equation}\label{eq:tail}
	\sum_{\run = 0}^{\infty} \step_{\run+1}\exof{\|\grad f(X_{\run})\|_2^2} < \infty.
	\end{equation}
\end{lemma}

\begin{remark}\quad
	Since $\sum_{\run = 0}^{\infty} \step_{\run+1} = \infty$, \cref{lem:tail} implies that
	$\lim\inf_{\run \rightarrow \infty} \exof{\|\grad f(X_{\run})\|_2^2} = 0$ (see Appendix A of~\cite{bertsekas1996neuro}).
Note that  the converse is not true: when a subsequence of $\exof{\|\grad f(X_{\run})\|_2^2}$ converges to $0$,
	the sum in Equation~\eqref{eq:tail} need not be finite. As a simple example, consider
	$\step_{\run+1} = \frac{1}{\run}$, and
	\begin{equation}
	\exof{\|\grad f(X_{\run})\|_2^2} =
	\begin{cases}
	\frac{1}{\run}, \text{ if } \run = 2^k\\
	1, \text{ otherwise}.
	\end{cases}
	\end{equation}
	Then the subsequence on indicies $2^{\run}$ converges to $0$, but the sum still diverges.
	Consequently, \cref{lem:tail} is stronger than subsequence convergence.
\end{remark}

\subsection{Bounding the Successive Differences}

However, \cref{lem:tail} is still not strong enough to guarantee that 
$\lim_{\run \rightarrow \infty} \exof{\|\grad f(X_{\run})\|_2^2} = 0$.
This is because the convergent sum given in Equation~\eqref{eq:tail} only limits the tail growth somewhat, but 
not completely. To demonstrate this point, let $c_t$ be the following boolean variable: 
\begin{equation}
c_{\run} =
\begin{cases}
1, \text{ if $\run$ contains the digit $9$ in its decimal expansion,} \\
0, \text{ otherwise}.
\end{cases}
\end{equation}
For instance, $c_9 = 1$, $c_{11} = 0$.
Now define 
$\step_{\run+1} = \frac{1}{\run}$, and
$	\exof{\|\grad f(X_{\run})\|_2^2}=
\begin{cases}
\frac{1}{\run}, \text{ if } c_{\run} = 1\\
1, \text{ if } c_{\run} = 0.
\end{cases}$
As first-year Berkeley Math PhDs delightfully found out during-- or painfully found out after -- their qualifying exam,  $\sum_{\run = 1}^{\infty} \step_{\run+1} 	\exof{\|\grad f(X_{\run})\|_2^2}< \infty$ (see Problem 1.3.24 in~\cite{de2012berkeley}), even though the limit $\exof{\|\grad f(X_{\run})\|_2^2}$ does not exist. This indicates that to obtain convergence of $\exof{\|\grad f(X_{\run})\|_2^2}$, we need to impose more stringent conditions to ensure its sufficient tail decay. Note that one issue that is revealed by the above counter-example is that the distance between two successive terms is always bounded away from $0$, no matter how large $\run$ is. This obviously makes it possible
for convergence to occur: a necessary condition for convergence is that the difference converges to $0$.
Note that intuitively, this pathological case cannot occur for gradient descent, because the step-size is shrinking to $0$, hence making the successive difference converge to $0$ (at least in expectation).
Consequently, 
we can bound the difference between every two successive terms in terms of a decreasing sequence that is converging to $0$.
This ensures that $\exof{\|\grad f(X_{\run})\|_2^2}$ cannot change two much from iteration to iteration. Further, the change between two successive terms will be vanishing.
This result is formalized in the following lemma (the proof is given in the appendix):

\begin{lemma}\label{lem:diff}\quad
	Under \cref{asm:basic,asm:delays},
	there exists a constant $C > 0$ such that for every $\run$, 
	$$\Big|\exof{\|\grad f(X_{\run+1})\|_2^2} - 	\exof{\|\grad f(X_{\run})\|_2^2}\Big| \le C \alpha_{\run+1}.$$
\end{lemma}

\subsection{Main Convergence Result}

Putting the above two characterizations together, we obtain:
\begin{theorem}\label{thm:nonconvex}\quad
	Let $\act_{\run}$ be the \ac{DASGD} iterates generated from \cref{alg:DASGD_un}.
	Under \cref{asm:basic,asm:delays}, if $\inf_{\actd \in \feas}\obj(\actd) > - \infty$, then
	\begin{equation}\label{eq:tail}
	\lim_{\run \to \infty} \exof{\|\grad f(X_{\run})\|_2^2} = 0.
	\end{equation}
\end{theorem}

\begin{remark}\quad
Three remarks are in order here. 
First, note that the condition $\inf_{\actd \in \feas}\obj(\actd) > - \infty$ means that the optimization problem~\eqref{eq:opt}
has a solution. This is necessary, because otherwise, a stationary point may not exist in the first place, and \ac{DASGD} (or even simple gradient descent) may continue decrease the objective value \textit{ad infinitum}.
Note that since $\obj$ is smooth, a minimum point is necessarily a stationary point.

Second, the above convergence is a fairly strong characterization of the fact that the gradient vanishes. In particular,
it means that the gradient of the \ac{DASGD} iterates converges to $0$ in mean square. Consequently, this implies that the norm of the gradient vanishes in expectation, and that the gradient converges to $0$ with high probability. Note further that 
if we strengthen Assumption~\ref{asm:basic} to require that all stochastic gradients are bounded almost surely (as opoosed to just bounded in second moments), then a similar analysis ensures almost sure convergence of
$\|\grad f(X_{\run})\|_2$. We omit the details due to space limitation.
Finally, that \cref{thm:nonconvex} is a direct consequence of \cref{lem:tail,lem:diff}
is a simple excercise in elementary series theory (in particular, \cref{lem:A3}).
\end{remark}

\section{Variationally Coherent Problems}
\label{sec:vc}

In this section, our goal is to establish global optimality convergence guarantees of \ac{DASGD} in as wide a class of optimization problems as possible.
Since global convergence cannot be expected to hold for all non-convex optimization problems (even without delays).
a standard structural assumption to make in the existing literature on the objectives (even in the no-delay case) is convexity.
Here we instead consider a much broader class of stochastic optimization problems than convex problems.
Further, we allow for constrained optimization; in particular, 
$\feas$ is assumed to be a convex and compact subset of $\mathbb{R}^{\vdim}$ throughout the section.
We first discuss the class of objectives and then present global convergence results and their analyses in the subsequent two subsections. We focus on the class of mean variationally coherent optimization problems~\cite{zhou2020,MZ19,MLZF+19,zhou2018distributed}, defined here as follows:

\begin{figure}[tbp]
	\centering
	\subfigure{\includegraphics[width=0.4\textwidth]{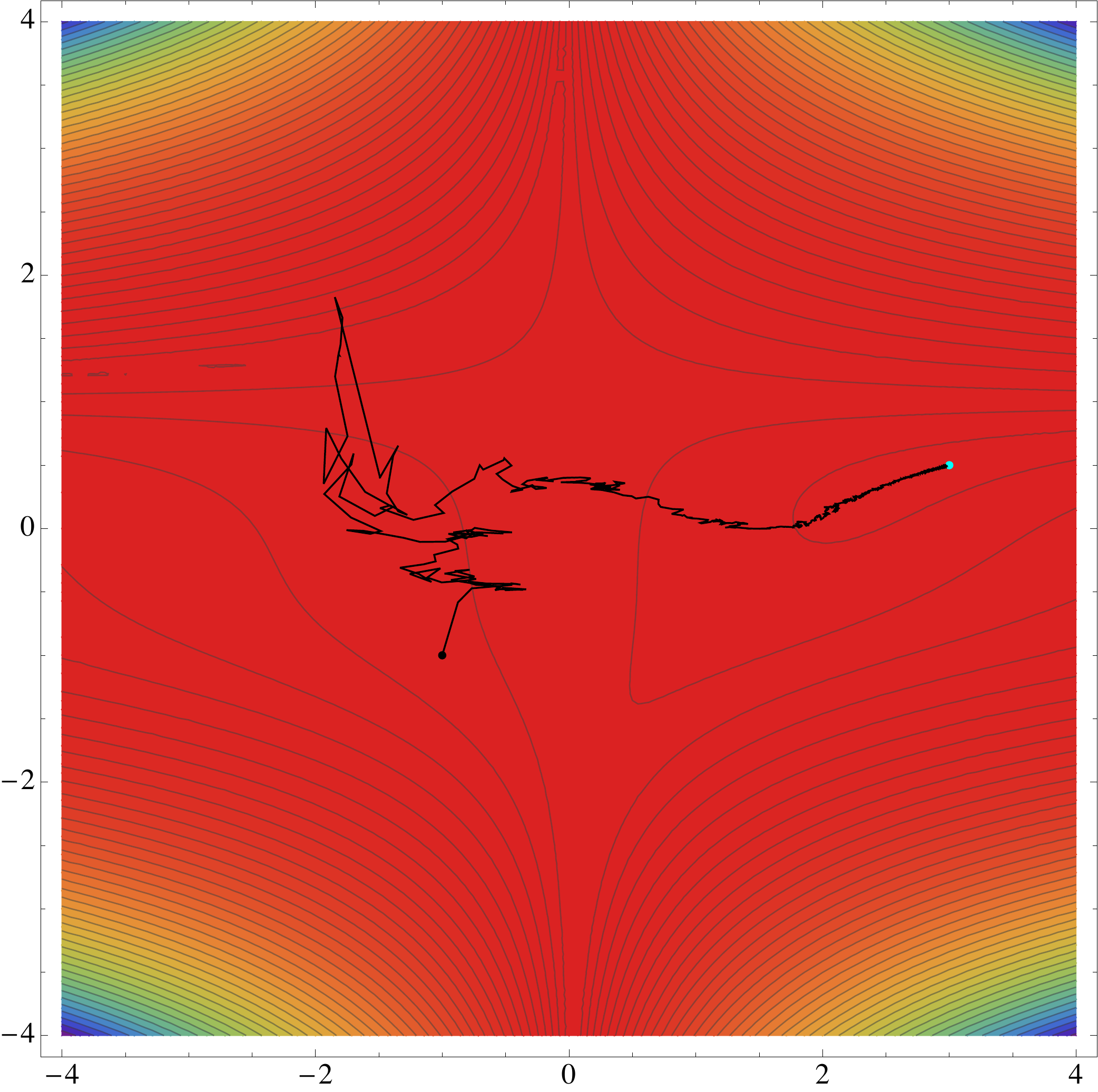}}
	\hfill
	\subfigure{\includegraphics[width=0.45\textwidth]{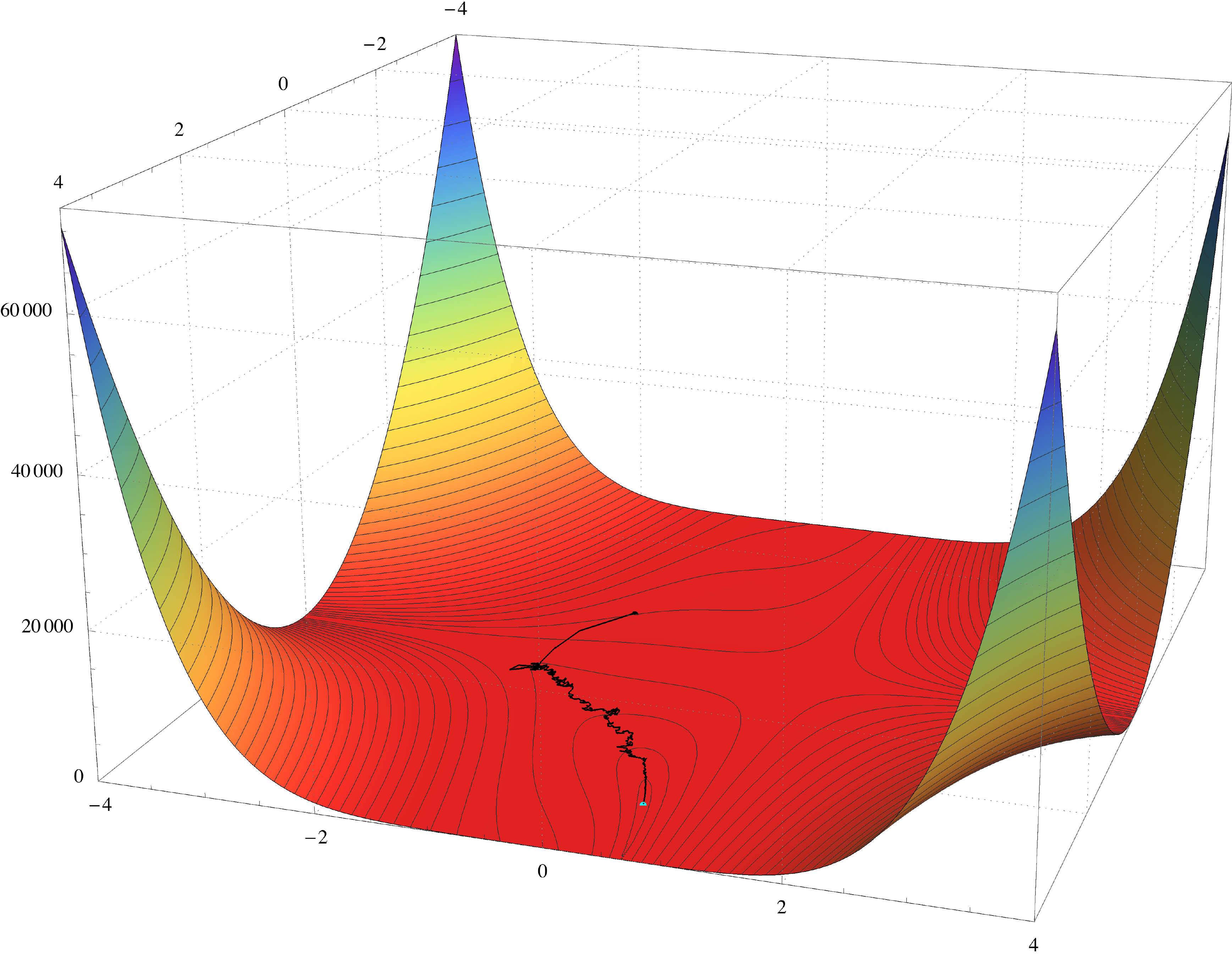}}%
	\\[1ex]
	\subfigure{\includegraphics[width=0.4\textwidth]{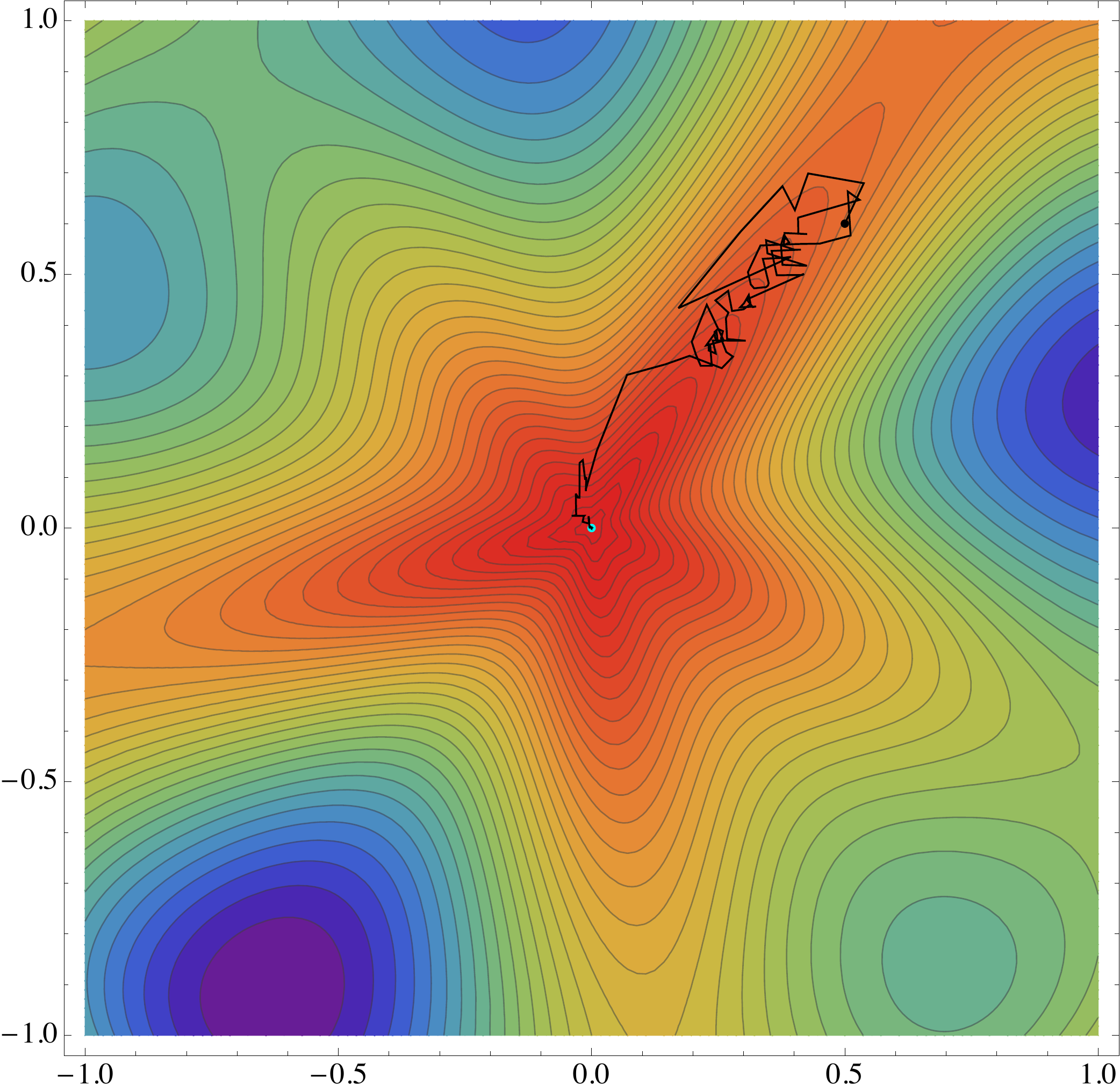}}
	\hfill
	\subfigure{\includegraphics[width=0.45\textwidth]{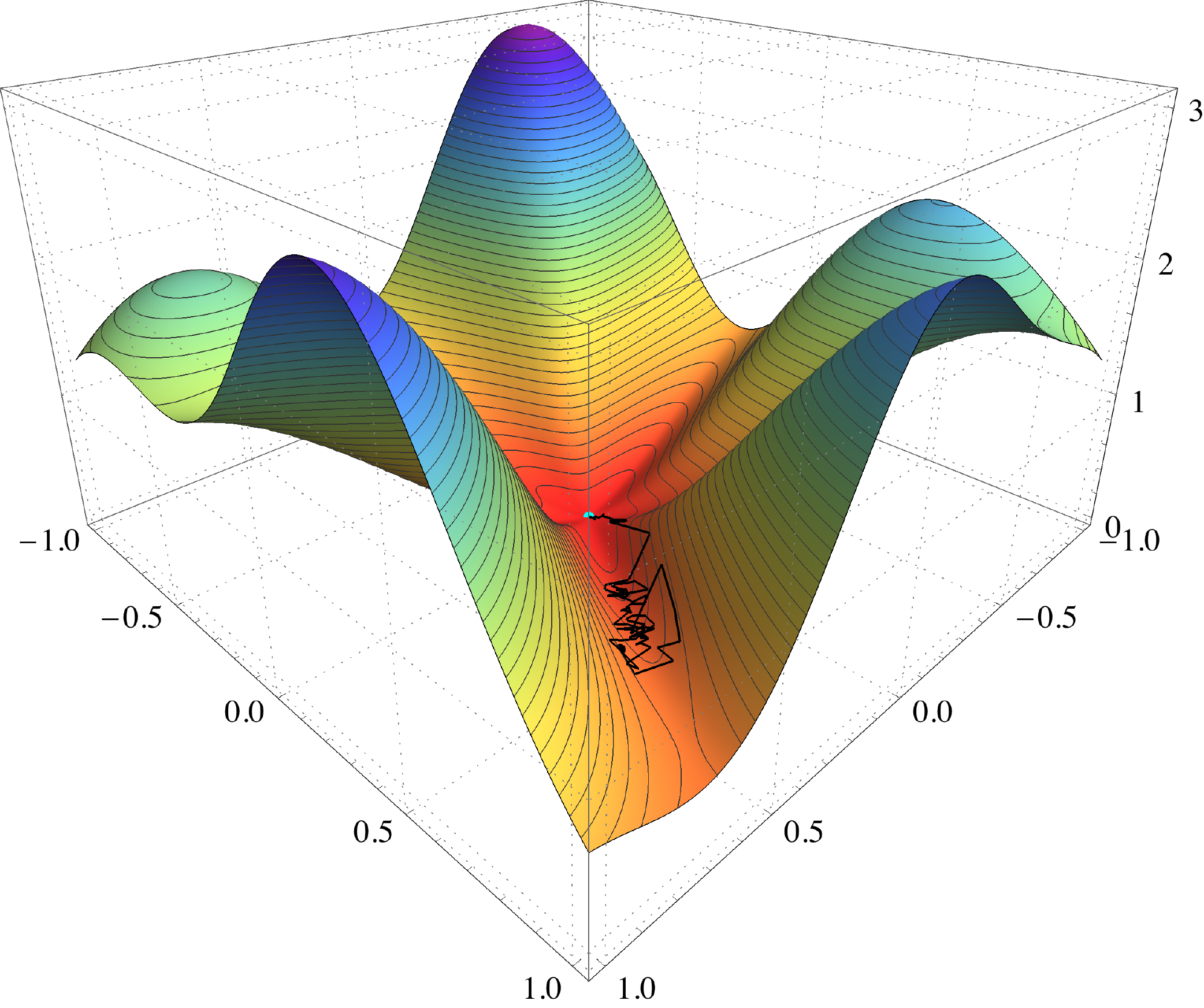}}%
	\caption{Examples of variationally coherent objectives:
		on the top row, the Beale function $\obj(x_{1},x_{2}) = (1.5 - x_{1} + x_{1}x_{2})^{2} + (2.25 - x_{1} + x_{1}x_{2}^{2})^{2} + (2.625 - x_{1} + x_{1}x_{2}^{3})^{2}$ over the benchmark domain $[-4,4]\times[-4,4]$;
		on the bottom row, the polar example $\obj(r,\theta) = (3 + \sin(5\theta) + \cos(3\theta)) r^{2} (5/3 - r)$ over the unit ball $0\leq r\leq1$.
		In both figures, the black curves indicate a sample trajectory of \ac{DASGD} with linearly growing delays.}
	\label{fig:plot-global}
\end{figure}


\begin{assumption}
	\label{asm:VC}\quad
	The optimization problem \eqref{eq:opt} is called \emph{variationally coherent in the mean} if
	\begin{equation}
	\label{eq:VC}
	\tag{VC}
	\exof{\braket{\grad\sobj(x;\sample)}{x-\sol}}
	\ge 0,
	\end{equation}
	for all $\sol\in\solset$ and all $\x \in \feas$ with equality only if $x\in\solset$.
\end{assumption} 

Note that we do not impose the ``if" condition for equality: if $x \in \solset$, then the equality may or may not hold.
By \cref{asm:basic}, we can interchange expectation and differentiation in \eqref{eq:VC} to obtain
\begin{equation}
\braket{\nabla\obj(x)}{x-\sol}
\ge 0,
\end{equation}
for all $x\in\feas$, $\sol\in\solset$.
As a result, mean variational coherence can be interpreted as an averaged coherence condition for the deterministic optimization problem with objective $\obj(x)$. Mean variationally coherent optimization problems include convex programs, pseudo-convex programs, non-degenerate quasi-convex programs and star-convex programs as special cases.  For instance, if $g$ is star-convex, then $\braket{\nabla \obj(x)}{x - \sol} \geq \obj(x) - \obj(\sol)$ for all $x\in\feas$, $\sol\in\solset$.
This is easily seen to be a special case of variational coherence because
$\braket{\nabla \obj(x)}{x - x^*} \ge \obj(x) - \obj(x^*) \ge 0$, with the last inequality strict unless $x \in \mathcal{X}^*$. Note that star-convex functions contain convex functions as a subclass (but not necessarily pseudo/quasi-convex functions). See~\cite{zhou2017stochastic, zhou2020} for a more detailed discussion on why these various classes of optimization problems are special cases.
\cref{fig:plot-global} also provides two more elaborate examples of a variationally coherent optimization problem that are not quasi-convex.
Put together, these examples indicate that variationally coherent objectives can have highly non-convex profiles. 
Nevertheless, for the class of variationally coherent functions, it is possible to establish almost sure convergence guarantees (which we do so in the subsequent sections).

\subsection{Deterministic Analysis: Convergence to Global Optima}
\label{sec:deterministic}

To streamline our presentation and build intuition along the way, we will begin with the deterministic case in this subsection, where there is no randomness in the calculation of a gradient update. 
In this case, \ac{DASGD} boils down to \acdef{DAGD}, as illustrated in \cref{alg:DAGD}:

\begin{algorithm}[H]
	\caption{Distributed Asynchronous Gradient Descent (\ac{DAGD})} 
	\label{alg:DAGD}
	\begin{algorithmic}[1]
		\REQUIRE Initial state $\scored_{0}\in \R^{\vdim}$, step-size sequence $\step_{\run}$
		\STATE $\run\leftarrow0$
		\REPEAT
		\STATE $\actd_{\run} = \mirror(\scored_{\run})$;
		\STATE $\scored_{\run+1} = \scored_{\run} - \step_{\run+1} \nabla \obj(\actd_{s(n)})$;
		\STATE $\run\leftarrow \run+1$;
		\UNTIL{end}
		\STATE \textbf{return} solution candidate $\actd_{\run}$
	\end{algorithmic}
\end{algorithm}

\subsubsection{Energy Function}
\label{subsec:energy}

We start with an energy function that measures on how ``optimal" the dual variable $\scored$ is: the smaller the energy (i.e. the closer it is to $0$), the better the dual variable. 

\begin{definition}
	\label{def:energy}\quad
	Let $\sol \in \solset$. Define the energy function $\energy\from\R^{\vdim} \to \R$ as follows:
	\begin{equation}
	\energy(\scored)
	= \inf_{\sol \in \solset} \energy_{\sol}(\scored), \text{ where }
	\energy_{\sol}(\scored)
	= \norm{\sol}_{2}^{2}
	- \norm{\mirror(\scored)}_{2}^{2}
	+  2\braket{\scored}{\mirror(\scored) - \sol}.
	\end{equation}
\end{definition}

Note that one can think of $\energy_{\sol}(\scored)$ as the energy of $\scored$ with respect to a fixed optimal solution $\sol$, while $\energy(\scored)$ is the best (smallest) energy for a given $\scored$. We next characterize a few of its useful properties.
\begin{lemma}\quad
	\label{lem:energy}
	For all $\scored\in\R^{\vdim}$, we have:
	\vspace{2ex}
	\begin{enumerate}
		\addtolength{\itemsep}{1ex}
		\item $\energy(\scored) \ge 0$ with equality if and only if $\mirror(\scored) \in \solset$.
		\item Let $\{\scored_{\run}\}_{\run=1}^{\infty}$  be a given sequence.
		If $\lim_{\run\to\infty} \energy(\scored_{\run}) = 0$, then
		$\mirror(\scored_{\run}) \rightarrow \solset$ as\footnote{Following the convention in point-set topology, a sequence $s_{\run}$ converges to a set $\mathcal{S}$ if $\text{dist}(s_{\run}, \mathcal{S}) \rightarrow 0$, with $\text{dist}(\cdot, \cdot)$ being the standard point-to-set distance function: $\text{dist}(s_{\run}, \mathcal{S}) \triangleq \inf_{s \in \mathcal{S}} \text{dist}(s_{\run}, s)$, where $\text{dist}(s_{\run}, s)$ is the Euclidean distance between $s_{\run}$ and $s$.} $\run\to\infty$.
	\end{enumerate}
\end{lemma}

\begin{remark}\label{rem:prop1}\quad
	The proof is given in the appendix, but it is helpful to make a few quick remarks.
	The first statement justifies the terminology of ``energy", as $\energy(\scored)$ is always non-negative.
	This energy function will also be the tool we use to establish an important component of the global convergence result.
	Further, given that $\energy(\scored) \ge 0$, it should also be clear that
	when $\mirror(\scored) \in \solset \subset \feas$, we can choose a particular $\sol = \mirror(\scored)$ such that $\energy(\scored) = 0$ (but that $\energy(\scored) = 0$ implies $\mirror(\scored) \in \solset$ is less obvious). 
    Statement 2 of the lemma provides us with a way to establish convergence to optimal solutions.
	If we can show that $\energy(\scored_{\run}) \to 0$, then $\actd_{\run} = \mirror(\scored_{\run}) \to \solset$.
	Nevertheless, as we shall see later, unlike many other Lyapunov functions in optimization, $\energy(\scored_{\run})$ does not decrease monotonically; consequently, it is difficult to directly establish $\energy(\scored_{\run}) \to 0$. In fact, a finer-grained analysis is required to characterize the convergence behavior of $\energy(\scored_{\run})$ (see \cref{subsubsec:main_deter}).
\end{remark}

We also assume that the converse of Statement 2 of the above lemma holds:

\begin{assumption}\label{assump:reci}\quad
	If $\lim_{\run\to\infty} \energy(\scored_{\run}) = 0$, then
	$\mirror(\scored_{\run}) \rightarrow \solset$ as $n \rightarrow \infty$.
\end{assumption} 

\Cref{assump:reci} can be seen as a primal-dual analogue
of the reciprocity conditions for the Bregman divergence. This assumption usually holds (e.g. when the feasible set $\feas$ is a polytope), unless  $\feas$ is pathological. 

\begin{lemma}\label{lem:energy_prop}
	\quad Fix any $\sol \in \solset$.
	\vspace{2ex}
	\begin{enumerate}
		\addtolength{\itemsep}{1ex}
		\item 
		$\|\mirror(\scored) - \hat{\scored}\|_2^2 - \|\mirror(\hat{\scored}) -\hat{\scored}\|_2^2 \le
		\|\scored-\hat{\scored}\|_2^2$, for any $\scored, \hat{\scored} \in \mathbb{R}^d$.
		\item 
		$\energy_{\sol}(\scored + \Delta \scored) - \energy_{\sol}(\scored) \le 2\langle   \Delta \scored, \mirror(\scored) - \actd^*  \rangle +
		\|\Delta \scored\|_2^2$, for any $\scored,\Delta \scored \in \mathbb{R}^d$.
	\end{enumerate}
\end{lemma}

\begin{remark}\label{rem:prop2}
	\quad
	The first statement of the above lemma serves as an important intermediate step in proving the second statement, 
	and is established by leveraging the envelop theorem and several important properties of Euclidean projection.
	To see that this is not trivial, consider the quantity $\|\mirror(\scored) - \hat{\scored}\|_2 - \|\mirror(\hat{\scored}) -\hat{\scored}\|_2$,
	which we know by triangle's inequality satisfies:
	\begin{equation}\label{eq:easy}
	\begin{split}
	&\|\mirror(\scored) - \hat{\scored}\|_2 - \|\mirror(\hat{\scored}) -\hat{\scored}\|_2 \le \| \mirror(\scored) -\mirror(\hat{\scored}) \|_2 \le \|\scored - \hat{\scored}\|_2,
	\end{split}
	\end{equation}
	where the last inequality follows from the fact that projection is a non-expansive map.
	However, this inequality is not sufficient for our purposes because in quantifying the perturbation $\energy(\scored + \Delta \scored) - \energy(\scored)$, we also need the squared term  $\|\Delta \scored\|_2^2$, which is not easily obtainable from Equation~\eqref{eq:easy}.
	In fact, a tighter analysis is needed to establish that $\|\scored -\hat{\scored}\|_2^2$ is an upper bound on $\|\mirror(\scored) - \hat{\scored}\|_2^2 - \|\mirror(\hat{\scored}) -\hat{\scored}\|_2^2$. 
\end{remark}

\subsubsection{Main Convergence Result}\label{subsubsec:main_deter}

An intermediate result, interesting on its own and useful also as a heavy-lifting tool for our convergence analysis is then provided by the following technical result:

\begin{proposition}\quad
	\label{prop:recur_d}
	Under \cref{asm:basic,asm:VC,asm:delays,assump:reci}, \ac{DAGD} admits a subsequence $\actd_{\run_{k}}$ that converges to $\solset$ as $k\to\infty$.
\end{proposition}

We highlight the main steps below and refer the reader to the appendix for the details:
\begin{enumerate}
	\item 
	Letting $b_{\run} = \v(x_{s(n)}) - \v(x_\run)$,
	we can rewrite the gradient update in \ac{DAGD} as:
	\begin{flalign}
	y_{\run+1}
	&= y_{\run}
	- \step_{\run+1} \v(x_{s(n)})
	\notag\\
	&= y_{\run}
	- \step_{\run+1} \v(x_\run)
	- \step_{\run+1} \{\v(x_{s(n)}) - \v(x_\run)\}
	\notag\\
	&= y_{\run}
	- \step_{\run+1} (\v(x_\run) + b_{\run}).
	\end{flalign}
	Recall here that  $s(\run)$ denotes the previous iteration count whose gradient becomes available only at the current iteration $\run$.
	By bounding the magnitude of $b_{\run}$ using the delay sequence through a careful analysis, we establish that under any one of the conditions in \cref{asm:delays}, $\lim_{\run \to \infty} \|b_{\run}\|_2 = 0$.
	The analysis here, particularly the one for the last three conditions, reveals the following pattern: as the magnitude of the delays gets larger and larger in the order of growth, one needs to use a more conservative step-size sequence in order to mitigate the damage done by the stale gradient information. Intuitively, smaller step-sizes are more helpful in larger delays because they carry a better ``amortization" effect that makes \ac{DAGD} more tolerant to delays.

	\item 
	With the defintion of $b_{\run}$, \ac{DAGD} can be written as:
	\begin{equation}
	\begin{aligned}
	\actd_{\run}
	&= \mirror(\scored_{\run}),\\
	\scored_{\run+1}
	&= y_{\run} - \step_{\run+1} (\v(x_\run) + b_{\run}).
	\end{aligned}
	\end{equation}
	We then use the energy function to study the behavior of $\scored_{\run}$ and $\actd_{\run}$. More specifically, we look at the quantity $\energy(\scored_{\run+1}) - \energy(\scored_{\run})$ and, using \cref{lem:energy_prop}, we bound this one-step change using the step size $\step_{\run}$, the $b_{\run}$ sequence and 
	the defining quantity $ \langle \v(x_\run), \actd_{\run} - \actd^*  \rangle$ of a variationally coherent function (as well as another term that will prove inconsequential).  
	We then telescope on $\energy(\scored_{\run+1}) - \energy(\scored_{\run})$ to
	obtain an upper bound for $\energy(\scored_{\run+1}) - \energy(\scored_0)$.
	Since the energy function is always non-negative, 
	$\energy(\scored_{\run+1}) - \energy(\scored_0)$ is at least $- \energy(\scored_0)$ 
	for every $\run$. Then, utilizing the fact that $b_{\run}$ converges to $0$ and that 
	$ \langle \v(x_\run), \actd_{\run} - \actd^*  \rangle$ is always positive (unless the iterate is exactly an optimal solution, in which case it is $0$), we show that the upper bound will approach $-\infty$ if $\act_{\run}$ only enters $\N(\solset, \epsilon)$, an open $\epsilon$-neighborhood of $\solset$, a finite number of times (for an arbitrary $\epsilon > 0$). This generates an immediate contradiction, and thereby establishes that $\act_{\run}$ will get arbitrarily close to $\solset$ for an infinite number of times.
	This then implies that there exists a subsequence of \ac{DAGD} iterates that converges to the solution set of \eqref{eq:opt}, i.e, $\actd_{\run_k} \to \solset$ as $k \to \infty$.
\end{enumerate}

\begin{theorem}
	\label{thm:convergence}\quad
	Under \cref{asm:basic,asm:VC,asm:delays,assump:reci}, the global state variable $\actd_{\run}$ of \ac{DAGD} \textpar{\cref{alg:DAGD}} converges to the solution set $\solset$ of \eqref{eq:opt}.
\end{theorem}

We give an outline of the proof below, referring to the appendix for the details.

Fix a $\delta > 0$.
Since $\actd_{\run_k} \to \actd^*$, as $k \to \infty$ (per \cref{prop:recur_d}),
we have $\energy(\scored_{\run_k}) \to 0$ as $k \to \infty$ per \cref{lem:energy}.
So we can pick an $\run$ that is sufficiently large and $\energy(\scored_{\run}) < \delta$.
Our goal here is to show that for $n$ large enough, once $\energy(\scored_{\run}) < \delta$, it will stay so forever: $\energy(\scored_{m}) < \delta, \forall m \ge n$.
Note in particular that the above statement is not true for all $\run$, but only for $\run$ large enough.

However, the behavior of $\energy(\scored_{\run})$ is not very regular: it can certainly increase from iteration to iteration for any $\run$. Nevertheless, we can precisely quantify how large this increment (if any) can be. This leads us to
break it down to two cases:
\vspace{2ex}
\begin{enumerate}
	\addtolength{\itemsep}{1ex}
	\item 
	Case 1: $\energy(\scored_{\run}) < \delta/2$.
	
	\item
	Case 2: $\delta/2 \leq \energy(\scored_{\run}) < \delta$.
\end{enumerate}

For Case 1, we show in the appendix that 
\begin{equation}
\energy{(\scored_{\run+1})} - \energy(\scored_{\run}) 
\leq 2BC_4\step_{\run+1}  +  2\step_{\run+1}^2(C_2 + B^2),
\end{equation}
for suitable constants $B$ and $C$‘s.
Now for $\run$ sufficiently large, we can make the right-hand arbitrarily small, and in particular, smaller
than $\delta/2$.
This means $\energy{(\scored_{\run+1})} \le \energy{(\scored_{\run})} + \frac{\delta}{2} < \delta$. Consequently, in this case, the energy stays within the $\delta$ bound in the next iteration.

For Case 2, we show in the appendix that
\begin{equation}
\energy{(\scored_{\run+1})} - \energy(\scored_{\run})
\leq -2\step_{\run+1} \bracks*{\frac{a}{2}  - \step_{\run+1}(C_2 + B^2)},
\end{equation}
where $a$ is a positive constant that depends only on $\delta$.
Again, since $\run$ is sufficiently large, we can make $\frac{a}{2}  - \step_{\run+1}(C_2 + B^2)$ positive, thereby making the right-hand side negative.
Consequently, $\energy{(\scored_{\run+1})} < \energy(\scored_{\run}) < \delta$.
Hence, again, the energy stays within the $\delta$ bound in the next iteration.

The key conclusion from the above is that, for large enough $\run$, once 
$\energy{(\scored_{\run})}$ is less than $\delta$,
$\energy(\scored_{\run+1})$ is less than $\delta$ as well and so are all the iterates afterwards. Since $\delta$ is arbitrary, it follows $\energy{(\scored_{\run})} \to 0$, and therefore $\actd_\run \to \solset$ by \cref{lem:energy}.

\subsection{Stochastic Analysis: Almost Sure Convergence to Global Optima}
\label{sec:stochastic}

Having established deterministic global convergence of \ac{DAGD}, we now proceed to study stochastic global convergence of \ac{DASGD}. Compared to the deterministic analysis,
the stochastic case is much more involved because randomness can lead to very volatile behavior in the presence of delays; in particular, the simple approach employed in \cref{thm:convergence} (to establish that once  $\energy{(\scored_{\run})}$ is less than $\delta$, it will always remain so) no longer works. To deal with both delays and noise, a much more sophisticated analysis framework needs to be developed, which requires several news ideas.
To streamline our presentation, we break the theoretical development into four subsections, each comprising an important component and step of the overall analysis.

\subsubsection{Recurrence of \ac{DASGD}}

Our first step lies in generalizing \cref{prop:recur_d} to the stochastic case. Specifically, in the presence of noise, we show that the iterates of \ac{DASGD} 
visit any neighorhood of $\solset$ infinitely often almost surely.

\begin{proposition}
	\label{prop:recur_s}\quad
	Under \cref{asm:basic,asm:VC,asm:delays,assump:reci}, \ac{DASGD} admits a subsequence $\act_{\run_{k}}$ that converges to $\solset$ almost surely: $\act_{\run_k} \to \solset$ with probability $1$ as $k \to \infty$.
\end{proposition}

	We outline the two main steps of the proof below, referring the reader to the appendix for the details.
	
	\begin{enumerate}
		\item 
		We begin by rewriting the gradient update step in \ac{DASGD} as:
		\begin{flalign}
		\frac{\score_{\run+1} - \score_{\run}}{\step_{\run+1}}
		&= -\nabla \sobj(\act_{s(\run)},\sample_{\run+1})
		\notag\\
		&= -\v(\act_\run)
		\notag\\
		&-\bracks{\v(\act_{s(\run)})- \v(\act_\run)}
		\notag\\
		&-\bracks{\nabla \sobj(\act_{s(\run)},\sample_{\run+1}) -  \v(\act_{s(\run)})}.
		\end{flalign}
		Letting $B_{\run} = \v(\act_{s(\run)})- \v(\act_\run)$ and $U_{\run+1} = \nabla \sobj(\act_{s(\run)},\sample_{\run+1}) -  \v(\act_{s(\run)})$, we can rewrite the \ac{DASGD} update as
		\begin{equation}
		\label{eq:DASGD_abbre}
		\score_{\run+1}
		= \score_{\run} - \step_{\run+1}\{\v(\act_\run) + B_{\run} + U_{\run+1}\}.
		\end{equation} 
		
		We then establish the following two facts in this step. 
		First, we verify that 
		$\sum_{r=0}^{\run} U_{\run+1}$ is a martingale adapted to $\score_1, \score_2 \dots, \score_{\run+1}$, where $\{U_{\run+1}\}_{\run=0}^{\infty}$ is a $L_2$-bounded martingale difference sequence. Second, we show that $\lim_{\run \to \infty} \|B_{\run}\|_2 = 0, { a.s.}.$
		
		The second claim is done by first giving an upper bound on $\|B_{\run}\|_2$ by writing $\v(\act_{s(\run)})- \v(\act_\run)$ as a sum of one-step changes ($\v(\act_{s(\run)})- \v(\act_{s(\run)+1}) + \v(\act_{s(\run)+1}) - \dots + \v(\act_{\run-1})- \v(\act_\run)$) and analyzing each such successive change. 
		We then break that upper bound into two parts, one deterministic and one stochastic. For the deterministic part, the same analysis in the proof of \cref{prop:recur_d} yields convergence to $0$.
		
		The stochastic part turns out to be the tail of a martingale.
		By leveraging the property of the step-size and a crucial property of martingale differences (two martingale differences at different time steps are uncorrelated), we establish that said martingale is $L_2$-bounded.
		Then, by applying a version of Doob's martingale convergence theorem, it follows that said martingale converges almost surely to a limit random variable with finite second moment (and hence almost surely finite).
		Consequently, writing the tail as a difference between two terms (each of which converges to the same limit variablewith probability $1$), we conclude that the tail converges to $0$ \as.

		
		\item 
		The full \ac{DASGD} update may then be written as
		\begin{flalign}
		\act_{\run}
		&= \mirror(\score_{\run})
		\notag\\
		\score_{\run+1}
		&= \score_{\run} - \step_{\run+1} \bracks{\v(\act_\run) + B_{\run} + U_{\run+1}}.
		\end{flalign}
		As in Step 2 of the proof of \cref{prop:recur_d}, we again bound the one-step change of the energy function 
		$\energy(\score_{\run+1}) - \energy(\score_{\run})$ and then telescope the differences.
		The two distinctions from the determinstic case are: 1) Everything is now a random variable. 2) We have three terms: in addition to the random gradient $\v(\act_\run)$ and the random drift
		$B_{\run}$, we also have a martingale term 
		$U_{\run+1}$. Since $B_{\run}$ converges to $0$ almost surely (as shown in the previous step),
		its effect can be shown to be negligible.
		Futher, the analysis utilizes law of large numbers for martingale as well as Doob's martingale convergence theorem to bound the effect of the various martingale terms and to establish that the final dominating term converges to $-\infty$ almost surely (which generates a contradiction since the energy function is always positive) unless a subsequence $\act_{\run_k}$ converges almost surely to $\solset$.
		\hfill$\blacksquare$
	\end{enumerate}	 

Even though recurrence, which can be seen as the counterpart of \cref{prop:recur_d}, holds as per the above proposition, the random iterates $\act_{\run}$ are much more irregular than their deterministic counterpart $\actd_{\run}$ in \ac{DAGD}. To deal with this complexity, we work with and characterize the sample trajectories ``generated" (to be made precise later) by $\act_{\run}$ (rather than individual iterates $\act_{\run}$).
To work towards this general direction, we first push the \ac{DASGD} update into a determinstic ordineary differential equation (ODE), as explained in the next subsection.

\subsubsection{Mean-Field Approximation of \ac{DASGD}}

We can rewrite the \ac{DASGD} update as:
\begin{flalign}
\label{eq:help}
\act_{\run}
&= \mirror(\score_{\run})
\notag\\
\score_{\run+1}
&= \score_{\run} - \step_{\run+1}\{\v(\act_\run) + B_{\run} + U_{\run+1}\}.
\end{flalign}
Written in this way, \ac{DASGD} can be viewed as a discretization of the ``mean-field'' ODE
\begin{flalign}
\label{eq:ode}
\actd
&= \mirror(\scored),
\notag\\
\dot\scored
&= -\v(\actd).
\end{flalign}
The intuition is that this ODE provides a ``mean" approximation of the \ac{DASGD} update, because in \eqref{eq:help},
the noise term $U_{\run+1}$ has $0$ mean, and the term $B_{\run}$ converges to $0$ (and therefore has
negilible effect in the long run).
Thes leaves only the term $\score_{\run+1} = \score_{\run} - \step_{\run+1}\v(\act_\run)$, which can be seen as a Euler discretization of the ODE.
(Of course, that Equation~\eqref{eq:ode} is a good-enough approximation of Equation~\eqref{eq:help} for global almost sure convergence purposes here will be rigorously justified later.) 

Next, writing Equation~\eqref{eq:ode} solely in terms of $\scored$ yields $\dot{\scored}= -\v(\mirror(\scored))$.
Since $\v$ and $\mirror$ are both Lipschitz continuous and $\feas$ is a compact set, the composition $\v\circ\mirror$ is itself Lipschitz continuous and bounded.
Standard results from the theory of dynamical systems then show that~\eqref{eq:ode} admits a unique global solution $\scored(t)$ for any initial condition $\scored(0)$.
On the other hand, since $\mirror$ is not a one-to-one map, it is not invertible;
consequently, there need not exist a unique solution trajectory for $\actd(t)$.
By this token, the rest of our analysis will focus on the trajetory of $\scored(t)$.

With the guarantee of the existence and uniqueness of the $\scored$ trajectory, let $P\from\R_{+}\times\R^{\vdim}\to\R^{\vdim}$ be the semiflow\footnote{See \cref{app:first} for a more rigorous definition. Furthermore, $\R_+$ is the set of non-negative reals.} of \eqref{eq:ode}, i.e., $P(t, y_0)$ denotes the state of \eqref{eq:ode} at time $t$ when the initial condition is $y_0$.
In other words, when viewed as a function of time, $P(\cdot, y_0)$ is the solution trajectory to $\dot{\scored}= -\v(\mirror(\scored))$.
It is worth pointing out that writing it in this double-argument form also allows us to interpret $P$ as a function of the initial condition: for a fixed $t$, $P(t, \cdot)$ gives different states at $t$ when the ODE starts from different initial conditions (in particular, $P(0, \scored) = \scored$). Both views will be useful later.

It turns out that with the energy function introduced here, 
$\energy(\flow(t, \scored))$ is always non-increasing. Furthermore, it is also decreasing at a meaingful rate.
We end this subsection with a ``sufficient decrease'' property of the mean dynamics \eqref{eq:ode} (the proof given in the appendix due to space limitation):
\begin{lemma}
	\label{lem:odeprop}\quad
	With notation as above, we have:
	\vspace{2ex}
	\begin{enumerate}
		\addtolength{\itemsep}{1ex}
		\item
		If $\mirror(\flow(t, \scored)) \notin \solset$,  then $\energy(\flow(t, \scored))$ is strictly decreasing in $t$ for all $y\in\R^{\vdim}$.
		\item
		For all $\delta > 0$, there exists some $T\equiv T(\delta) > 0$ such that, for all $t\geq T$, we have
		\begin{equation}
		\hspace{-2em}
		\sup\nolimits_{\scored}
		\setdef{\energy(\flow(t, \scored)) - \energy(\scored)}{\energy(\flow(t,\scored)) \ge  \delta/2}
		\le -\delta/2.
		\end{equation}
	\end{enumerate}
	\end{lemma}	

\cref{lem:odeprop} essentially says $\energy(\flow(t, \scored))$ is strictly and uniformly (across all $\scored$) decreasing at a non-vanishing rate. 
More specifically, to give some intuition of the second part of \cref{lem:odeprop}, note that $\energy(\flow(t, \scored)) - \energy(\scored)$
is the energy change at time $t$ when starting at $\scored$.
The first part says this quantity is always non-negative. While the second part
says provided $\energy(\flow(t,\scored)) \ge  \delta/2$,
the decrease in energy will be at least $\frac{\delta}{2}$, no matter what the initial point $\scored$ is. If $\energy(\flow(t,\scored)) \ge  \delta/2$ does not hold,
that means the energy at time $t$ is already really small (i.e. $\energy(\flow(t,\scored)) <  \delta/2$). Consequently, the mantra of the above lemma can be stated succinctly as follows: either the energy is already close to $0$, or the energy will decrease towards $0$. 

In fact, by some additional analysis, one can further show\footnote{Although this is an interesting conclusion, we do not prove it here because we are mainly concerned with establishing convergence of the \ac{DASGD} iterates, rather than the ODE solution trajectory.} that \cref{lem:odeprop} implies $  \flow(t, \scored) \to \solset, \forall y$ as $t \to \infty$. 
Now, despite being an interesting result on its own (which establishes that the continuous dynamics of \ac{DASGD} converges to $\solset$), it is still 
some distance away from our final desideratum: our goal is to establish (almost sure) convergence of the discrete-time process in \ac{DASGD}.
So unless we can somehow relate the discrete-time iterates to the continuous-time trajectories, the convergence of \ac{DASGD} is still uncertain.
We fulfill this taks in the next subsection.

\subsubsection{Relating \ac{DASGD} Iterates to ODE Trajectories}

Our goal here is to establish a quantitative connection between the 
\ac{DASGD} iterates and the ODE trajectory studied in the previous subsection.
Our general idea is that if we show the trajectory generated by the discrete-time iterates of \ac{DASGD} is \textbf{path-by-path} ``close" to the continuous-time trajectory $\flow(t, \scored)$, then likely almost-sure convergence of the \ac{DASGD} iterates can be guaranteed as well.

To be more specific, there are two things that need to be more precisely defined from the preceding high-level discussion.
First, what does it mean to be a trajectory generated by the discrete-time iterates of \ac{DASGD}?
Second, what does it mean for this trajectory to be ``close" to the ODE trajectory?
In general, the answers to these questions can vary depend on the specific goal at hand.
In the current context, our goal is to establish global almost sure convergence (a very strong result).
Consequently, we need to choose the answers rather judiciously: on the one hand, the answers must be stringent enough to ensure
global almost sure convergence in the end (for instance, for the second question, a fairly strong notion of ``closeness" is needed); on the other hand, they must also be flexible enough to fit in the current context.

As it turns out, the answer to the first question is rather intuitive: (perhaps) the simplest way to generate a continuous trajectory from a sequence of discrete points is the \textit{affine interpolation}: connect the iterates $\score_0, \score_1,\dots, \score_{\run}$ at times $0, \alpha_1, \dots, \sum_{r=1}^{\run=1}\alpha_r$.
We call this curve the affine interpolation curve of \ac{DASGD} and denote it by $A(t)$. 
Note that $A(t)$ is a random curve because the \ac{DASGD} iterates $\score_0, \score_1,\dots, \score_{\run}$ are random. To avoid confusion, we summarize the three different objects discussed so far:
\vspace{2ex}
\begin{enumerate}
	\addtolength{\itemsep}{1ex}
	\item The \ac{DASGD} iterates $\score_0, \score_1,\dots, \score_{\run}$.
	\item The affine interpolation curve $A(t)$ of $\score_{\run}$.
	\item The  flow $P(t, \scored)$ of the ODE \eqref{eq:ode}.
\end{enumerate}

The answer to the second question lies in the notion of an \acdef{APT} , a concept introduced by \citet{BH96} and \citet{BS00}.
Specifically, in our current context, a continuous curve $s(t)$ is considered close
to ODE solution $P(t, \scored)$ if the following holds:

\begin{definition}\quad
	A continuous function $s: \mathbb{R}_+ \to \R^\vdim$ is an \ac{APT} for $\flow$ if for every $T > 0$, 
	\begin{equation}
	\lim_{t \to \infty} \sup_{0 \le h \le T} d(s(t+h), \flow(h, s(t))) = 0,
	\end{equation} 	
	where $d(\cdot, \cdot)$ is the Euclidean metric\footnote{As should be obvious from the definition, APTs can be defined more generally in metric spaces in exactly the same way.} in $\R^\vdim$.
\end{definition}

Intuitively, the definition matches exactly the naming:
$s$ is an \ac{APT} for $\flow$ if, for sufficiently large $t$, the flow lines of $P$ remain arbitrarily close to $s(t)$ over a time window of any (fixed) length.
More precisely, for each fixed $T > 0$, one can find a large enough $t_0$, such that for all $t > t_0$,
the curve $s(t+h)$ approximates the trajectory $\flow(h, s(t))$ on the interval $h \in [0, T]$ with any predetermined degree of accuracy. 

With this definition in place, to push through the agenda, we need to establish that 
$A(t)$, the affine interpolation curve of the \ac{DASGD} iterates, is an APT for the flow $P(t, \scored)$ induced by the ODE~\eqref{eq:ode}. 
More precisely, we establish that $A(t)$ is an APT for the flow $P(t, \scored)$ almost surely, because as mentioned before, $A(t)$ is a random curve.

\begin{lemma}\label{lem:apt}\quad
Let	$A(t)$ be the random affine interpolation curve generated from the \ac{DASGD}
iterates.
Then $A(t)$ is an \ac{APT} of $P(t,\scored)$ almost surely.
\end{lemma}

Note that this result means any resulting affine interpolation curve of \ac{DASGD} is close to the ODE trajectory. This also forms the basis for reasoning convergence on a path-by-path scale.
However, more work still remains to be done 
because, unfortunately, the condition that $A(t)$ is an \ac{APT} for $P$ almost surely does not guarantee that the discrete-time iterates of \ac{DASGD} converge to $\solset$ (for many counterexamples in general dynamical systems, see \citet{Ben99}). 
In other words, the notion of \ac{APT} is not sharp enough to ensure direct convergence result.
We fulfill this final gap in the next subsection.

\subsection{Main Convergence Result}

Even though $A(t)$ being an \ac{APT} for $P$ almost surely does not itself guarantee that the discrete-time iterates of \ac{DASGD} converge to $\solset$,
we can use the energy function to further sharpen this result. Specifically, we use $E(A(t))$ to further control the behavior 
of the affine interpolation curve.
In fact, the advantage of working with the affine interpolation curve $A(t)$ is that 
once we show $\energy(A(t))$ is bounded by some $\delta$ almost surely from some point onwards, then we know
$\energy(\score_{\run})$ is bounded by $\delta$ almost surely (also from some point onwards): this is because the 
discrete-time iterates and the affine curve coincide at discrete time points.
Consequently, we focus on boudning $\energy(A(t))$, which will enable us to
 obtain our main convergence result:

\begin{theorem}
	\label{thm:main}\quad
	Under \cref{asm:basic,asm:VC,asm:delays,assump:reci}, the global state variable $\act_{\run}$ of \ac{DASGD} \textpar{\cref{alg:DASGD}} converges \as to the solution set $\solset$ of \eqref{eq:opt}.
\end{theorem}

Again, we only give an outline of the proof below.
	By \cref{prop:recur_s}, $\act_{\run}$ gets arbitrarily close to $\solset$ infinitely often.	
	Thus, it suffices to show that, if $\act_{\run}$ ever gets $\epsilon$-close to $\solset$, all the ensuing iterates are $\epsilon$-close to $\solset$ \as.
The way we show this ``trapping" property is to use the energy function. Specifically, we consider $\energy(A(t))$ and show that
	no matter how small $\epsilon$ is, for all sufficiently large $t_0$, if $\energy(A(t_0))$ is less than $\epsilon$ for some $t_0$, then $\energy(A(t)) < \epsilon, \forall t> t_0$.
	This would then complete the proof because $A(t)$ actually contains all the \ac{DASGD} iterates, and hence if 
	$\energy(A(t)) < \epsilon, \forall t> t_0$, then $\energy(\score_{\run}) < \epsilon$ for all sufficiently large $\run$. Furthermore, since $A(t)$ contains all the iterates, the hypothesis  that `` if $\energy(A(t_0))$ is less than $\epsilon$ for some $t_0$" will be satisfied due to \cref{prop:recur_s}.

	We expand on one more layer of detail and defer the rest into appendix.
	To obtain control $\energy(A(t))$, we control two things: the energy on the ODE path $\energy(\flow(t, y))$ and the discrepancy between $\energy(\flow(t, y))$ and $\energy(A(t))$.
	The former can be made arbitrarily small as a result of \cref{lem:odeprop} (we have a direct handle on how the ODE path would behave). The latter can also be made arbitrarily small as a result of \cref{lem:apt}:
	since $A(t)$ is an \ac{APT} for $\flow$, the two paths are close. Therefore, the discrepancy between $\energy(\flow)$ and $\energy(A)$ should also be vanishingly small.
	Consequently, since $\energy(A(t)) = \energy(\flow(t, y)) +  \{\energy(A(t))- \energy(\flow(t, y))\}$,
	and both terms on the right can be made arbitrarily small, so can $\energy(A(t))$ be made arbitrarily small.

\section{Discussion}
\label{sec:discussion}

\newcommand{\ros}{\obj_{\mathrm{Ros}}}

We end the paper with a short simulation discussion that reveals an interesting practical observation.
Specifically, we test the convergence of \cref{alg:DASGD} against a Rosenbrock test function with $\vdim=1001$ degrees of freedom, a standard non-convex global optimization benchmark \cite{Ros60}.
Specifically, we consider the objective
\begin{equation}
\ros(x)
= \sum_{i=1}^{1000} \bracks{1000(x_{i+1} - x_{i}^{2})^{2} + (1-x_{i})^{2}},
\end{equation}
with $x_{i}\in[0,2]$, $i=1,\dots,1001$.
The global minimum of $\ros$ is located at $(1,\dotsc,1)$, at the end of a very thin and very flat parabolic valley which is notoriously difficult for first-order methods to traverse \cite{Ros60}.
Since the minimum of the Rosenbrock function is known, \eqref{eq:VC} is easily checked over the problem's feasible region.

	\begin{figure*}[tbp]
	\footnotesize
	\subfigure[Convergence with no delays between gradient updates]{\label{fig:conv-synchronous}%
		\includegraphics[width=.482\textwidth]{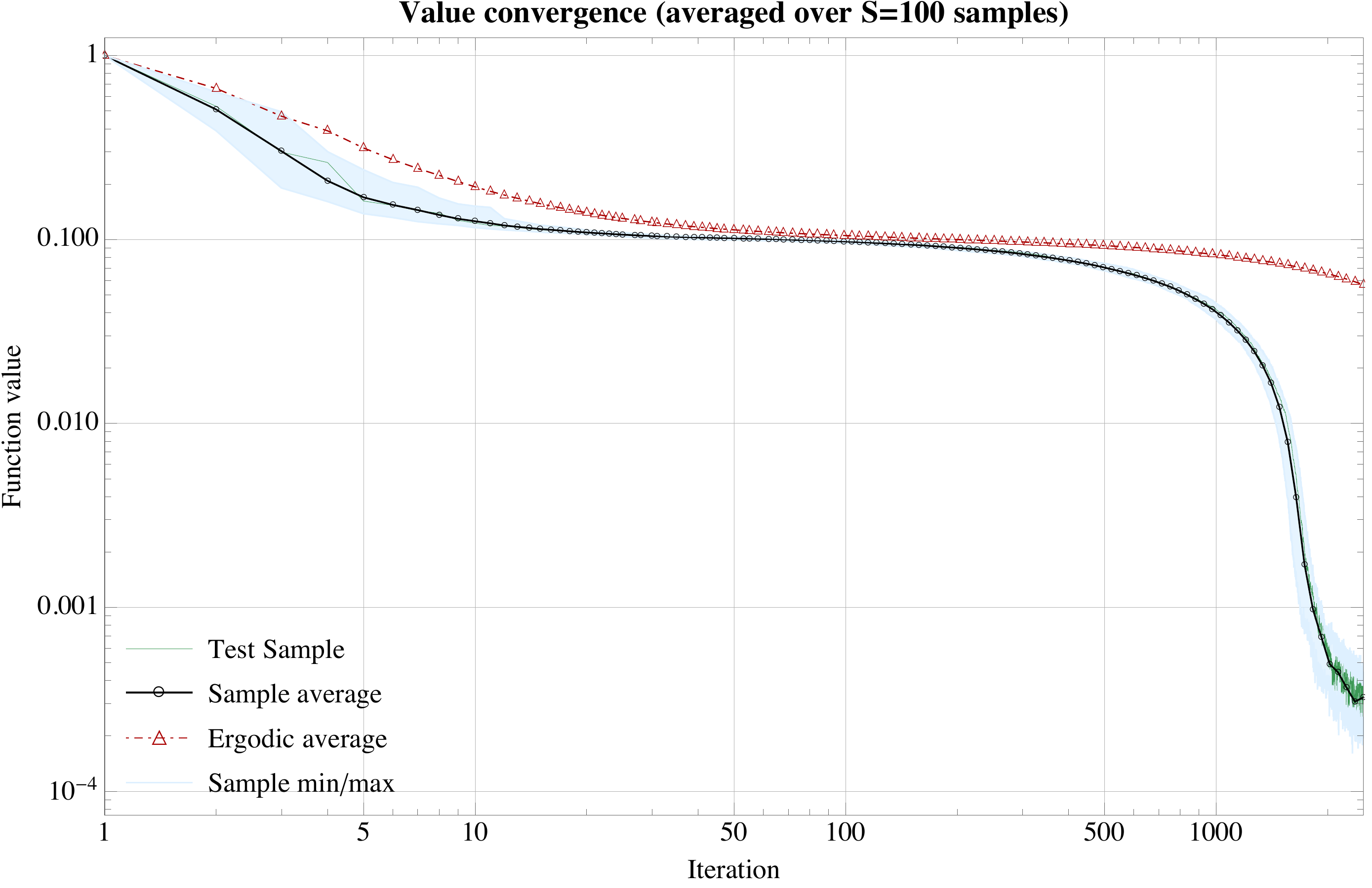}}
	\hfill
	\subfigure[Convergence with linearly growing delays]{\label{fig:conv-asynchronous}%
		\includegraphics[width=.48\textwidth]{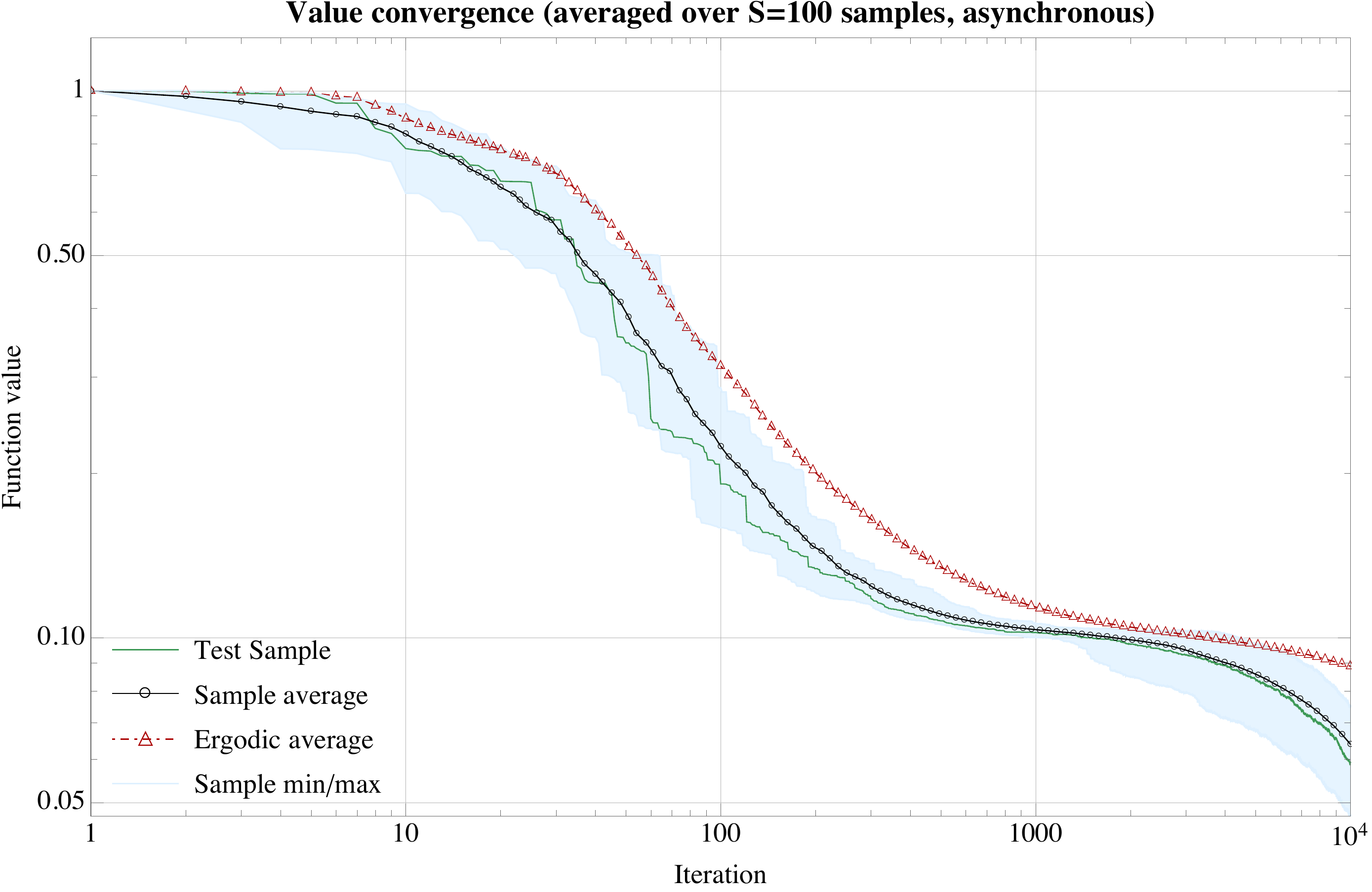}}
	\caption{Value convergence in a non-convex stochastic optimization problem with $\vdim=1001$ degrees of freedom.}
	\label{fig:conv}
\end{figure*}

For our numerical experiments, we considered
\begin{inparaenum}
	[\itshape a\upshape)]
	\item
	a synchronous update schedule as a baseline;
	and
	\item
	an asynchronous master-worker framework with random delays that scale as $d_{\run} = \Theta(\run)$.
\end{inparaenum}
In both cases, \cref{alg:DASGD} was run with a decreasing step-size of the form $\step_{\run} \propto 1/(\run \log\run)$ and stochastic gradients drawn from a standard multivariate Gaussian distribution (i.e., zero mean and identity covariance matrix).

Our results are shown in \cref{fig:conv}.
Starting from a random (but otherwise fixed) initial condition, we ran $S=10^{5}$ realizations of \ac{DASGD} (with and without delays).
We then plotted a randomly chosen trajectory (``test sample'' in \cref{fig:conv}), the sample average, and the min/max over all samples at every update epoch.
For comparison purposes, we also plotted the value of the so-called ``ergodic average''
\begin{equation}
\bar\act_{\run}
= \frac{\sum_{\iRun=1}^{\run} \step_{\iRun} \act_{\iRun}}{\sum_{\iRun=1}^{\run} \step_{\iRun}},
\end{equation}
which is often used in the analysis of \ac{DASGD} in the convex case (see e.g., \citealp{NIPS2011_4247}).
Even though this averaging leads to very robust convergence rate estimates in the convex case, we see here that it performs worse than the worst realization of \ac{DASGD}.
The reason for this is the lack of convexity:
due to the ridges and talwegs of the Rosenbrock function, Jensen's inequality fails dramatically to produce an improvement over $\act_{\run}$ (and, in fact, causes delays as it causes $\act_{\run}$ to deviate from its gradient path).
Consequently, this simple simulation indicates that establishing convergence of the iterate $\act_{\run}$ itself is not only theoretically stronger (and hence more difficult) than convergence of the ergodic average, but also more practically relevant.

\numberwithin{lemma}{section}		
\numberwithin{proposition}{section}		
\numberwithin{corollary}{section}		
\numberwithin{equation}{section}		
\appendix

\section{Auxiliary Results}
\label{app:first}

We collect here in one place all the auxiliary results in the existing literature that will
be used in our proofs subsequently. 
The first one is a well-known characterization of convex sets and the projection operator given in~\cite{Nes04}:

\begin{lemma}\label{lem:aux1}\quad
	Let $\feas$ be a compact and convex subset of $\mathbb{R}^d$. Then for any $\actd \in \feas, \scored \in \mathbb{R}^d$:
	\begin{equation}
	\inner{\mirror(\scored) - \actd}{\mirror(\scored) -\scored}
	\leq 0.
	\end{equation}
	
\end{lemma}

The second one is an $L_{p}$-bounded martingale convergence theorem given in~\cite{HH80}:

\begin{lemma} \quad
	\label{thm:mg_convergence2}
	Let $S_n$ be a martingale adapted to the filtration $\mathcal{S}_n$.
	If for some $p \ge 1$, $\sup_{n \ge 0} \mathbf{E}[|S_n|^p] < \infty$, then $S_n$ converges almost surely to a limiting random variable $S_{\infty}$ with $\mathbf{E}[|S_{\infty}|^p] < \infty$.
\end{lemma}

\begin{remark}\quad
	Note that $\mathbf{E}[|S_{\infty}|^p] < \infty \text{ for } p \ge 1$ obviously implies $S_{\infty}$ is finite almost surely.
\end{remark}

The third one is the classical envelope theorem (see~\cite{carter2001foundations}).

\begin{lemma}\label{lem:env}\quad
	Let $f: \mathbb{R}^n \times \mathbb{R}^m \rightarrow \mathbb{R}$ be a continuously differentiable function.
	Let $U$ be a compact set and consider the problem 
	\begin{equation}\label{eq:env_opt}
	\max_{x \in U} f(x, \theta).
	\end{equation}
	Let $x^*: \mathcal{O} \rightarrow \mathbb{R}^m$ be a continuous function defined on an open set $\mathcal{O} \subset \mathbb{R}^m$ such that for each $\theta \in \mathcal{O}$,
	$x^*(\theta)$ solves the problem in Equation~\ref{eq:env_opt}.
	Define $V: \mathbb{R}^m \rightarrow \mathbb{R}$ where $V(\theta) = f(x^*(\theta), \theta)$.
	Then $V(\theta)$ is differentiable on $\mathcal{O}$ and:
	\begin{equation}
	\grad V(\theta) = \grad f(x^*(\theta), \theta).
	\end{equation}
\end{lemma}

The fourth one is an elementary sequence result (see~\cite{bertsekas1995nonlinear}).
\begin{lemma}\label{lem:A3}\quad
	Let $a_n, b_n$ be two non-negative sequences such that 
	$\sum_{n=1}^\infty a_n =\infty , \sum_{n=1}^\infty 
	a_nb_n < \infty$. 
	If there exists a real number $K > 0$ such that $|b_{n+1} - b_n| \leq Ka_n$. Then, $\lim_{n \rightarrow \infty} b_n =0$.
\end{lemma}	

The fifth one is law of large numbers for martingales given in~\cite{HH80}:
\begin{lemma}
	\label{thm:mg_convergence1}\quad
	Let $M_{\run} = \sum_{\iRun=0}^{\run} d_{\iRun}$ be a martingale adapted to $(\filter_{\run})_{\run=0}^{\infty}$ and let $(u_{\run})_{\run=0}^{\infty}$ be a nondecreasing sequence of positive numbers with $\lim_{\run\to\infty} u_{\run} = \infty$.
	If $\sum_{\run=0}^{\infty} u_{\run}^{-p} \exof{\abs{d_{\iRun}}^{p}\given\filter_{\run}} <\infty$ for some $p\in[1,2]$ \as, then:
	\begin{equation}
	\label{eq:LLN}
	\lim_{\run\to\infty} \frac{M_{\run}}{u_{\run}} = 0
	\quad
	\as
	\end{equation}
\end{lemma}

Finally, we recall the standard notion of semiflow. 
\begin{definition}\quad
	A semiflow $\flow$ on a metric space $(M, d)$ is a continuous map $\flow: \R_+ \times M \to M$:
	$$(t, x) \to \flow_t(x),$$
	such that the semi-group properties hold: $\flow_0$ = identity, $\flow_{t+s} = \flow_t \circ \flow_s$ for all $(t,s) \in \R_+ \times \R_+$.
\end{definition}

\begin{remark}\label{rem:semi_flow}\quad
	A standard way to induce a semiflow is via an \ac{ODE}.
	Specifically, if $F: \R^m \to \R^m$ is a continuous function and if the following \ac{ODE} has a unique solution trajectory for each initial point $\tilde{x} \in \R^m$:
	\begin{eqnarray*}
		\frac{dx}{dt}=&F(x),\\
		x(0)&=\tilde{x},
	\end{eqnarray*}
	then $\flow_t(\tilde{x})$ defined by the solution trajectory $x(t) \in \R^m$ as follows is a semiflow: $\flow_t(\tilde{x}) \triangleq x(t)$ with $x(0) = \tilde{x}$. We say $\flow$ defined in this way is the semiflow induced by the corresponding \ac{ODE}.
\end{remark}

\section{Technical Proofs}

\subsection{General Nonconvex Objectives}

\subsubsection{Proof of Lemma~\ref{lem:tail}}
{\em Proof:}
We start by rewriting the delayed gradient update $\act_{\run+1} = \act_{\run} - \step_{\run+1} \nabla \sobj(\act_{s(\run)},\sample_{\run+1})$ in Algorithm~\ref{alg:DASGD_un} in two forms as follows, both of which will be used subsequently:
\begin{equation}
\begin{split}
&\act_{\run+1} = \act_{\run} - \step_{\run+1}\Bigg(\grad \obj(\act_{s(\run)}) + \grad \sobj(\act_{s(\run)},\sample_{\run+1}) -  \grad \obj(\act_{s(\run)}) \Bigg)\\
& =\act_{\run} - \step_{\run+1}\Bigg( \grad \obj(\act_{\run}) + \grad \obj(\act_{s(\run)}) - \grad \obj(\act_{\run})  + \grad \sobj(\act_{s(\run)},\sample_{\run+1}) -  \grad \obj(\act_{s(\run)}) \Bigg).
\end{split}
\end{equation}

Denoted $B \triangleq \sup_{\actd \in \feas} \exof{\|\nabla \sobj(x,\sample)\|_2^2}$ per Assumption~\ref{asm:basic}.
Since $\grad f(x) = \exof{\grad \sobj(\actd, \sample)}$ is Lipschitz per Assumption~\ref{asm:basic}, letting $L$ be the Lipschitz constant, we have:
\begin{equation}\label{eq:Lip}
\begin{split}
&\obj(\act_{\run+1}) - \obj(\act_{\run}) 
\leq  \langle  \grad \obj(\act_{\run}),  \act_{\run+1} - \act_{\run}\rangle + \frac{L}{2}\|\act_{\run+1} - \act_{\run} \|_2^2\\
& = - \step_{\run+1}\langle  \grad \obj(\act_{\run}),   \nabla \sobj(\act_{s(\run)},\sample_{\run+1})\rangle + \frac{L}{2}\|\step_{\run+1} \nabla \sobj(\act_{s(\run)},\sample_{\run+1}) \|_2^2\\
& = - \step_{\run+1}\langle  \grad \obj(\act_{\run}),  \grad \obj(\act_{\run}) + \grad \obj(\act_{s(\run)}) - \grad \obj(\act_{\run})  + \grad \sobj(\act_{s(\run)},\sample_{\run+1}) -  \grad \obj(\act_{s(\run)}) \rangle \\
&+ \frac{L}{2}\step_{\run+1}^2\|  \grad \obj(\act_{s(\run)})  + \grad \sobj(\act_{s(\run)},\sample_{\run+1}) -  \grad \obj(\act_{s(\run)})  \|_2^2\\
& = - \step_{\run+1} \|\grad \obj(\act_{\run})\|_2^2 -  \step_{\run+1}\langle  \grad \obj(\act_{\run}),   \grad \obj(\act_{s(\run)}) - \grad \obj(\act_{\run})\rangle
	\notag\\
	&-  \step_{\run+1}\langle  \grad \obj(\act_{\run}), \grad \sobj(\act_{s(\run)},\sample_{\run+1}) -  \grad \obj(\act_{s(\run)}) \rangle \\
&+ \frac{L}{2}\step_{\run+1}^2\|  \grad \obj(\act_{s(\run)})  + \grad \sobj(\act_{s(\run)},\sample_{\run+1}) -  \grad \obj(\act_{s(\run)})  \|_2^2\\
& \leq - \step_{\run+1} \|\grad \obj(\act_{\run})\|_2^2 +  \step_{\run+1}\| \grad \obj(\act_{\run})\|_2  \| \grad \obj(\act_{s(\run)}) - \grad \obj(\act_{\run})\|_2 
	\notag\\
	&- \step_{\run+1}\langle  \grad \obj(\act_{\run}), \grad \sobj(\act_{s(\run)},\sample_{\run+1}) -  \grad \obj(\act_{s(\run)}) \rangle \\
&+ \frac{L}{2}\step_{\run+1}^2 \Big\{2 \|\grad \obj(\act_{s(\run)})\|_2^2  + 2\|\grad \sobj(\act_{s(\run)},\sample_{\run+1}) -  \grad \obj(\act_{s(\run)})  \|_2^2\Big\}\\
& \leq - \step_{\run+1} \|\grad \obj(\act_{\run})\|_2^2 +  \sqrt{B}\step_{\run+1}\| \grad \obj(\act_{s(\run)}) - \grad \obj(\act_{\run})\|_2
	\notag\\
	&-  \step_{\run+1}\langle  \grad \obj(\act_{\run}), \grad \sobj(\act_{s(\run)},\sample_{\run+1}) -  \grad \obj(\act_{s(\run)}) \rangle \\
& \leq L\step_{\run+1}^2 \Big\{ \sqrt{B}  + \|\grad \sobj(\act_{s(\run)},\sample_{\run+1}) -  \grad \obj(\act_{s(\run)})  \|_2^2\Big\},\\
\end{split}
\end{equation}
where in the last inequality follows because by Jensen's inequality, we have:
$$\| \grad \obj(\act_{\run})\|_2^2   \leq \sup_{\actd \in \feas} \|\exof{\nabla \sobj(x,\sample)}\|_2^2  \leq \sup_{\actd \in \feas} \exof{\|\nabla \sobj(x,\sample)\|_2^2} \leq B.$$ 

Denote the filtration generated by $\act_0, \act_1, \dots, \act_{\run}$ to be
$\mathcal{F}_{\run}$. We take the expectation of both sides of Equation~\eqref{eq:Lip} and obtain:
\begin{equation}\label{eq:Lip_ex}
\begin{split}
&\exof{\obj(\act_{\run+1}) - \obj(\act_{\run})}
\leq - \step_{\run+1} \exof{\|\grad \obj(\act_{\run})\|_2^2} +  \sqrt{B}\step_{\run+1} \exof{\| \grad \obj(\act_{s(\run)}) - \grad \obj(\act_{\run})\|_2} \\
& -  \step_{\run+1} \exof{\langle  \grad \obj(\act_{\run}), \grad \sobj(\act_{s(\run)},\sample_{\run+1}) -  \grad \obj(\act_{s(\run)}) \rangle} +  L\sqrt{B}\step_{\run+1}^2
	\notag\\
	&+ L\step_{\run+1}^2\exof{\|\grad \sobj(\act_{s(\run)},\sample_{\run+1}) -  \grad \obj(\act_{s(\run)})  \|_2^2}\\
& = - \step_{\run+1} \exof{\|\grad \obj(\act_{\run})\|_2^2} +  \sqrt{B}\step_{\run+1} \exof{\| \grad \obj(\act_{s(\run)}) - \grad \obj(\act_{\run})\|_2}
	\notag\\
	&-  \step_{\run+1} 
\ex \Bigg\{\ex \Big\{\langle  \grad \obj(\act_{\run}), \grad \sobj(\act_{s(\run)},\sample_{\run+1}) -  \grad \obj(\act_{s(\run)}) \rangle \Big | \mathcal{F}_{\run} \Big\} \Bigg\}\\
&+  L\sqrt{B}\step_{\run+1}^2   
+ L\step_{\run+1}^2\exof{\|\grad \sobj(\act_{s(\run)},\sample_{\run+1}) -  \grad \obj(\act_{s(\run)})  \|_2^2}\\
& =- \step_{\run+1} \exof{\|\grad \obj(\act_{\run})\|_2^2} +  \sqrt{B}\step_{\run+1} \exof{\| \grad \obj(\act_{s(\run)}) - \grad \obj(\act_{\run})\|_2}
	\notag\\
	&-  \step_{\run+1} 
\ex \Bigg\{\langle  \grad \obj(\act_{\run}), \ex \Big\{\grad \sobj(\act_{s(\run)},\sample_{\run+1}) -  \grad \obj(\act_{s(\run)})  \Big | \mathcal{F}_{\run}\Big\}\rangle  \Bigg\}\\
&+  L\sqrt{B}\step_{\run+1}^2   
+ L\step_{\run+1}^2\exof{\|\grad \sobj(\act_{s(\run)},\sample_{\run+1}) -  \grad \obj(\act_{s(\run)})  \|_2^2}\\
& =- \step_{\run+1} \exof{\|\grad \obj(\act_{\run})\|_2^2} +  \sqrt{B}\step_{\run+1} \exof{\| \grad \obj(\act_{s(\run)}) - \grad \obj(\act_{\run})\|_2}  
+  L\sqrt{B}\step_{\run+1}^2
	\notag\\
	&+ L\step_{\run+1}^2\exof{\|\grad \sobj(\act_{s(\run)},\sample_{\run+1}) -  \grad \obj(\act_{s(\run)})  \|_2^2}\\
& \leq - \step_{\run+1} \exof{\|\grad \obj(\act_{\run})\|_2^2} +  \sqrt{B}\step_{\run+1} \exof{\| \grad \obj(\act_{s(\run)}) - \grad \obj(\act_{\run})\|_2}  
+  L\sqrt{B}\step_{\run+1}^2   
+ 4LB\step_{\run+1}^2,\\
\end{split}
\end{equation}
where the third equality follows from $\ex \Big\{\grad \sobj(\act_{s(\run)},\sample_{\run+1}) -  \grad \obj(\act_{s(\run)})  \Big | \mathcal{F}_{\run}\Big \} = 0$, since $\sample_{\run +1}$ is independent of $\mathcal{F}_{\run}$, and
the last inequality follows from 
$\exof{\|\grad \sobj(\act_{s(\run)},\sample_{\run+1}) -  \grad \obj(\act_{s(\run)})  \|_2^2} \leq 
\sup_{\actd \in \feas} \exof{\|\grad \sobj(\actd,\sample_{\run+1}) -  \grad \obj(\actd)  \|_2^2}  \leq 4B$.
Since $\grad \obj$ is $L$-Liptchiz continuous, we have:
\begin{flalign}
\label{eq:expansion}
&\|\v(\act_{s(\run)}) - \v(\act_\run)\|_2 \leq L \|\act_{s(\run)} - \act_\run\|_2
	\notag\\
	&= L \Big\| \act_{s(\run)} - \act_{{s(\run)}+1} + \act_{{s(\run)}+1} - \act_{{s(\run)}+2} + \dots +  \act_{\run-1} - \act_{\run}  \Big\|_2
\notag\\
&= L\Big\|\sum_{r = s(\run)}^{\run-1} \Big\{\act_{r} - \act_{r+1} \Big\}\Big\|_2
= L\Big\|\sum_{r = s(\run)}^{\run-1} \step_{r+1} \nabla \sobj(\act_{s(r)},\sample_{r+1}) \Big\|_2
\notag \leq L \sum_{r = s(\run)}^{\run-1} \step_{r+1} \Big\|\nabla \sobj(\act_{s(r)},\sample_{r+1})\Big\|_2
\notag.\\
\end{flalign}
Taking the expectation of both sides of the above equation then yields:
\begin{flalign}
\label{eq:expansion1}
&\exof{\|\v(\act_{s(\run)}) - \v(\act_\run)\|_2}
\leq L \sum_{r = s(\run)}^{\run-1} \step_{r+1} \exof{\|\nabla \sobj(\act_{s(r)},\sample_{r+1})\|_2}
	\notag\\
	&\leq L \sum_{r = s(\run)}^{\run-1} \step_{r+1} \sup_{\actd \in \feas} \exof{\|\nabla \sobj(x,\sample_{r+1})\|_2} \\
& \leq L\bar{B} \sum_{r = s(\run)}^{\run-1} \step_{r+1},
\end{flalign}
where $\bar{B}  \triangleq \sup_{\actd \in \feas} \exof{\|\nabla \sobj(x,\sample_{r+1})\|_2} < \infty$ per Remark~\ref{rem:basic}.
Combining Equation~\eqref{eq:expansion1} with Equation~\eqref{eq:Lip_ex} then yields:
\begin{equation}\label{eq:Lip_ex2}
\begin{split}
&\exof{\obj(\act_{\run+1})} - \exof{\obj(\act_{\run})} \leq - \step_{\run+1} \exof{\|\grad \obj(\act_{\run})\|_2^2} +  L\bar{B}\sqrt{B}  \step_{\run+1}\sum_{r = s(\run)}^{\run-1} \step_{r+1}
+  L\sqrt{B}\step_{\run+1}^2   
+ 4LB\step_{\run+1}^2.\\
\end{split}
\end{equation}
Telescoping Equation~\eqref{eq:Lip_ex2} then yields:
\begin{equation}\label{eq:Lip_ex3}
\begin{split}
& -\infty < \inf_{\actd \in \feas} \obj(\actd)- \obj(\act_{0}) 
\le
\exof{\obj(\act_{T+1})} - \obj(\act_{0}) =
\sum_{\run=0}^T\Big\{\exof{\obj(\act_{\run+1})} - \exof{\obj(\act_{\run})} \Big\}\\
&\leq - \sum_{\run=0}^T\step_{\run+1} \exof{\|\grad \obj(\act_{\run})\|_2^2} +  L\bar{B}\sqrt{B} \sum_{\run=0}^T \Big(\step_{\run+1}\sum_{r = s(\run)}^{\run-1} \step_{r+1}\Big)
+  (L\sqrt{B}+ 4LB) \sum_{\run = 0}^T\step_{\run+1}^2.\\
\end{split}
\end{equation}
Taking $T \to \infty$, the above equation yields:
\begin{equation}\label{eq:Lip_ex3}
\begin{split}
& -\infty < 
- \sum_{\run=0}^{\infty}\step_{\run+1} \exof{\|\grad \obj(\act_{\run})\|_2^2} +  L\bar{B}\sqrt{B} \sum_{\run=0}^{\infty} \Big(\step_{\run+1}\sum_{r = s(\run)}^{\run-1} \step_{r+1}\Big)
+  (L\sqrt{B}+ 4LB) \sum_{\run = 0}^{\infty}\step_{\run+1}^2.\\
\end{split}
\end{equation}
Note that in all cases in Assumption~\ref{asm:delays}, we have $ \sum_{\run = 0}^{\infty}\step_{\run+1}^2 < \infty$.
We next proceed to bound $\sum_{\run=0}^{\infty} \Big(\step_{\run+1}\sum_{r = s(\run)}^{\run-1} \step_{r+1}\Big)$ and show that
it is finite in each of the cases in Assumption~\ref{asm:delays}.

\begin{enumerate}
	\item In the bounded delays case, since $\sup_{\run} d_{\run} \leq D$, it follows that $s(n) + D \ge n$ and hence:
	\begin{equation}
	\begin{split}
	&\sum_{\run=0}^{\infty} \Big(\step_{\run+1}\sum_{r = s(\run)}^{\run-1} \step_{r+1}\Big) \leq 
	\sum_{\run=0}^{\infty} \Big(\step_{\run+1}\sum_{r = \run - 1 -D}^{\run} \step_{r+1}\Big) \leq 
	(D+1)\sum_{\run=0}^{\infty} \Big(\step_{\run+1}  \max_{r \in \{\run - 1 -D, \dots, \run\}   } \step_{r+1}\Big) \\
	&\le
	(D+1)\sum_{\run=0}^{\infty} \Big( \max_{r \in \{\run - 1 -D, \dots, \run\}   } \step_{r+1} \cdot \max_{r \in \{\run - 1 -D, \dots, \run\}   } \step_{r+1}\Big) = (D+1)\sum_{\run=0}^{\infty} \Big(\max_{r \in \{\run - 1 -D, \dots, \run\}   } \step_{r+1}^2\Big)
	\notag\\
	&\leq  (D+1)\sum_{\run=0}^{\infty} \Big(\sum_{r =\run - 1 -D}^{\run} \step_{r+1}^2\Big)\\
	& = (D+1)^2 \sum_{\run = 0}^{\infty}\step_{\run+1}^2 < \infty,
	\end{split}
	\end{equation}
	where all the terms $\alpha_r$ are defined to be $0$ when $r$ drops below $0$.
	\item In the sublinearly growing delays case, since $d_{\run} = O(\run^p)$, it follows that $s(n) + K s^p(\run) \ge \run$ for some positive number $K$, which further implies that $s(n) + K s(\run) \ge s(n) + K s^p(\run) \ge \run$, thereby leading to
	$s(n) \ge \frac{n}{K+1}$.
	Consequently, we have:
	\begin{equation}
	\begin{split}
	&\sum_{\run=0}^{\infty} \Big(\step_{\run+1}\sum_{r = s(\run)}^{\run-1} \step_{r+1}\Big) \leq 
	\sum_{\run=0}^{\infty} \Big(\frac{1}{\run} \sum_{r = s(\run) }^{\run} \frac{1}{r}\Big) \leq 
	\sum_{\run=0}^{\infty} \Big(\frac{1}{\run}   \frac{Ks^p(\run)}{s(\run)}\Big)= K\sum_{\run=0}^{\infty} \frac{1}{\run}  s^{p-1}(\run) \leq K\sum_{\run=0}^{\infty} \frac{1}{\run}   (\frac{n}{K+1})^{p-1}\\
	& \leq  (K+1)^{-p}\sum_{\run=0}^{\infty} \run^{p-2}
	< \infty.
	\end{split}
	\end{equation}
	\item In the linearly growing delays case, since $d_{\run} = O(\run)$, it follows that $s(n) + K s(\run) \ge \run$ for some positive number $K$ and hence again
	$s(n) \ge \frac{n}{K+1}$.
	Consequently, we have:
\begin{flalign}
\sum_{\run=0}^{\infty} \Big(\step_{\run+1}\sum_{r = s(\run)}^{\run-1} \step_{r+1}\Big)
	&\leq  \sum_{\run=0}^{\infty} \Big(\frac{1}{\run \log \run} \sum_{r = s(\run) }^{\run} \frac{1}{r \log r}\Big)
	\leq \sum_{\run=0}^{\infty} \Big(\frac{1}{\run \log \run} \sum_{r = s(\run) }^{s(\run) + Ks(\run)} \frac{1}{r \log r}\Big)
	\notag\\
	&\leq  \sum_{\run=0}^{\infty} \Big(\frac{1}{\run \log \run}     \int_ {s(n)}^{s(\run) + K s(\run)} \frac{1}{r\log r} dr \Big)
	= \sum_{\run=0}^{\infty} \Big(\frac{1}{\run \log \run} \log \frac{\log (s(\run) + K s(\run))}{\log s(\run)} \Big)
	\notag\\
	&= \sum_{\run=0}^{\infty} \Big(\frac{1}{\run \log \run} \log \frac{\log (K+1) + \log s(\run) }{\log s(\run)}            \Big)
	= \sum_{\run=0}^{\infty} \Big(\frac{1}{\run \log \run} \log (1 + \frac{\log (K+1)}{\log s(\run)}  )          \Big)
	\notag\\
	& \leq \sum_{\run=0}^{\infty} \Big(\frac{1}{\run \log \run} \frac{\log (K+1)}{\log s(\run)} \Big)
	\leq \sum_{\run=0}^{\infty} \Big(\frac{1}{\run \log \run} \frac{\log (K+1)}{\log \run - \log (K+1)}
	\notag\\
	& \sim \overline{K}    \sum_{\run=0}^{\infty} \frac{1}{\run (\log \run)^2}     
	< \infty.
\end{flalign}
\item In the polynomially growing delays case, since $d_{\run} = O(\run^q)$, it follows that $s(n) + K s^q(\run) \ge \run$ for some positive number $K$. Note that in this case, $s(n) = \Omega (n^{\frac{1}{2q}})$, because otherwise, $s(n) + K s^q(\run) =  O (n^{\frac{1}{2q}} + K n^{\frac{1}{2}}) = o(n)$, which is a contradiction.
	Consequently, we have:
	\begin{flalign}
	\sum_{\run=0}^{\infty} \Big(\step_{\run+1}\sum_{r = s(\run)}^{\run-1} \step_{r+1}\Big)
	&\leq  \sum_{\run=0}^{\infty} \Big(\frac{1}{\run \log \run \log\log \run} \sum_{r = s(\run) }^{\run} \frac{1}{r \log r \log\log r}\Big)
	\notag\\
	&\leq \sum_{\run=0}^{\infty} \Big(\frac{1}{\run \log \run \log \log \run} \sum_{r = s(\run) }^{s(\run) + Ks^q(\run)} \frac{1}{r \log r \log\log r}\Big)
	\notag\\
	&\leq  \sum_{\run=0}^{\infty} \Big(\frac{1}{\run \log \run \log\log\run}  \int_ {s(n)}^{s(\run) + K s^q(\run)} \frac{1}{r\log r \log\log r} dr \Big)
	\notag\\
	&\leq \sum_{\run=0}^{\infty} \Big(\frac{1}{\run \log \run  \log\log\run} \log \frac{\log((K+1)\log s(\run) + a\log s(\run))}{\log \log s(\run)}     \Big)
	\notag\\
	& \leq \sum_{\run=0}^{\infty} \Big(\frac{1}{\run \log \run  \log\log\run} \log (1 + \frac{\log(K+1+a) }{\log \log s(\run)}  )   \Big)
	\notag\\
	&\leq \sum_{\run=0}^{\infty} \Big(\frac{1}{\run \log \run  \log\log\run} \frac{\log(K+1+a) }{\log \log s(\run)}    \Big)
	\notag\\
	& \leq \sum_{\run=0}^{\infty} \Big(\frac{1}{\run \log \run  \log\log\run} \frac{\log(K+1+a) }{\log \log \frac{n^{\frac{1}{2q}}}{K^\prime} }   \Big)
	\notag\\
	&\sim \overline{K}    \sum_{\run=0}^{\infty} \frac{1}{\run \log\run (\log\log \run)^2}     
	< \infty.
	\end{flalign}
\end{enumerate}

Consequently, in each of the above 4 cases, we have $\sum_{\run=0}^{\infty} \Big(\step_{\run+1}\sum_{r = s(\run)}^{\run-1} \step_{r+1}\Big) < \infty$. Therefore, Equation~\eqref{eq:Lip_ex3} yields
\begin{equation}\label{eq:Lip_ex4}
\begin{split}
& -\infty < 
- \sum_{\run=0}^{\infty}\step_{\run+1} \exof{\|\grad \obj(\act_{\run})\|_2^2} +  L\bar{B}\sqrt{B} \sum_{\run=0}^{\infty} \Big(\step_{\run+1}\sum_{r = s(\run)}^{\run-1} \step_{r+1}\Big)
+  (L\sqrt{B}+ 4LB) \sum_{\run = 0}^{\infty}\step_{\run+1}^2 \\
&\leq - \sum_{\run=0}^{\infty}\step_{\run+1} \exof{\|\grad \obj(\act_{\run})\|_2^2} + \overline{C},\\
\end{split}
\end{equation}
for some finite constant $\overline{C}$. Reversing the above inequality immediately yields the result:
$$\sum_{\run=0}^{\infty}\step_{\run+1} \exof{\|\grad \obj(\act_{\run})\|_2^2} < \infty.$$
\hfill$\blacksquare$

\subsubsection{Proof of Lemma~\ref{lem:diff}}

{\em Proof:}
We first recall a useful fact: for any two vectors $\mathbf{a}, \mathbf{b}$ and any finite-dimensional vector norm $\|\cdot\|$, 
\begin{equation}\label{eq:help1}
\Big|(\|a\| + \|b\|)(\|a\|-\|b\|)\Big|\leq \|a+b\|\|a-b\|.
\end{equation}
Using this fact, we can expect bound the term in question as follows:
\begin{align*}
& \bigg|\mathbb{E}[\|\nabla \obj(\act_{\run+1})\|_2^2] -
\mathbb{E}[\|\nabla \obj(\act_{\run})\|_2^2]\bigg|  = \Bigg| 
\mathbb{E} \Bigg[
\bigg(\|\nabla \obj(\act_{\run+1})\|_2 +
\|\nabla \obj(\act_{\run})\|_2 \bigg)   
\bigg(\|\nabla \obj(\act_{\run+1})\|_2 -
\|\nabla \obj(\act_{\run})\|_2 \bigg)   
\Bigg] \Bigg| \\
& \leq  
\mathbb{E} \Bigg[
\bigg|\|\nabla \obj(\act_{\run+1})\|_2 +
\|\nabla \obj(\act_{\run})\|_2 \bigg| \cdot   
\bigg|\|\nabla \obj(\act_{\run+1})\|_2 -
\|\nabla \obj(\act_{\run})\|_2 \bigg|   
\Bigg] \\
& \leq 
\mathbb{E} \Bigg[
\bigg\|\nabla \obj(\act_{\run+1}) +
\nabla \obj(\act_{\run})\bigg\|_2  \cdot   
\bigg\|\nabla \obj(\act_{\run+1}) -
\nabla \obj(\act_{\run})\bigg\|_2    
\Bigg] 
\leq \mathbb{E} \Bigg[
2\sup_{\actd \in \feas}\|\nabla \obj(\actd) \|_2 \cdot   
\bigg\|\nabla \obj(\act_{\run+1}) -
\nabla \obj(\act_{\run})\bigg\|_2    
\Bigg]\\
& \leq 2\sqrt{B}\mathbb{E} \Bigg[
\bigg\|\nabla \obj(\act_{\run+1}) -
\nabla \obj(\act_{\run})\bigg\|_2    
\Bigg] 
\leq 2\sqrt{B}\mathbb{E} \Bigg[
L\bigg\|\act_{\run+1} -
\act_{\run}\bigg\|_2    
\Bigg] 
= 2L\sqrt{B}\mathbb{E} \Bigg[
\bigg\|\act_{\run+1} -
\act_{\run}\bigg\|_2    
\Bigg] \\
&= 2L\sqrt{B}\mathbb{E} \Bigg[
\bigg\|   \step_{\run+1} \nabla \sobj(\act_{s(\run)},\sample_{\run+1})         \bigg\|_2    
\Bigg] \le
2L\sqrt{B}\step_{\run+1}\sup_{\actd \in \feas}\mathbb{E} \Bigg[
\bigg\|    \nabla \sobj(\actd,\sample_{\run+1})         \bigg\|_2    
\Bigg]
	\notag\\
	&= 2L\sqrt{B}\step_{\run+1}\sup_{\actd \in \feas}\mathbb{E} \Bigg[
\sqrt{\bigg\|   \nabla \sobj(\actd,\sample_{\run+1})         \bigg\|_2 ^2   }
\Bigg]  \\
& \leq 2L\sqrt{B}\step_{\run+1}\sup_{\actd \in \feas} 
\sqrt{\mathbb{E}\Bigg[\bigg\|   \nabla \sobj(\actd,\sample_{\run+1})         \bigg\|_2 ^2  \Bigg] }
	\notag\\
	&\leq 2L\sqrt{B}\step_{\run+1}\sqrt{B} =2LB\step_{\run+1},
\end{align*}
where the first inequality is an application of Jensen's inequality,
the second inequality follows from Equation~\eqref{eq:help1}
the fifth inequality follows from Lipschitz continuity and the second-to-last inequality follows from another application of Jensen's inequality and that the squre root function is concave.
\hfill$\blacksquare$

\subsection{Variationally Coherent Objectives}
We define the following constants that will be handy later:
\begin{enumerate}
	\item $C_ 1 = \sup_{\actd \in \feas}\exof{\dnorm{\nabla\sobj(\actd;\sample)}}$.
	\item $C_ 2 = \sup_{\actd \in \feas}\exof{\dnorm{\nabla\sobj(\actd;\sample)}^{2}}$.
	\item  $C_3$ is the Lipschitz constant for $\grad\obj(\actd) (= \exof{\nabla\sobj(x;\sample)}):$
	\begin{equation}
	\|\grad\obj(\actd) - \grad\obj(\actd^{\prime})\|_2 \leq C_3 \|\actd - \actd^{\prime}\|_2.
	\end{equation}
	\item $C_4 = \sup_{x, x^{\prime} \in \feas} \|x - x^{\prime}\|_2.$
\end{enumerate}



\subsubsection{Proof of Lemma~\ref{lem:energy}}

	Per the definition of the energy function, we have:
	\begin{equation}
	\begin{split}
	&\energy_{\sol}(\scored) - \|\mirror(\scored) - \sol\|_2^2 = \|\sol\|^2_2 -  \|\mirror(\scored)\|_2^2 + 
	2\langle\scored, \mirror(\scored) - \sol\rangle - \Big\{\|\mirror(\scored) \|_2^2 - 2\langle \mirror(\scored), \sol \rangle + \|\sol\|_2^2 \Big\} \\
	&= -2\|\mirror(\scored) \|_2^2 + 2\langle \scored - \mirror(\scored), \mirror(\scored) - \sol\rangle + 
	2\langle \mirror(\scored), \mirror(\scored) - \sol\rangle +
	2\langle \mirror(\scored), \sol \rangle \\
	& = 2\langle \scored - \mirror(\scored), \mirror(\scored) - \sol\rangle \ge 0,
	\end{split}
	\end{equation}
	where the last inequality follows from Lemma~\ref{lem:aux1}.
	Consequently, $\energy_{\sol}(\scored) \ge \|\mirror(\scored) - \sol\|_2^2 \ge 0$.
	Taking the infimum over $\solset$ yields $\energy(\scored) \ge 0$.
	
	Next, we establish that $\energy(\scored) = 0$ if and only if $\mirror(\scored) \in \solset$. The if part is already established in Remark~\ref{rem:prop1}.
	It suffices to show that $\energy(\scored) = 0$ implies $\mirror(\scored) \in \solset$.
	To see this, observe that if $\energy_{\sol}(\scored) = 0$, we must have $\|\mirror(\scored) - \sol\|_2^2 = 0$, therefore implying
	$\mirror(\scored) = \sol.$ Since $\energy_{\sol}(\scored)$ is a continuous function of $\sol$ for each fixed $\scored$, and since $\solset$ is a compact set, $\inf_{\sol \in \solset} \energy_{\sol}(\scored)$ must be achieved by
	some $z^* \in \solset$. Namely $\energy(\scored) = \energy_{z^*}(\scored)$.
	Consequently, by the preceding observation, $\energy(\scored) = 0$ implies $\mirror(\scored) = z^* \in \solset$.
	
	For the second statement, suppose on the contrary $E(\scored_{\run}) \to 0$ but
	$\mirror(\scored_{\run})$ does not converge to $\solset$. Then there must exist  a subsequence $\run_k$ such that $\mirror(\scored_{\run_k})$ is bounded away from $\solset$. 
	Denoting by $\N(\solset, \epsilon)$ the $\epsilon$-open ball around $\solset$ (i.e. $\N(\solset, \epsilon) \triangleq \{x\in\R^{\vdim} \mid \text{dist}(x, \solset) < \epsilon\}$).
	Then for the subsequence $\run_k$, there must exist some positive $\epsilon$ such that $\mirror(\scored_{\run})$ remains in $\feas \cap \N^c(\solset, \epsilon)$.
	Since $\feas \cap \N^c(\solset, \epsilon)$ is an intersection of two closed sets, it is itself closed; it is also bounded since $\feas$ is bounded: hence it is compact. Further, since for each fixed $\sol$, $E_{\sol}(\scored)$ is a continuous function 
	of $\mirror(\scored)$, $E_{\sol}(\scored)$ must achieve the minimum value on the compact set $\feas \cap \N^c(\solset, \epsilon)$, where the minimum value $a_{\sol}$ is positive per the first statement:
	$$E_{\sol}(\scored) \ge a_{\sol} > 0, \forall \mirror(\scored) \in \feas \cap \N^c(\solset, \epsilon), \forall \sol \in \solset$$
	 Finally, since $E_{\sol}(\scored)$ is continuous in $\sol$, it must achieve the minimum value (over $\sol$) on the compact set $\solset$, where the minimum value $a$ (corresponding to some $E_{\sol}(\scored)$) must again be positive:
	 $$E(\scored) = \inf_{\sol \in \solset} E_{\sol}(\scored) \ge \inf_{\sol \in \solset} a_{\sol} = a > 0.$$
	
	Consequently, $E(\scored_{\run_k}) \ge a > 0$ since $\mirror(\scored_{\run_k}) \in \feas \cap \N^c(\solset, \epsilon), \forall k$.
	However, on this subsequence, the energy function still converges to $0$ by assumption: $E(\scored_{\run_k}) \to 0$, which immediately yields a contradiction.
	The claim is hence established.
\hfill $\blacksquare$

\subsubsection{Proof of Lemma~\ref{lem:energy_prop}}
	We first prove the first claim. By expanding it, we have:
	\begin{equation}\label{eq:iden}
	\begin{split}
	&\|\mirror(\scored) - \hat{\scored}\|_2^2 - \|\mirror(\hat{\scored}) -\hat{\scored}\|_2^2 =
	\|\mirror(\scored) - \scored + \scored - \hat{\scored}\|_2^2 - \|\mirror(\hat{\scored}) -\hat{\scored}\|_2^2  \\
	& = \|\scored - \hat{\scored}\|_2^2 + \|\mirror(\scored) - \scored\|_2^2 + 2\langle \mirror(\scored) - \scored , \scored - \hat{\scored} \rangle- \|\mirror(\hat{\scored}) -\hat{\scored}\|_2^2 \\
	& = \|\scored - \hat{\scored}\|_2^2 -\Big\{ \|\mirror(\hat{\scored}) -\hat{\scored}\|_2^2  - \|\mirror(\scored) - \scored\|_2^2 - 2\langle \scored - \mirror(\scored)  ,  \hat{\scored} -\scored  \rangle\Big\}. \\
	\end{split}
	\end{equation}
	Now define the function $f(\actd, \scored) = \|\actd- \scored\|_2^2$. It follows easily that the solution to the problem $\max_{\actd \in \feas} \|\actd- \scored\|_2^2$  is $x^*(\scored) = \mirror(\scored)$.
	Consequently, by Lemma~\ref{lem:env}, $V(\scored) = f(x^*(\scored), \scored)$ is a differential function in $\scored$ and its derivative can be computed explicitly as follows:
	\begin{equation}\label{eq:gradient}
	\grad V(\scored) =  \grad f(\sol (\scored), \scored) = 2 (\scored - x^*(\scored)) = 2(\scored - \mirror(\scored)).
	\end{equation}
	
	Futher, since for each $\actd \in \feas$, $f(\actd, \scored)$ is a convex function in $\scored$, and taking the maximum preserves convexity, we have $V(\scored)$ is also a convex function in $\scored$.
	This means that $$V(\hat{\scored}) - V(\scored) - \langle \grad V(\scored), \hat{\scored} -\scored\rangle \ge 0.$$
	By Equation~\eqref{eq:gradient} and that $V(\scored) = f(x^*(\scored), \scored) = \|\mirror(\scored) -\actd\|_2^2$, the above equation becomes:
	$$ \|\mirror(\hat{\scored}) -\hat{\scored}\|_2^2  - \|\mirror(\scored) - \scored\|_2^2 - 2\langle \scored - \mirror(\scored)  ,  \hat{\scored} -\scored \rangle \ge 0.$$
	Consequently, Equation~\eqref{eq:env_opt} then immediately yields:
	$$\|\mirror(\scored) - \hat{\scored}\|_2^2 - \|\mirror(\hat{\scored}) -\hat{\scored}\|_2^2 \le
	\|\scored-\hat{\scored}\|_2^2.$$
	
	We now prove the second part. Expanding using the definition of the Lyapunov function (and skipping some tedious algebra in between), we have:
	\begin{flalign}
	\label{eq:energy_change}
	\energy_{\sol}(\scored + \Delta) - \energy_{\sol}(\scored)
	&= \|\sol\|^2_2 -  \|\mirror(\scored + \Delta)\|_2^2 + 
	2\langle\scored + \Delta, \mirror(\scored + \Delta) - \sol\rangle
	\notag\\
	&- \Big\{ \|\sol\|^2_2 -  \|\mirror(\scored)\|_2^2 + 
	2\langle\scored, \mirror(\scored) - \sol\rangle  \Big\}
	\notag\\
	& = \|\mirror(\scored)\|_2^2  - \|\mirror(\scored + \Delta)\|_2^2  - 2\langle\scored, \mirror(\scored) - \sol\rangle  +  2\langle \scored + \Delta, \mirror(\scored + \Delta) - \sol \rangle
	\notag\\
	& = \|\mirror(\scored)\|_2^2  - \|\mirror(\scored + \Delta)\|_2^2 +
	2\langle \y, \mirror(\scored + \Delta) - \mirror(\scored) \rangle
	\notag\\
	&\hfill
	+2\langle  \Delta, \mirror(\scored) - \sol + \mirror(\scored + \Delta) -\mirror(\scored)  \rangle
	\notag\\
	&= 2\langle  \Delta, \mirror(\scored) - \sol  \rangle +
	2\langle \scored + \Delta, \mirror(\scored + \Delta) - \mirror(\scored) \rangle
	+ \|\mirror(\scored)\|_2^2  - \|\mirror(\scored + \Delta)\|_2^2
	\notag\\
	& =  2\langle  \Delta,  \mirror(\scored)- \sol  \rangle +
	\|\mirror(\scored) - (\scored + \Delta)\|_2^2 - \|\mirror(\scored + \Delta) - (\scored + \Delta)\|_2^2
	\notag\\
	& \leq 2\langle  \Delta,  \mirror(\scored) - \sol  \rangle +
	\|\Delta\|_2^2
	\end{flalign}
where the last equality follows from completing the squares and the last inequality follows from the first part of the lemma.
\hfill $\blacksquare$

\subsubsection{Proof of Proposition~\ref{prop:recur_d}}
	We provide the details for all the steps.
	\begin{enumerate}
		\item 
		Defining $b_{\run} = \v(x_{s(n)}) - \v(x_\run)$,
		we can rewrite the gradient update in \ac{DAGD} as:
		\begin{flalign}
		y_{\run+1}
		&= y_{\run} - \alpha_{\run+1} \v(x_{s(n)})
		\notag\\
		&= y_{\run} - \alpha_{\run+1} \v(x_\run)  - \alpha_{\run+1} \{\v(x_{s(n)}) - \v(x_\run)\}
		\notag\\
		&= y_{\run} - \alpha_{\run+1} (\v(x_\run) + b_{\run}).
		\end{flalign}
		Recall here once again that  $s(\run)$ denotes the previous iteration count whose gradient becomes available only at the current iteration $\run$.
	To establish the claim, we start by expanding $b_{\run}$ as follows:
		\begin{flalign}
		\|b_{\run}\|_2 
		&= \|\v(\actd_{s(\run)}) - \v(\actd_\run)\|_2 \leq C_3 \|\actd_{s(\run)} - \actd_\run\|_2
		\notag\\
		& = C_3 \|\mirror(\scored_{s(\run)}) - \mirror(\scored_\run)\|_2  \leq  C_3 \|\scored_{s(\run)} - \scored_\run\|_2
		\notag\\
		& \leq  C_3 \Big\{ \|\scored_{s(\run)} - \scored_{{s(\run)}+1}\|_2 + \|\scored_{{s(\run)}+1} - \scored_{{s(\run)}+2}\|_2 + \dots +  \|\scored_{\run-1} - \scored_{\run}\|_2  \Big\}
		\notag\\
		&= C_3 \sum_{r = s(\run)}^{\run-1} \|\step_{r+1} \v(\actd_{s(r)}) \|_2 \leq  C_3\sup_{\actd \in \feas} \|\v(\actd)\|_2 \sum_{r = s(n)}^{n-1} \step_{r+1} = C_3V_{\max} \sum_{r = s(n)}^{n-1} \step_{r+1}.
		\end{flalign}
		We now consider two cases, depending on whether the delays are bounded or not.
		\begin{enumerate}
			\item
			If $\{\step_{\run}\}_{\run = 1}^{\infty}$ and $d_{\run} \leq D, \forall \run$ satisfy \cref{asm:delays}, then $d_{s(\run)} = \run - s(\run) \leq D$. Consequently, 
			\begin{equation}
			0
			< \sum_{r = s(n)}^{n-1} \step_{r+1} = \sum_{r = s(n)+1}^{n} \step_{r}\leq \sum_{r = n-D}^{n} \step_{r}
			\leq D \max_{r \in \{n-D, \dots, n\}} \step_{r} \to 0
			\quad\text{ as $n\to\infty$}, 
			\end{equation}
			where the limit approaching $0$ follows from $\lim_{\run \to \infty} \step_{\run} = 0$, which itself is a consequence of \cref{asm:delays}. This implies $\lim_{\run \to \infty} C_3V_{\max} \sum_{r = s(n)}^{n-1} \step_{r+1} = 0$ and consequently, $\lim_{\run \to \infty} \|b_{\run}\|_2 = 0$.
			
			\item 
			We consider each of the three conditions in turn.
			
			When $\step_{\run-1} = \frac{1}{\run}$  and $d_{\run}= o(\run)$,
			we have $\run - s(\run) \leq K o(s(\run))$ for some universal constant $K > 0$, which means $\run \leq s(\run) + K o(s(\run))$. 
			Consequently, we have:
			\begin{flalign}
			0
			&< \sum_{r = s(n)}^{n-1} \step_{r+1} = \sum_{r = s(n) }^{s(\run) + K o(s(\run))} \step_{r}
			\leq \int_ {s(n)}^{s(\run) + K o(s(\run))} \frac{1}{r} dr
			\notag\\
			&= \log (s(\run) + K o(s(\run))) - \log s(\run)
			= \log \frac{Ko(s(\run)) + s(\run)}{s(\run)} \to \log 1 = 0
			\quad \text{as $\run\to\infty$},
			\end{flalign}
			where the last limit follows from $s(\run) \to \infty$ as $\run \to \infty$ because
			$\run \leq s(\run) + K o(s(\run))$. 
			
			When $\step_{\run-1} = \frac{1}{\run \log \run }$  and $d_{\run}= O(\run)$,
			it is easy to verify (by integration) that this particular choice of sequence satisfies $\sum_{\run=\start}^{\infty} \step_{\run}^{2}
			<\infty,
			\sum_{\run=\start}^{\infty} \step_{\run}
			= \infty.$
			Since $d_{\run}= O(\run)$, we have $\run - s(\run) \leq K s(\run)$ for some universal constant $K > 0$, which means $\run \leq s(\run) + K s(\run)$. 
			Consequently, we have:
			\begin{flalign}
			0
			&< \sum_{r = s(n)}^{n-1} \step_{r+1} = \sum_{r = s(n) }^{s(\run) + K s(\run)} \step_{r}
			\leq \int_ {s(n)}^{s(\run) + K s(\run)} \frac{1}{r\log r} dr
			\notag\\
			&= \log \frac{\log (s(\run) + K s(\run))}{\log s(\run)}
			= \log \frac{\log (K+1) + \log s(\run) }{\log s(\run)} \to 0
			\quad \text{as $\run\to\infty$},
			\end{flalign}
			where the last limit follows from $s(\run) \to \infty$ as $\run \to \infty$ because
			$\run \leq s(\run) + K s(\run)$. 
			
			When $\step_{\run-1} = \frac{1}{\run \log \run \log\log\run}$  and $d_{\run}= O(\run^a), a>1$,
			it is again easy to verify (by integration) that this particular choice of sequence satisfies \cref{asm:delays}.
			Since $d_{\run}= O(\run^a)$, we have $\run - s(\run) \leq K s(\run)^a$ for some universal constant $K > 0$, which means $\run \leq s(\run) + K s(\run)^a$. 
			Consequently, we have:
			\begin{flalign}
			0 < \sum_{r = s(n)}^{n-1} \step_{r+1} = \sum_{r = s(n) }^{s(\run) + K s(\run)^a} \step_{r}
			&\leq \int_ {s(n)}^{s(\run) + K s(\run)^a} \frac{1}{r\log r \log\log r}dr
			\notag\\
			&= \log \frac{\log\log (s(\run) + K s(\run)^a)}{\log\log s(\run)}
			\notag\\
			&\leq  \log \frac{\log\log (s(\run)^a + K s(\run)^a)}{\log\log s(\run)}
			\notag\\
			&< \log \frac{\log((K+1)\log s(\run) + a\log s(\run))}{\log \log s(\run)}
			\notag\\
			&= \log \frac{\log(K+1+a) + \log\log s(\run)}{\log \log s(\run)} \to 0
			\quad
			\text{as $\run \to \infty$},
			\end{flalign}
			where the last limit follows from the fact that $s(\run) \to \infty$ as $\run \to \infty$ (again, because
			$\run \leq s(\run) + K s(\run)^a$). 
			
			
		\end{enumerate}
		\item 
		With the defintion of $b_{\run}$, \ac{DAGD} can be written as:
		\begin{equation}
		\begin{aligned}
		\actd_{\run}
		&= \mirror(\scored_{\run}),
		\\
		\scored_{\run+1}
		&= y_{\run} - \alpha_{\run+1} (\v(x_\run) + b_{\run}).
		\end{aligned}
		\end{equation}
To prove the claim, we start by fixing an arbitrary $\sol \in \solset$ and applying Lemma~\ref{lem:energy_prop} to bound the energy change in a single gradient update as follows:
		\begin{equation}\label{eq:energy_change}
		\begin{split}
		\energy_{\actd^*}(\scored_{\run+1}) - \energy_{\actd^*}(\scored_{\run})
		&\leq 2\langle \scored_{\run+1} - \scored_{\run}, \mirror(\scored_{\run}) -\sol \rangle
		 + 
		 \|\scored_{\run}-\scored_{\run+1}\|_2^2\\
		 &= -2\langle  \alpha_{\run+1} (\v(x_\run) + b_{\run}), \actd_{\run} - \actd^*  \rangle +
		\|\scored_{\run}-\scored_{\run+1}\|_2^2. \\
		\end{split}
		\end{equation}

		Now telescoping the above inequality yields:
		\begin{flalign}
		\label{eq:telescoping}
		\energy_{\actd^*}(\scored_{\run+1}) - \energy_{\actd^*}(\scored_0)
		&= \sum_{r=0}^{\run} \{\energy_{\actd^*}( \scored_{r+1}) - \energy_{\actd^*}(\scored_{r})\}
		\notag\\
		& \leq \sum_{r=0}^{\run} \{-2 \alpha_{r+1} \langle  \v(x_r) + b_r, \actd_{r} - \actd^*  \rangle +
		\alpha_{r+1}^2\| \v(x_r) + b_{r}\|_2^2 \}
		\notag\\
		&\leq -2 \sum_{r=0}^{\run} \alpha_{r+1}\{ \langle  \v(x_r), \actd_{r} - \actd^*  \rangle  - \|b_r\|_2 \|\actd_{r} - \actd^*  \|_2 \}+
		2\sum_{r=0}^{\run}\alpha_{r+1}^2\{\| \v(x_r)\|_2^2  + \|b_{r}\|_2^2 \}
		\notag\\
		& \leq -2 \sum_{r=0}^{\run} \alpha_{r+1}\{ \langle  \v(x_r), \actd_{r} - \actd^*  \rangle  - C_4\|b_r\|_2  \}+
		2\sum_{r=0}^{\run}\alpha_{r+1}^2\{C_2  + B \},
		\end{flalign}
		where the last inequality follows from the fact that $b_{\run}$'s must be bounded
		(since $\lim_{\run \to \infty} \|b_{\run}\|_2 = 0$) and hence let $B \triangleq  \sup_{\run} \|b_{\run}\|_2$.
		By \cref{asm:delays}, we have $2\sum_{r=0}^{\run}\alpha_{r+1}^2\{C_2  + B \} = \overline{B} < \infty$. Now fix any positive number $\epsilon$. Assume for contradiction purposes
		$\actd_{\run}$ only enters $\N(\actd^*, \epsilon)$ a finite number of times and let $t_1$ be the last time 
		this occurs. This means that for all $ r > t_1$, $\actd_r $ is outside the open set $\N(\solset, \epsilon)$.
		Therefore, since a continuous function always achieves its minimum on a compact set, we have:
		$\langle  \v(x_r), \actd_{r} - \actd^*  \rangle \ge \min_{\actd \in \feas -  \N(\solset, \epsilon)}\langle  \v(x), \actd - \actd^*  \rangle \triangleq a > 0, \forall r > t_1$ (note that $a$ depends on $\epsilon$).
		Further, since $b_r \to 0$ as $r \to \infty$, pick $t_2$ such that $\|b_r\|_2 < \frac{a}{2C_4}, \forall r \ge t_2$. Denoting $t = \max (t_1, t_2)$, we 
		continue the chain of inequalities in Equation~\eqref{eq:telescoping} below:
		\begin{equation}\label{eq:telescoping2}
		\begin{split}
		- \energy_{\actd^*}(\scored_0)
		\notag\\
		&\leq \energy_{\actd^*}(\scored_{\run+1}) - \energy_{\actd^*}(\scored_0)  \le
		-2 \sum_{r=0}^{t} \alpha_{r+1}\{ \langle  \v(x_r), \actd_{r} - \actd^*  \rangle  - C_4\|b_r\|_2  \} \\
		&  - 2 \sum_{r=t+1}^{\run} \alpha_{r+1}\{ \langle  \v(x_r), \actd_{r} - \actd^*  \rangle  - C_4\|b_r\|_2  \} +
		2\sum_{r=0}^{\run}\alpha_{r+1}^2\{C_2  + B \} \\
		& \leq  -2 \sum_{r=0}^{t} \alpha_{r+1}\{ \langle  \v(x_r), \actd_{r} - \actd^*  \rangle  - C_4\|b_r\|_2  \} 
		-2 \sum_{r=t+1}^{\run} \alpha_{r+1}\{ a - C_4\|b_r\|_2  \} + \overline{B} \\
		& \leq 2 C_4\sum_{r=0}^{t} \alpha_{r+1} |b_r\|_2   + \overline{B} - 2 \sum_{r=t+1}^{\run} \alpha_{r+1}\{ a - \frac{a}{2}  \}  \\
		& = \overline{\overline{B}} -  a\sum_{r=t+1}^{\run} \alpha_{r+1}  
		\to -\infty, \text{ as } \run \to \infty
		\end{split}
		\end{equation}
		where the first inequality follows from the energy function always being positive (\cref{lem:energy}), the second-to-last inequality follows from variational coherence and the limit on the last line follows from
		\cref{asm:delays} and that
		$\overline{\overline{B}} \triangleq  2 C_4\sum_{r=0}^{t} \alpha_{r+1} |b_r\|_2   + \overline{B}$ is just some finite constant.  
		This yields an immediate contradiction and the claim is therefore established.
	\end{enumerate}
\hfill$\blacksquare$

\subsubsection{Proof of Theorem~\ref{thm:convergence}}

	Fix a given $\delta > 0$.
	Since $\step_{\run} \to 0, b_\run \to 0$ as $\run \to \infty$, for any $a > 0$, we can pick an $N$ large enough (depending on $\delta$ and $a$) such that $\forall n \ge N$, the following three statements all hold:
	\begin{equation}\label{eq:threshold}
	\begin{split}
	& 2BC_4\step_{\run+1}  +  2\step_{\run+1}^2(C_2 + B^2) \leq \frac{\delta}{2},\\
	& C_4 \| b_\run\|_2 \leq \frac{a}{2},\\
	&  \step_{\run+1}(C_2 + B^2) < \frac{a}{2} .
	\end{split}
	\end{equation}

	We show that under either of the following two (exhaustive) possibilities, if $\energy(\scored_{\run})$ is less than $\delta$,
	$\energy(\scored_{\run+1})$ is less than $\delta$ as well, where $n \ge N$.
	\begin{enumerate}
		\item 
		Case 1: $\energy(\scored_{\run}) < \frac{\delta}{2}$.
		\item
		Case 2: $\frac{\delta}{2} \leq \energy(\scored_{\run}) < \delta$.
	\end{enumerate}
	
	Under Case 1, it follows from Equation~\eqref{eq:energy_change}:
	\begin{equation}
	\begin{split}
	&\energy_{\actd^*}(\scored_{\run+1}) -
	\energy_{\actd^*}(\scored_{\run}) \leq -2\step_{\run+1} \langle \v(\actd_\run) + b_\run, x_\run - x^*\rangle +  \step_{\run+1}^2\|\v(\actd_\run) + b_\run \|_2^2 \\
	& \leq -2\step_{\run+1} \langle b_\run, x_\run - x^*\rangle +  \step_{\run+1}^2\|\v(\actd_\run) + b_\run \|_2^2 \\
	& \leq 2\step_{\run+1} \| b_\run\|_2 \|x_\run - x^*\|_2+  2\step_{\run+1}^2(C_2 + B^2) \\
	&\leq 
	2BC_4\step_{\run+1}  +  2\step_{\run+1}^2(C_2 + B^2) 
	 \leq \frac{\delta}{2},
	\end{split}
	\end{equation}
	where the second inequality follows from variational coherence.
	Taking the infimum of $\sol$ over $\solset$ then yields:
	 $\energy(\scored_{\run+1}) -
	 \energy(\scored_{\run}) \leq \frac{\delta}{2}.$
	This then implies that $\energy{(\actd^*, \scored_{\run+1})} \leq 
	\energy(\actd^*, \scored_{\run}) +  \frac{\delta}{2} < \delta$.
	
	Under Case 2, \cref{eq:energy_change} readily yields:
	\begin{equation}
	\begin{split}
	& \energy_{\actd^*}(\scored_{\run+1}) -
	\energy_{\actd^*}(\scored_{\run}) \leq -2\step_{\run+1} \langle \v(\actd_\run) + b_\run, x_\run - x^*\rangle +  \step_{\run+1}^2\|\v(\actd_\run) + b_\run \|_2^2 \\
	& = -2\step_{\run+1} \langle \v(\actd_\run), x_\run - x^*\rangle -2\step_{\run+1} \langle b_\run, x_\run - x^*\rangle +  \step_{\run+1}^2\|\v(\actd_\run) + b_\run \|_2^2 \\
	& \leq -2\step_{\run+1} a + 2\step_{\run+1} \| b_\run\|_2 \|x_\run - x^*\|_2+  2\step_{\run+1}^2(C_2 + B^2) \\
	& \leq -2\step_{\run+1} \Big\{a - C_4 \| b_\run\|_2  - \step_{\run+1}(C_2 + B^2)\Big\} \\
	& \leq -2\step_{\run+1} \Big\{a - \frac{a}{2}  - \step_{\run+1}(C_2 + B^2)\Big\} \\
	& = -2\step_{\run+1} \Big\{ \frac{a}{2}  - \step_{\run+1}(C_2 + B^2)\Big\} \\
	& < 0,
	\end{split}
	\end{equation}
	where the second inequality follows from $\langle \v(\actd_\run), x_\run - x^*\rangle \ge a$ under Case 2\footnote{Here $a$ is a constant that only depends on $\delta$: if $\energy{(\actd^*, y)} \ge \frac{\delta}{2}$, then by part 2 of \cref{lem:energy}, $\mirror(y)$ must be outside an $\epsilon$-neighborhood of $\actd^*$, for some $\epsilon > 0$. On this neighborhood, the strictly positive continuous function $\langle \v(\actd), x - x^*\rangle$ must achieve a minimum value $a > 0$.}. Taking the infimum of $\sol$ over $\solset$ then yields:
	$\energy(\scored_{\run+1}) \le
	\energy(\scored_{\run}) .$
	This then implies that $\energy{( \scored_{\run+1})} \leq 
	\energy(\scored_{\run}) <  \frac{\delta}{2}$.
	
	Consequently, putting the above two cases together therefore yields $\energy{(\scored_{\run+1})} < \delta$.
\hfill{$\blacksquare$}

\subsubsection{Proof of Proposition~\ref{prop:recur_s}}

Using the definitions introduced in the main text, we rewrite the gradient update in \ac{DASGD} as:
		\begin{equation}\label{eq:DSLGD_abbre}
		\score_{\run+1} = \score_{\run} - \step_{\run+1}\{\v(\act_\run) + B_{\run} + U_{\run+1}\}.
		\end{equation} 
		
		\begin{enumerate}
			\item
		To see that 
		$\sum_{r=0}^{\run} U_{\run+1}$ is a martingale adapted to $\score_0, \score_1 \dots, \score_{\run+1}$,
		first note that, by defintion, $B_{\run}$ is adapted to $\score_0, \score_1 \dots, \score_{\run}$ (since $\act_{\run}$ is a deterministic function of $\score_{\run}$) and $\score_{\run+1}, \score_{\run}, B_{\run}$ together determine $U_{\run+1}$.
		We then check that their first moments are bounded: 
		\begin{equation}
		\begin{split}
		&\exof{\|\sum_{r=0}^ \run \|U_{r+1} \|_2} \leq   \sum_{r=0}^\run \exof{ \|U_{r+1} \|_2} = \sum_{r=0}^\run \exof{\| \nabla \sobj(\act_{s(r)},\sample_{r+1}) -  \v(\act_{s(r)})\|_2}\\
		& \leq \sum_{r=0}^\run \Big\{\exof{\| \nabla \sobj(\act_{s(r)},\sample_{r+1})\|_2} +  \exof{\|\v(\act_{s(r)})\|_2}\Big\} \\
		&\leq 
		\sum_{r=0}^\run \Big\{\sup_{\actd \in \feas} \exof{\| \nabla \sobj(\actd,\sample)\|_2} +  \sup_{\actd \in \feas}\|\v(\actd)\|_2\Big\}  \\
		&= \sum_{r=0}^\run \Big\{\sup_{\actd \in \feas} \exof{\| \nabla \sobj(\actd,\sample)\|_2} +  \sup_{\actd \in \feas}\|\exof{\nabla \sobj(\actd, \sample)}\|_2\Big\}\\
		& \leq \sum_{r=0}^\run 2 \sup_{\actd \in \feas} \exof{\| \nabla \sobj(\actd,\sample)\|_2} = \sum_{r=0}^\run 2 C_1
		= 2(n+1) C_1 < \infty,
		\end{split}
		\end{equation}
		where the last inequality follows from Jensen's inequality (since $\|\cdot\|_2$ is a convex function).
		Finally, the martingale property holds because:
		$\exof{\sum_{r=0}^ \run U_{r+1} \mid \score_1, \dots, \score_{\run+1}} = 
		\exof{\nabla \sobj(\act_{s(\run)},\sample_{\run+1}) -  \v(\act_{s(\run)}) \mid \score_1, \dots, \score_{\run}}+ \sum_{r=0}^ {\run-1} U_{r+1}   = \sum_{r=0}^ {\run-1} U_{r+1}.$, since $\sample_{\run+1}$ is \textit{iid}.
		Therefore, $\sum_{r=0}^{\run} U_{r+1}$ is a martingale, and $U_{\run+1}$ is a martingale difference sequence adapted to $\score_0, \score_1 \dots, \score_{\run+1}$.

		Next, we show that $\lim_{\run \to \infty} \|B_{\run}\|_2 = 0, { a.s.}.$ By definition, we can expand $B_{\run}$ as follows:
		\begin{flalign}
		\label{eq:expansion}
		\|B_{\run}\|_2
		&= \|\v(\act_{s(\run)}) - \v(\act_\run)\|_2 \leq C_3 \|\act_{s(\run)} - \act_\run\|_2
		= C_3 \|\mirror(\score_{s(\run)}) - \mirror(\score_\run)\|_2
		\notag\\
		&\leq  C_3 \|\score_{s(\run)} - \score_\run\|_2
		= C_3 \Big\| \score_{s(\run)} - \score_{{s(\run)}+1} + \score_{{s(\run)}+1} - \score_{{s(\run)}+2} + \dots +  \score_{\run-1} - \score_{\run}  \Big\|_2
		\notag\\
		&= C_3\Big\| \sum_{r = s(\run)}^{\run-1} \Big\{\score_{r} - \score_{r+1} \Big\}\Big\|_2
		= C_3\Big\|\sum_{r = s(\run)}^{\run-1} \step_{r+1} \nabla \sobj(\act_{s(r)},\sample_{r+1}) \Big\|_2
		\notag\\
		&= C_3\Big\|\sum_{r = s(\run)}^{\run-1} \step_{r+1} \Big\{\v(\act_{s(r)}) + \nabla \sobj(\act_{s(r)},\sample_{r+1}) -  \v(\act_{s(r)})\Big\}\Big\|_2
		\notag\\
		&= C_3\Big\|\sum_{r = s(\run)}^{\run-1} \step_{r+1} \v(\act_{s(r)}) +\sum_{r = s(\run)}^{\run-1} \step_{r+1} U_{r+1}\Big\|_2
		\notag\\
		&\leq  C_3 \sum_{r = s(\run)}^{\run-1} \step_{r+1} \|\v(\act_{s(r)})\|_2 +  C_3\Big\| \sum_{r = s(\run)}^{\run-1}\step_{r+1}U_{r+1}\Big\|_2
		\notag\\
		&\leq C_3C_1 \sum_{r = s(\run)}^{\run-1} \step_{r+1} + C_3 \|\sum_{r = s(\run)}^{\run-1} \step_{r+1} U_{r+1}\|_2
		\notag\\
		& = C_3C_1 \sum_{r = s(\run)}^{\run-1} \step_{r+1} + C_3 \|\sum_{r = 0}^{\run-1} \step_{r+1} U_{r+1} -\sum_{r = 0}^{s(\run)-1} \step_{r+1} U_{r+1}\|_2,
		\end{flalign}
		where the first inequality follows from $\grad\obj$ being Liptichz-continuous (Assumption~\ref{assump:0.6}) and the second inequality follows from $\mirror$ is a non-expansive map.
		
		By the same analysis as in the deterministic case, 
		the first part of the last line of Equation~\eqref{eq:expansion} converges to $0$ (under each one of the conditions on step-size and delays in \cref{asm:delays}):
		\begin{equation}\label{eq:limit1}
		\lim_{\run \to \infty} C_3C_1 \sum_{r = s(n)}^{n-1} \step_{r+1} = 0.
		\end{equation}

		We then analyze the limit of $\|\sum_{r = 0}^{\run-1} \step_{r+1} U_{r+1} -\sum_{r = 0}^{s(\run)-1} \step_{r+1} U_{r+1}\|_2$. 
		Define:
		$$ M_{\run} = \sum_{r = 0}^{\run-1} \step_{r+1} U_{r+1}.$$
		Since $U_{r+1}$'s are martingale differences, $M_{\run}$ is a martingale.
		Further, in each of the three conditions, $\sum_{\run=\start}^{\infty} \step_{\run}^{2}
		<\infty$. This implies that $M_{\run}$ is an $L_2$-bounded martingale because:
		\begin{equation}\label{eq:reason}
		\begin{split}
		& \sup_{\run} \exof{\|M_{\run}\|_2^2} = \sup_{\run} \exof{\Big\|\sum_{r = 0}^{\run-1} \step_{r+1} U_{r+1}\Big\|_2^2}
		= \sup_{\run} \exof{\langle \sum_{r = 0}^{\run-1} \step_{r+1} U_{r+1}, \sum_{r = 0}^{\run-1} \step_{r+1} U_{r+1}\rangle}\\
		&= \sup_{\run} \exof{\sum_{i,j} \langle  \step_{i+1} U_{i+1},  \step_{j+1} U_{j+1}\rangle} = \sup_{\run} \sum_{r =0}^{\run-1} \exof{\langle  \step_{r+1} U_{r+1},  \step_{r+1} U_{r+1}\rangle} \\
		& =  \sup_{\run} \sum_{r =0}^{\run-1} \step_{r+1}^2\exof{\| U_{r+1}\|_2^2}  \leq \sup_{\run} 4C_2\sum_{r =0}^{\run-1}  \step_{r+1}^2 \leq 4C_2 \sum_{r =0}^{\infty}  \step_{r+1}^2 < \infty,
		\end{split}
		\end{equation}
		where the last inequality in the second line follows from the martingale property as follows:
		\begin{equation}
		\begin{split}
		&\exof{\langle  \step_{i+1} U_{i+1},  \step_{j+1} U_{j+1}\rangle} =
		\step_{i+1} \step_{j+1} \exof{\langle U_{i+1},  U_{j+1}\rangle} \\
		&= \step_{i+1} \step_{j+1} \exof{\exof{\langle U_{i+1},  U_{j+1}\rangle \mid Y_0, Y_1, \dots, Y_{i+1}}} \\
		& = \step_{i+1} \step_{j+1} \exof{\langle U_{i+1}, \exof{U_{j+1} \mid Y_0, Y_1, \dots, Y_{i+1}} \rangle} =
		\step_{i+1} \step_{j+1} \exof{\langle U_{i+1}, 0\rangle} = 0,
		\end{split}
		\end{equation}
		where the  second equality follows from the tower property (and without loss of generality, we have assumed $i < j$, the third equality follows from $U_{i+1}$ is adapted to $  Y_0, Y_1, \dots, Y_{i+1}$ and the second-to-last equality follows
		from $U_{\run +1}$ is a martingale difference.
		Consequently, all the cross terms in the second line of Equation~\eqref{eq:reason} are $0$.
		Therefore, by \cref{thm:mg_convergence2}, by taking $p=2$
		$\lim_{\run \to \infty} M_{\run} = M_{\infty}, \text{ a.s.} $, where $M_{\infty}$ has finite second-moment.
		Further, since in all three cases $s(\run) \to \infty$ as $\run \to \infty$ (because there is at most a polynomial lag between $s(\run)$ and $\run$), we have $\lim_{\run \to \infty} M_{s(\run)} = M_{\infty}, \text{ a.s.} $. Therefore 
		$$\lim_{\run \to \infty} \Big\{\sum_{r = 0}^{\run-1} \step_{r+1} U_{r+1} -\sum_{r = 0}^{s(\run)-1} \step_{r+1} U_{r+1}\Big\}= \lim_{\run \to \infty} \Big\{M_{\run} - M_{s(\run)}\Big\}= 0, \text{ a.s.},$$
		thereby implying:
		\begin{equation}\label{eq:limit2}
		\lim_{\run \to \infty} C_3 \Big\|\sum_{r = 0}^{\run-1} \step_{r+1} U_{r+1} -\sum_{r = 0}^{s(\run)-1} \step_{r+1} U_{r+1}\Big\|_2 = 0.
		\end{equation}
		Combining Equation~\eqref{eq:limit1} and Equation~\eqref{eq:limit2} yields
		$\lim_{\run \to \infty} \|B_{\run}\|_2 = 0, { a.s.}.$

		\item The full \ac{DASGD} update is then:
		\begin{eqnarray}
		& \act_{\run} = \mirror(\score_{\run}) \\
		& \score_{\run+1} = \score_{\run} - \step_{\run+1}\{\v(\act_\run) + B_{\run} + U_{\run+1}\}. 
		\end{eqnarray}
		We now bound the one-step change of the energy function 
		$\energy(\solset, \score_{\run+1}) - \energy(\solset, \score_{\run})$ (which is now a random quantity) and then telescope the differences.
		
		Pick an arbitrary $\sol \in \solset$ and apply Lemma~\ref{lem:energy_prop}, we have:
		\begin{equation}\label{eq:energy_changes}
		\begin{split}
		&\energy_{\actd^*}(\score_{\run+1}) - \energy_{\actd^*}(\score_{\run})  \leq -2\alpha_{\run+1} \langle  \v(\act_\run) + B_{\run} + U_{\run+1}, \act_{\run} - \actd^*  \rangle +
		\|\score_{\run}-\score_{\run+1}\|_2^2 \\
		& = -2 \alpha_{\run+1} \langle  \v(\act_\run) + B_{\run} + U_{\run+1}, \act_{\run} - \actd^*  \rangle +
		\alpha_{\run+1}^2 \|  \v(\act_\run) + B_{\run} + U_{\run+1}  \|_2^2 \\
		& \leq -2 \alpha_{\run+1} \langle  \v(\act_\run) + B_{\run} + U_{\run+1}, \act_{\run} - \actd^*  \rangle +
		3\alpha_{\run+1}^2 \Big\{\|  \v(\act_\run)\|_2^2 + \|B_{\run} \|_2^2+ \|U_{\run+1}  \|_2^2 \Big\} \\
		& \leq -2 \alpha_{\run+1} \langle  \v(\act_\run) + B_{\run} + U_{\run+1}, \act_{\run} - \actd^*  \rangle +
		3\alpha_{\run+1}^2 \Big\{C_2 + \|B_{\run}\|_2^2+ \|U_{\run+1}  \|_2^2 \Big\}.
		\end{split}
		\end{equation}
		
		For contradiction purposes assume $\act_{\run}$ enters $\N(\solset, \epsilon)$ only a finite number of times with positive probability. By starting the sequence at a later index if necessary, we can without loss of generality
		$\act_{\run}$ never enters $\N(\actd^*, \epsilon)$ with positive probability.  Then on this event (of $\act_{\run}$ never entering $\N(\solset, \epsilon)$), we have $ \langle  \v(\act_\run), \act_{\run} - \actd^*  \rangle \ge a > 0$ as before. Telescoping Equation~\eqref{eq:energy_changes} then yields:
		\begin{equation}\label{eq:telescopings}
		\begin{split}
		& -\infty <  - \energy_{\actd^*}(\score_0)  \leq \energy_{\actd^*}(\score_{\run+1}) - \energy_{\actd^*}(\score_0)  = \sum_{r=0}^{\run} \{\energy_{\actd^*}(\score_{r+1}) - \energy_{\actd^*}(\score_{r})\}  \\
		& \leq -2\sum_{r=0}^{\run} \alpha_{\run+1} \langle  \v(\act_\run) + B_{\run} + U_{\run+1}, \act_{\run} - \actd^*  \rangle + 3\sum_{r=0}^{\run} \alpha_{\run+1}^2 \Big\{C_2 + \|B_{\run}\|_2^2+ \|U_{\run+1}  \|_2^2 \Big\}\\
		& \leq -2\sum_{r=0}^{\run} \alpha_{\run+1} \Big\{a +\langle B_{\run} + U_{\run+1}, \act_{\run} - \actd^*  \rangle \Big\}+
		3\sum_{r=0}^{\run} \alpha_{\run+1}^2 \Big\{C_2 + \|B_{\run}\|_2^2+ \|U_{\run+1}  \|_2^2 \Big\}\\
		& \to -\infty \text{ a.s. } \text{ as } \run \to \infty.
		\end{split}
		\end{equation}
		
		We justify the last-line limit of Equation~\eqref{eq:telescopings} by looking at each of its components in turn:
		\begin{enumerate}
			\item Since $\sum_{r=0}^{\run} \alpha_{\run+1}^2 < \infty$, and per the previous step, $\lim_{\run \to \infty} \|B_{\run}\|_2^2 = 0, \text{ a.s.}$, we have $3\sum_{r=0}^{\infty} \alpha_{\run+1}^2 \Big\{C_2 + \|B_{\run}\|_2^2\Big\} = C, \text{ a.s.}$,
			for some constant $C < \infty$.
			
			\item $\sum_{r=0}^{\run} \alpha_{\run+1}^2 \|U_{\run+1}  \|_2^2$ is submartingale that is $L_1$ bounded since:
			\begin{equation}
			\begin{split}
			&\sup_{\run}\exof{\sum_{r=0}^{\run} \alpha_{\run+1}^2 \|U_{\run+1}  \|_2^2} \leq \sup_{\run} \sum_{r=0}^{\run} \alpha_{\run+1}^2 \exof{\|U_{\run+1}  \|_2^2} \leq \sup_{\run} \sum_{r=0}^{\run} \alpha_{\run+1}^2 \exof{\|U_{\run+1}  \|_2^2} \\
			& = 
			\sup_{\run} \sum_{r=0}^{\run} \alpha_{\run+1}^2 \exof{\| \nabla \sobj(\act_{s(\run)},\sample_{\run+1}) -  \v(\act_{s(\run)})       \|_2^2} 
			\notag\\
			&\leq 2\sup_{\run} \sum_{r=0}^{\run} \alpha_{\run+1}^2\Big\{ \exof{\|\grad \sobj(\act_{s(\run)}, \sample_{\run+1}) \|_2^2} + \exof{\|\grad \obj(\act_{s(\run)}) \|_2^2} \Big\} \\
			& \leq 2\sup_{\run} \sum_{r=0}^{\run} \alpha_{\run+1}^2\Big\{ \sup_{\x \in \feas}\exof{\|\grad \sobj(x, \sample) \|_2^2} + \sup_{\x \in \feas}\|\grad \obj(\x) \|_2^2 \Big\} \\
			& \leq  2\sup_{\run} \sum_{r=0}^{\run} 2C_2 \alpha_{\run+1}^2 < \infty.
			\end{split}
			\end{equation}
			
			Consequently, by martingale convergence theorem (Lemma~\ref{thm:mg_convergence2} by taking $p=1$), $3\sum_{r=0}^{\run} \alpha_{\run+1}^2 \|U_{\run+1}  \|_2^2 \to R, \text {a.s.}$, for some
			random variable $R$ that is almost surely finite (in fact $\exof{|R|} < \infty$).
			
			\item Since $\|B_{\run}\|_2$ converges to $0$ almost surely, its average also converges to $0$ almost surely: $$\sum_{\run= 0}^{\infty} \frac{\step_{\run+1}\|B_{\run}\|_2}{\sum_{r=1}^{\run} \step_{r+1}} = 0, \text{ a.s.},$$
			there by implying that $$\sum_{\run= 0}^{\infty} \frac{\step_{\run+1}\langle B_{\run}, \act_{\run} - \actd^*\rangle}{\sum_{r=1}^{\run} \step_{r+1}} = 0, \text{ a.s.},$$
			since $|\langle B_{\run}, \act_{\run} - \actd^*\rangle| \leq \|B_{\run}\|_2\|\act_{\run} - \actd^*\|_2 \leq C_4 \|B_{\run}\|_2$.
			
			In addition,  $\step_{\run+1}\langle U_{\run+1}, \act_{\run} - \actd^*  \rangle$ is a martingale difference that is $L_2$ bounded because $\step_{\run+1}$ is square summable and
			$$\exof{\|\langle U_{\run+1}, \act_{\run} - \actd^*  \rangle\|_2^2} \leq \exof{\|U_{\run+1}\|_2^2 \|\act_{\run} - \actd^* \|_2^2} \leq C_4 \exof{\|U_{\run+1}\|_2^2} \leq 4C_4C_2 < \infty.$$ 
			
			Consequently, by applying Lemma~\ref{thm:mg_convergence1} with $p = 2$ and $u_n = \sum_{r=1}^{\run} \step_{r+1}$ (which is not summable), law of large number for martingales therefore implies:
			$$\sum_{\run= 0}^{\infty} \frac{\step_{\run+1}\langle U_{\run+1}, \act_{\run} - \actd^*\rangle}{\sum_{r=1}^{\run} \step_{r+1}} = 0, \text{ a.s.}$$
			Combining the above two limits, we have
			$$\lim_{\run \to \infty} 
			\frac{\sum_{r=0}^{\run} \alpha_{\run+1}\langle B_{\run} + U_{\run+1}, \act_{\run} - \actd^*  \rangle}{\sum_{r=0}^{\run} \step_{r+1}} = 0, \text { a.s.}$$
			
			Consequently, $-\sum_{r=0}^{\run} \alpha_{\run+1} \Big\{a +\langle B_{\run} + U_{\run+1}, \act_{\run} - \actd^*  \rangle \Big\} =-\{\sum_{r=0}^{\run} \alpha_{\run+1}\} \Big\{a +  \frac{\sum_{r=0}^{\run} \alpha_{\run+1}\langle B_{\run} + U_{\run+1}, \act_{\run} - \actd^*  \rangle}{\sum_{r=0}^{\run} \step_{r+1}}\Big\} \to -\infty$,
			as $\run \to -\infty$.
		\end{enumerate}
	\end{enumerate}

\subsubsection{Proof of Lemma~\ref{lem:odeprop}}

	The first claim follows by computing the derivative of the energy function with respect
	to time (for notational simplicity, here we just use $\scored(t)$ to denote $\flow(t, \by))$:
	\begin{equation}
	\label{eq:dFench}
	\begin{split}
	&\frac{d}{dt} \fench_{\sol}(y(t)) =\frac{d}{dt} \Big\{ \|\sol\|^2_2 -  \|\mirror(\scored(t))\|_2^2 + 
	2\langle\scored(t), \mirror(\scored(t)) - \sol\rangle \Big\} \\
	& = \frac{d}{dt} \Big\{  - \|\mirror(\scored(t)) - \scored(t)\|_2^2 + \|\scored(t)\|_2^2+
	2\langle\scored(t),  - \sol\rangle \Big\} \\
	& =    2\langle \dot{\scored}(t), \mirror(\scored(t)) - \scored(t)\rangle + 2\langle \scored(t), \dot{\scored}(t)\rangle+
	2\langle\dot{\scored}(t),  - \sol\rangle \Big\} \\
	& =    - 2\langle \grad\obj(\actd(t)), \actd(t) - \scored(t)\rangle - 2\langle \grad\obj(\actd(t)),  \scored(t)\rangle-
	2\langle \grad\obj(\actd(t)),  - \sol\rangle \Big\} \\
	&= -\braket{\grad\obj(\actd(t))}{x(t) - \sol}
	\leq 0,
	\end{split}
	\end{equation}
	where the last inequality is strict unless $\mirror(y(t)) = x(t) = \sol$. Take the infimum over $\sol$ then yields the result: in particular, if $\mirror(y(t)) = x(t) \notin \solset$, then $\frac{d}{dt} \fench_{\sol}(y(t)) \leq -\epsilon < 0, \forall \sol \in \solset$, hence yielding the strict part of the inequality.
	Note also that even though $\mirror(\scored(t)) - \scored(t)$ is not differentiable, 
	$\|\mirror(\scored(t)) - \scored(t)\|_2^2$ is; and in computing its derivative, we applied the envelope theorem as given in Lemma~\ref{lem:env}.

	For the second claim, consider all $\scored$ that satisfy $\energy (\flow(t, \scored)) >  \frac{\delta}{2}$. Fix any $\sol \in \solset$.
	By the monotonicity property in the first part of the lemma, it follows that
	$\energy_{\actd^*}(\flow(s, \scored)) >  \frac{\delta}{2}, \forall 0\leq s \leq t$.
	Consequently, $P(s, y)$ must be outside some $\epsilon$ neighborhood of $\solset$ for $0 \leq s \leq t$, for otherwise, $\energy_{\actd^*} (\flow(t, \scored))$ would be $0$ for at least some $\sol \in \solset$, which is a contradiction to
	$\energy (\flow(t, \scored)) >  \frac{\delta}{2}$.
	
	This means that there exists some positive constant $a(\delta)$ such that $\forall 0\leq s \leq t, \forall \sol \in \solset$:
	\begin{equation}
	\label{eq:dFench-far}
	\frac{d}{ds} \fench_{\sol}(\flow(s, y))
	= -\braket{\grad\obj(x(s))}{x(s) - \sol}
	\leq -a(\delta).
	\quad
	\end{equation}
	
	Consequently, pick $T(\delta) = \frac{\delta}{2a(\delta)}$,
	Equation~\eqref{eq:dFench-far} implies that for any $t > T(\delta)$:
	\begin{equation}
	\label{eq:SA6}
	\fench_{\sol}(\flow(t, y)) \leq \fench_{\sol}(\flow(T(\delta), y)) 
	\leq \fench_{\sol}(y) - T(\delta) a(\delta)
	\leq \fench_{\sol}(y) - \frac{\delta}{2}.
	\end{equation}
	Taking the supremum over $\sol \in \solset$ then yields:
	\begin{equation}\label{eq:SA7}
	\fench(\flow(t, y))
	\leq \fench(y) - \frac{\delta}{2}.
	\end{equation}
	Since Equation~\eqref{eq:SA7} is true for any $\scored$, taking $\sup$ over $\scored$ establishes the claim.
	\hfill$\blacksquare$

\subsection{Proof of Lemma~\ref{lem:apt}}
First, from previous analysis, we already know 
that $U_{\run+1}$ is a martingale difference sequence with $\sup_{\run}\exof{\|U_{\run+1}\|_2^2} < \infty$ and that $\step_{\run}$ is a square-summable sequence. Recall also $\lim_{\run \to \infty} B_{\run} = 0$ \as.
Now, by combining Proposition 4.1 and Proposition 4.2 in \citet{Ben99},  
we obtain the following sufficient condition for APT:

The affine interpolation curve of the iterates generated by the difference equation
	$\score_{\run+1} = \score_{\run} - \step_{\run+1}\{G(\act_\run) + U_{\run+1}\}$ is an \ac{APT} for the solution to the ODE $\dot{\scored} = -G(\scored)$ if the following three conditions \textbf{all} hold:
	\vspace{2ex}
	\begin{enumerate}
		\item $G$ is Lipschitz continuous and bounded.
			\footnote{This condition is also sufficient for the existence and uniqueness of the ODE solution}.
		\item  $U_{\run+1}$ is a martingale difference sequence with $\sup_{\run}\exof{\|U_{\run+1}\|_2^p} < \infty$
		and $\sum_{\run = 0}^{\infty} \step_{\run+1}^{1 +\frac{p}{2}} < \infty$ for some $p >  2$.
	\end{enumerate}

\vspace{2ex}
By using a similar analysis as in~\citet{Ben99} (which we omit here due to space limitation), we can show that
if the following conditions all hold:
	\begin{enumerate}
	\item $G$ is Lipschitz continuous and bounded.
	\item $\lim_{\run \to \infty} B_{\run} = 0$ \as.
	\item  $U_{\run+1}$ is a martingale difference sequence with $\sup_{\run}\exof{\|U_{\run+1}\|_2^p} < \infty$
	and $\sum_{\run = 0}^{\infty} \step_{\run+1}^{1 +\frac{p}{2}} < \infty$ for some $p >  2$,
\end{enumerate}
then, the affine interpolation curve of the iterates generated by the difference equation
$\score_{\run+1} = \score_{\run} - \step_{\run+1}\{G(\act_\run) + B_{\run} + U_{\run+1}\}$ is an \ac{APT} for the solution to the ODE $\dot{\scored} = -G(\scored)$.

Take $p=2$ and recall $\v(\mirror(\cdot))$ is Lipschitz continuous and bounded.
The above list of three conditions are thus all verified, thereby yielding the result.

\subsubsection{Proof of Theorem~\ref{thm:main}}
	By Proposition~\ref{prop:recur_s}, $\score_{\run}$ will get arbitrarily close to $\solset$ infinitely often.	
	It then suffices to show that, after long enough iterations, if $\score_{\run}$ ever gets $\epsilon$-close to $\solset$, all the ensuing iterates will be $\epsilon$-close to $\solset$ almost surely.
	The way we show this ``trapping" property is to use the energy function. Specifically, we consider $\energy(\sol,A(t))$ and show that
	no matter how small $\epsilon$ is, for all sufficiently large $t$, if $\energy(\sol,A(t_0))$ is less than $\epsilon$ for some $t_0$, then $\energy(\sol,A(t)) < \epsilon, \forall t> t_0$.
	This would then complete the proof because $A(t)$ actually contains all the \ac{DASGD} iterates, and hence if 
	$\energy(\sol,A(t)) < \epsilon, \forall t> t_0$, then $\energy(\sol,\score_{\run}) < \epsilon$ for all sufficiently large $\run$. Furthermore, since $A(t)$ contains all the iterates, the hypothesis  that `` if $\energy(\sol,A(t_0))$ is less than $\epsilon$ for some $t_0$" will be satisfied due to Prop~\ref{prop:recur_s}.

We now flesh out more details of the proof. Fix any $\epsilon > 0$.
	Since $A(t)$ is an asymptotic pseudotrajectory for $\flow$,  we have:
	\begin{equation}
	\lim_{t \to \infty} \sup_{0 \leq h \leq T} \dnorm{\score(t+h) - \flow(h,\score(t))}
	= 0.
	\end{equation}
	Consequently, for any $\delta > 0$,	there exists some $\tau(\delta,T)$ such that $\dnorm{\score(t+h) - \flow(h,\score(t))} < \delta$ for all $t\geq\tau$ and all $h \in [0, T]$. We therefore have the following 
	chain of inequalities:
	\begin{flalign}
	&\fench_{\sol}(A(t+h))
	= \fench_{\sol}(\flow(h,A(t)) + A(t+h) -\flow(h,A(t)))\\
	&\leq \fench_{\sol}(\flow(h,A(t)))
	+ \braket{A(t+h) - \flow(h,A(t))}{\mirror(\flow(h,A(t))) - \sol}
	+ \frac{1}{2} \dnorm{A(t+h) -\flow(h,A(t))}^{2} 
	\notag\\
	&\leq \fench_{\sol}(\flow(h,A(t)))
	+ C_4\delta
	+ \frac{1}{2}\delta^2
	= \fench_{\sol}(\flow(h,A(t))) + \frac{\eps}{2},
	\end{flalign}
	where in the last step we have choosen $\delta$ small enough such that $ C_4\delta
	+ \frac{1}{2}\delta^2 =\frac{\eps}{2}$.
	
	Now by Proposition~\ref{prop:recur_s}, there exists some $\tau_{0}$ such that
	$\fench_{\sol}(A(\tau_{0}))
	< \frac{\eps}{2}.$
	Our goal is to establish that $\fench_{\sol}(A(\tau_{0} + h)) < \eps$ for all $h\in[0, \infty)$.
	To that end, partition the $[0, \infty)$ into disjoint time intervals of the form $[(n-1) T_{\eps},nT_{\eps})$ for some appropriate $T_{\eps}$.
	By \cref{lem:odeprop}, we have:
	\begin{equation}
	\label{eq:bound1}
	\fench_{\sol}(\flow(h,A(\tau_{0})))
	\leq \fench_{\sol}(\flow(0, A(\tau_{0})))
	= \fench_{\sol}(A(\tau_{0})) <\frac{\eps}{2}
	\quad
	\text{for all $h \geq 0$}.
	\end{equation}
	Consequently:
	\begin{equation}
	\label{eq:bound2}
	\fench_{\sol}(A(\tau_{0}+h))
	< \fench_{\sol}(\flow(h,A(\tau_{0}))) + \frac{\eps}{2}
	< \frac{\eps}{2} + \frac{\eps}{2}
	= \eps,
	\end{equation}
	where the last inequality is a consequence of \eqref{eq:bound1}.
	
	Now, assume inductively that \cref{eq:bound2} holds for all $h \in [(n-1)T_{\eps}, nT_{\eps})$ for some $n\geq1$.
	Then, for all $h\in [(n-1)T_{\eps}, nT_{\eps})$, we have:
	\begin{flalign}
	\label{eq:bound3}
	&\fench_{\sol}(A(\tau_{0} + T_{\eps} + h))
	< \fench_{\sol}(\flow(T_{\eps},A(\tau_{0}+h))) + \frac{\eps}{2}
	\leq \max\braces*{\frac{\eps}{2}, \fench_{\sol}(A(\tau_{0}+h)) -\frac{\eps}{2}} + \frac{\eps}{2}
	\notag\\
	&\leq \frac{\eps}{2} + \frac{\eps}{2}
	= \eps.
	\end{flalign}
	Consequently, \cref{eq:bound2} holds for all $h \in [nT_{\eps}, (n+1)T_{\eps})$.
This completes the induction. Taking the infimum over $\sol$ then completes our proof.

\section*{Acknowledgments}
{\small
%
%
The authors are grateful to an associate editor and two anonymous referees for their suggestions and remarks, and to John Tsitsiklis for his insightful input on previous work.
The first author in particular enjoyed the conversations with John on this topic during the winter of 2019 when visiting MIT.

Zhengyuan Zhou was partially supported by IBM Goldstine fellowship. 
P.~Mertikopoulos was partially supported by the COST Action CA16228 ``European Network for Game Theory'' (GAMENET).
P.~Mertikopoulos also received financial support from
the French National Research Agency (ANR) in the framework of
the ``Investissements d'avenir'' program (ANR-15-IDEX-02),
the LabEx PERSYVAL (ANR-11-LABX-0025-01),
MIAI@Grenoble Alpes (ANR-19-P3IA-0003),
and the grants ORACLESS (ANR-16-CE33-0004) and ALIAS (ANR-19-CE48-0018-01).}

\bibliographystyle{ormsv080}
\bibliography{references}

\end{document}